\documentclass[12pt,twoside,english]{article}
\usepackage[T1]{fontenc}
\usepackage{units}
\usepackage{amsmath}
\usepackage{amssymb}
\usepackage{graphicx}

\makeatletter

\providecommand{\tabularnewline}{\\}

\newcommand{\lyxaddress}[1]{
\par {\raggedright #1
\vspace{1.4em}
\noindent\par}
}

\usepackage{a4wide}
\usepackage{graphicx}
\usepackage{fancyhdr}
\usepackage{amsmath}
\usepackage{amssymb}
\usepackage{subfigure}
\usepackage{setspace}
\usepackage{verbatim} 
\numberwithin{equation}{section}

\newcommand{\vek}[1]{\mathchoice{\displaystyle\boldsymbol#1}
{\textstyle\boldsymbol#1}{\scriptstyle\boldsymbol#1}
{\scriptscriptstyle\boldsymbol#1}}
\newcommand{\mat}[1]{\mathchoice{\displaystyle\mathbf#1}
{\textstyle\mathbf#1}{\scriptstyle\mathbf#1}
{\scriptscriptstyle\mathbf#1}}

\makeatother

\usepackage{babel}
\begin{document}

\title{Higher-order meshing of implicit geometries---part III:\\Conformal
Decomposition FEM (CDFEM)}

\author{T.P. Fries}

\maketitle

\lyxaddress{\begin{center}
Institute of Structural Analysis\\
Graz University of Technology\\
Lessingstr. 25/II, 8010 Graz, Austria\\
\texttt{www.ifb.tugraz.at}\\
\texttt{fries@tugraz.at}
\end{center}}
\begin{abstract}
A higher-order accurate finite element method is proposed which uses
automatically generated meshes based on implicit level-set data for
the description of boundaries and interfaces in two and three dimensions.
The method is an alternative for fictitious domain and extended finite
element methods. The domain of interest is immersed in a background
mesh composed by higher-order elements. The zero-level sets are identified
and meshed followed by a decomposition of the cut background elements
into conforming sub-elements. Adaptivity is a crucial ingredient of
the method to guarantee the success of the mesh generation. It ensures
the successful decomposition of cut elements and enables improved
geometry descriptions and approximations. It is confirmed that higher-order
accurate results with optimal convergence rates are achieved with
the proposed conformal decomposition finite element method (CDFEM).

\end{abstract}
\newpage{}\tableofcontents{}\newpage{}

\section{Introduction\label{sec:Introduction}}

The $p$-version of the finite element method ($p$-FEM) enables a
higher-order accurate and efficient approximation of boundary value
problems (BVPs) in engineering, natural sciences, and related fields
\cite{Solin_2003a,Bathe_1996a,Belytschko_2000b,Zienkiewicz_2000d,Szabo_2004a}.
Two crucial requirements are needed for the successful application
of the $p$-FEM: (i) The geometry must be accurately represented by
a mesh composed of higher-order elements and (ii) the solution of
the BVP should be sufficiently smooth. Both requirements are not easily
met. For (i), \emph{curved} boundaries and interfaces in the domain
of interest may render the mesh generation difficult, in particular
in three dimensions and with elements of higher orders. Even more
so when frequent mesh manipulations are desired, for instance, in
the context of moving interfaces (interface tracking) or mesh refinements
in adaptivity and convergence studies. The original geometry is often
generated based on Computer Aided Geometric Design (CAGD or CAD) and
the interplay with the analysis tool, i.e., the $p$-FEM, is not easily
established and hardly automated. Concerning (ii), the smoothness
of the involved fields in the BVP, it is noted that discontinuities
(e.g., in the material parameters) and singularities (e.g., in the
stress field of a structure) are frequently present. For the successful
application of the $p$-FEM in these cases, it is again crucial to
provide suitable meshes, i.e., those which conform to the discontinuities
and are refined at the singularities. It is thus seen that a lot of
effort is associated to generating higher-order accurate meshes as
properties of the geometry and the approximated solutions must both
be considered.

Herein, the focus is on the \emph{automatic}, higher-order accurate
generation of conforming meshes based on implicitly defined geometries.
The domain of interest is completely immersed in a background mesh.
The boundary of the domain and interfaces therein, for example, between
different materials, are defined by (several) level-set functions
\cite{Osher_2003a,Osher_2001a,Sethian_1999b}. For each level-set
function, the elements cut by the zero-level set are decomposed into
conforming, higher-order sub-elements. Therefore, the zero-level set
is first identified and meshed by interface elements (reconstruction)
and then customized mappings generate the sub-elements (decomposition).
This follows previous works of the author in \cite{Fries_2015a,Fries_2016a,Fries_2016b}
where the resulting meshes are used in the context of integration
and interpolation in implicitly defined domains. However, in elements
where the decomposition fails, e.g., due to very complex level-set
data, (isolated) recursive refinements were suggested and hanging
nodes are a natural consequence. Herein, we wish to use the generated
meshes in the context of approximating BVPs and hanging nodes shall
be avoided. The quality of the generated sub-elements becomes an important
issue in this context. Following \cite{Loehnert_2014a,Kramer_2017a},
node manipulations in the background mesh are suggested to ensure
suitable, shape-regular elements. One may also possibly use stabilizations
similar to those suggested in \cite{Burman_2010a,Burman_2012a,Burman_2014a,Hansbo_2002a}.

Adaptive refinements of the background mesh are suggested in order
to (i) refine elements where the decomposition failed, (ii) improve
the geometry description driven by the curvature of the involved level-set
functions near the zero-level sets, (iii) improve the approximation
of the BVP, for example, based on error indicators. Because ``good''
meshes must consider both, the geometry \emph{and} the involved (sought)
fields of the BVP, adaptivity is a natural ingredient for \emph{automatic}
mesh generation without any user intervention. Hence, the suggested
procedure follows the isogeometric paradigm \cite{Hughes_2005a,Marussig_2014a}
to fully integrate design and analysis, however, for implicit geometries
rather than based on NURBS as in CAGD.

The fully automatic generation of meshes based on implicit level-set
data is gaining increasing attention. We emphasize the work of \cite{Noble_2010a}
in a low-order context for moving interfaces which coined the name
CDFEM. A higher-order extension of this work in two dimensions is
found in \cite{Omerovic_2016a} without adaptive refinements and measures
to avoid ill-shaped elements. This is the first work where the CDFEM
is extended to higher-order consistently in two and three dimensions,
including adaptivity and node manipulations to ensure the regularity
of the resulting elements. The resulting method is stable and efficient.

The decomposition of elements is frequently employed in the context
of ``fictitious domain methods'' (FDMs) such as the unfitted or
cut finite element method \cite{Burman_2010a,Burman_2012a,Burman_2014a,Hansbo_2002a},\textbf{
}finite cell method \cite{Abedian_2013b,Duester_2008a,Parvizian_2007a,Schillinger_2014a,Schillinger_2014b},
Cartesian grid method \cite{Uzgoren_2009a,Ye_1999a},\textbf{ }immersed
interface method \cite{Leveque_1994a}, virtual boundary method \cite{Saiki_1996a},
embedded domain method \cite{Lohner_2007a,Neittaanmaki_1995a} etc.
The important difference between the CDFEM and FDMs is that the first
uses the shape functions of the decomposed elements in the conforming
mesh as the approximation basis whereas the second employs the shape
functions of the original background mesh and uses the sub-elements
for integration purposes only. Integration in cut elements using element
decompositions is suggested in \cite{Abedian_2013a,Moumnassi_2011a,Dreau_2010a}
using \emph{polygonal} sub-cells together with recursive refinements.
Curved sub-cells based on higher-order elements are, e.g., used in
\cite{Legay_2005a,Cheng_2009a,Fries_2015a,Fries_2016b} and typically
lead to much less integration points. The integration schemes based
on element decompositions are also frequently employed in the context
of the extended or generalized finite element methods (XFEM/GFEM),
see e.g., \cite{Belytschko_1999a,Moes_1999a,Fries_2009b} for the
XFEM and \cite{Strouboulis_2000a,Strouboulis_2000b} for the GFEM.
They consider for inner-element jumps and kinks by adding enrichment
functions based on the partition of unity concept \cite{Babuska_1994a,Babuska_1997a,Melenk_1996a}.
Again, in these methods the decomposed sub-elements are only used
as integration cells without using the implied shape functions for
the approximation of the BVP. 

It is emphasized that the proposed higher-order accurate CDFEM may
be seen as an alternative for FDMs where \emph{boundaries} are defined
implicitly \emph{and} the XFEM where \emph{interfaces} are defined
implicitly. The numerical results show typical applications of FDMs
and the XFEM in two and three dimensions and higher-order convergence
rates are achieved.

The paper is organized as follows: In Section \ref{X_Preliminaries},
the concept of background meshes and their interaction with the implicitly
defined boundaries and interfaces based on level-set functions is
introduced and the procedure to decompose cut elements into conforming
sub-elements based on \cite{Fries_2015a,Fries_2016b} is summarized.
Section \ref{X_Adaptivity} details the adaptive refinement strategy
in elements where the decomposition fails and, in addition, to improve
the geometry representation and approximation properties. From the
resulting set of elements, the generation of a finite element mesh
including the connectivity information is outlined in Section \ref{X_MeshGeneration}.
A node manipulation scheme to ensure the shape-regularity of the generated
elements is described there as well. The proposed higher-order CDFEM
is very general, however, herein it is applied in the context of solid
mechanics with the governing equations as given in Section \ref{X_Governing-Equations-of-Lin-Elast}.
Numerical results in two and three dimensions are presented in Section
\ref{X_NumericalResults} where typical applications of FDMs and XFEM
are considered with the proposed higher-order CDFEM. Finally, the
paper ends with a summary and conclusion in Section \ref{X_Conclusions}.

\section{Preliminaries\label{X_Preliminaries}}

\begin{figure}
\centering

\subfigure[$k=1$]{\includegraphics[height=4cm]{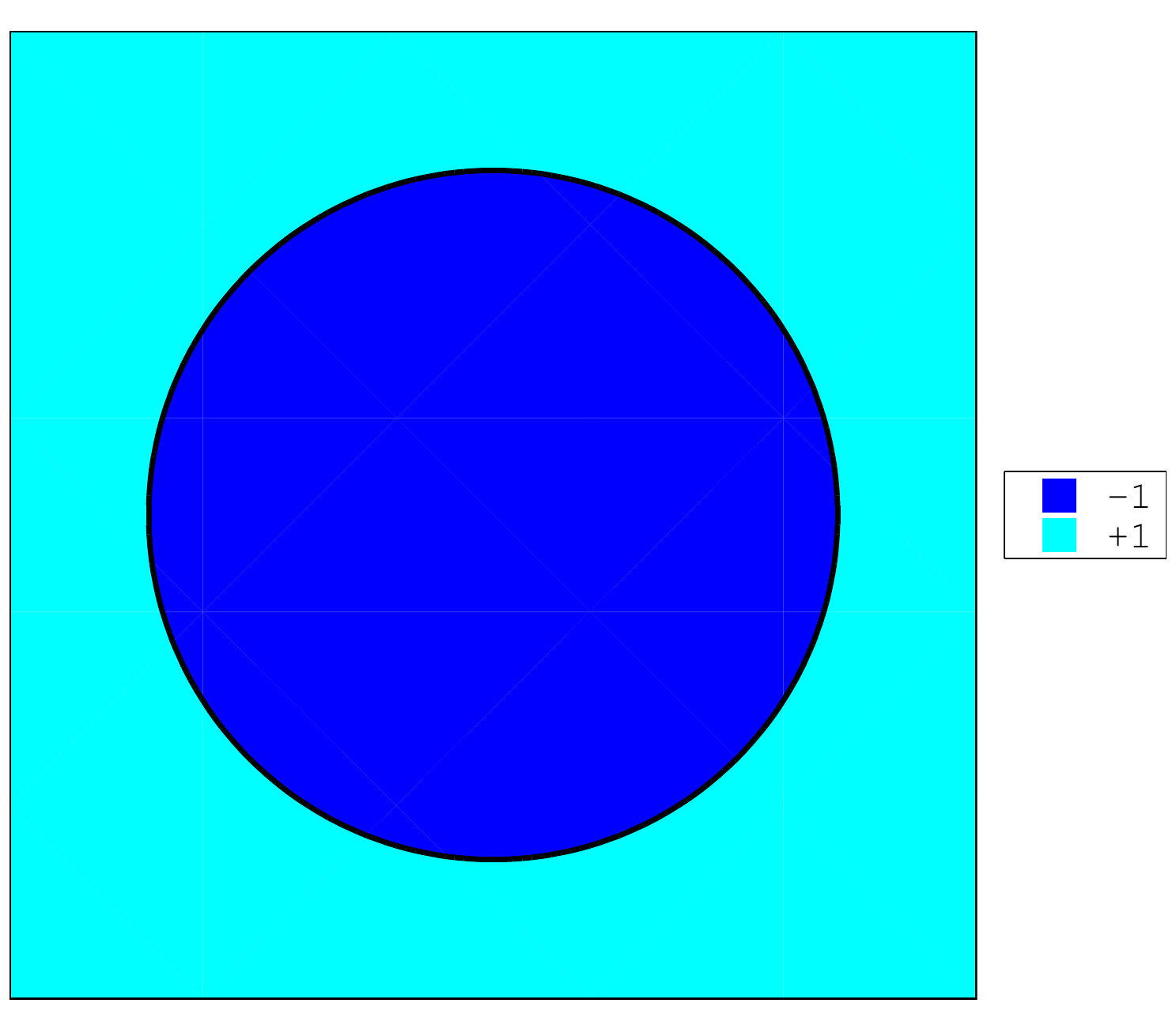}}\qquad\subfigure[$k=2$]{\includegraphics[height=4cm]{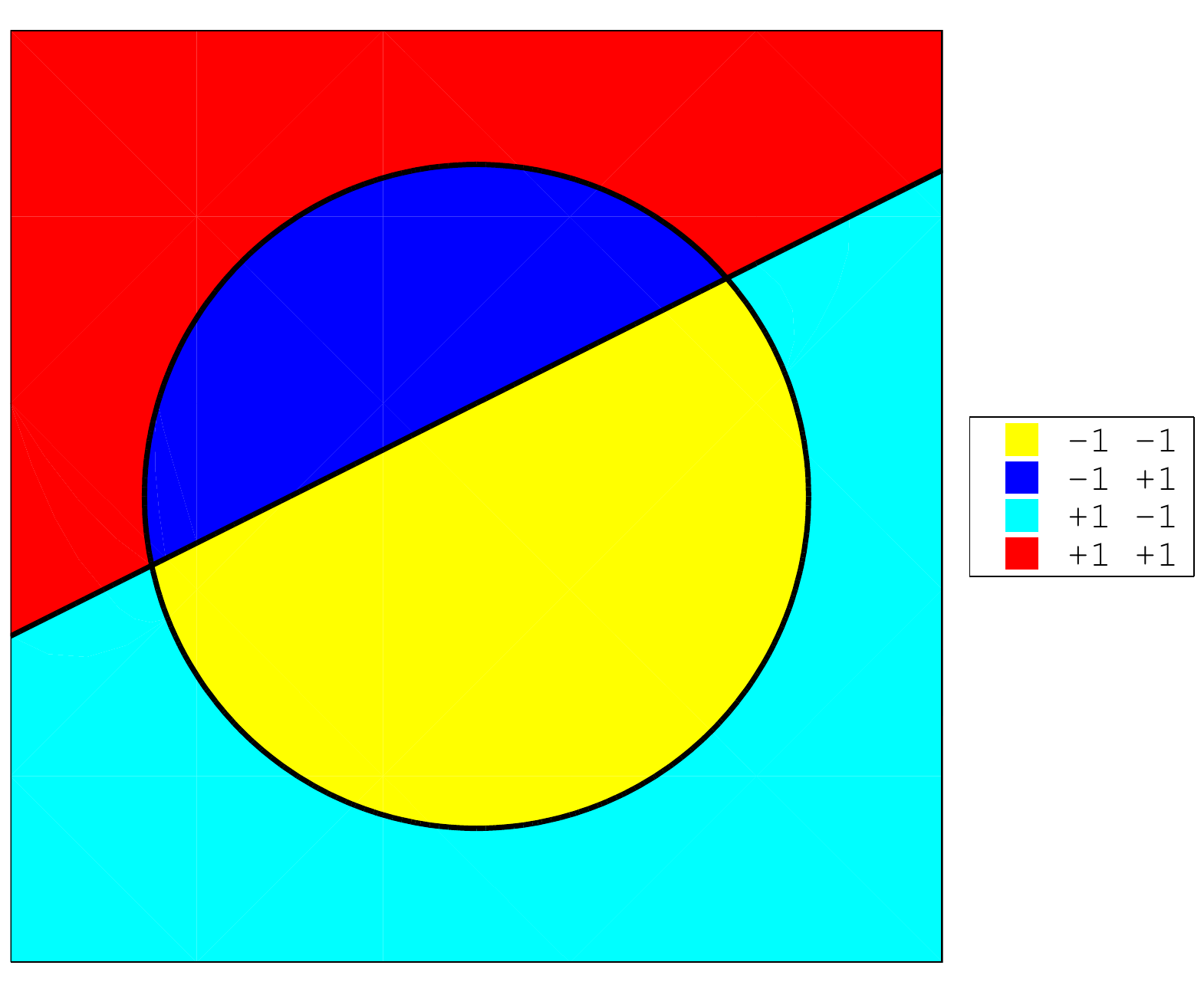}}

\subfigure[$k=3$]{\includegraphics[height=4cm]{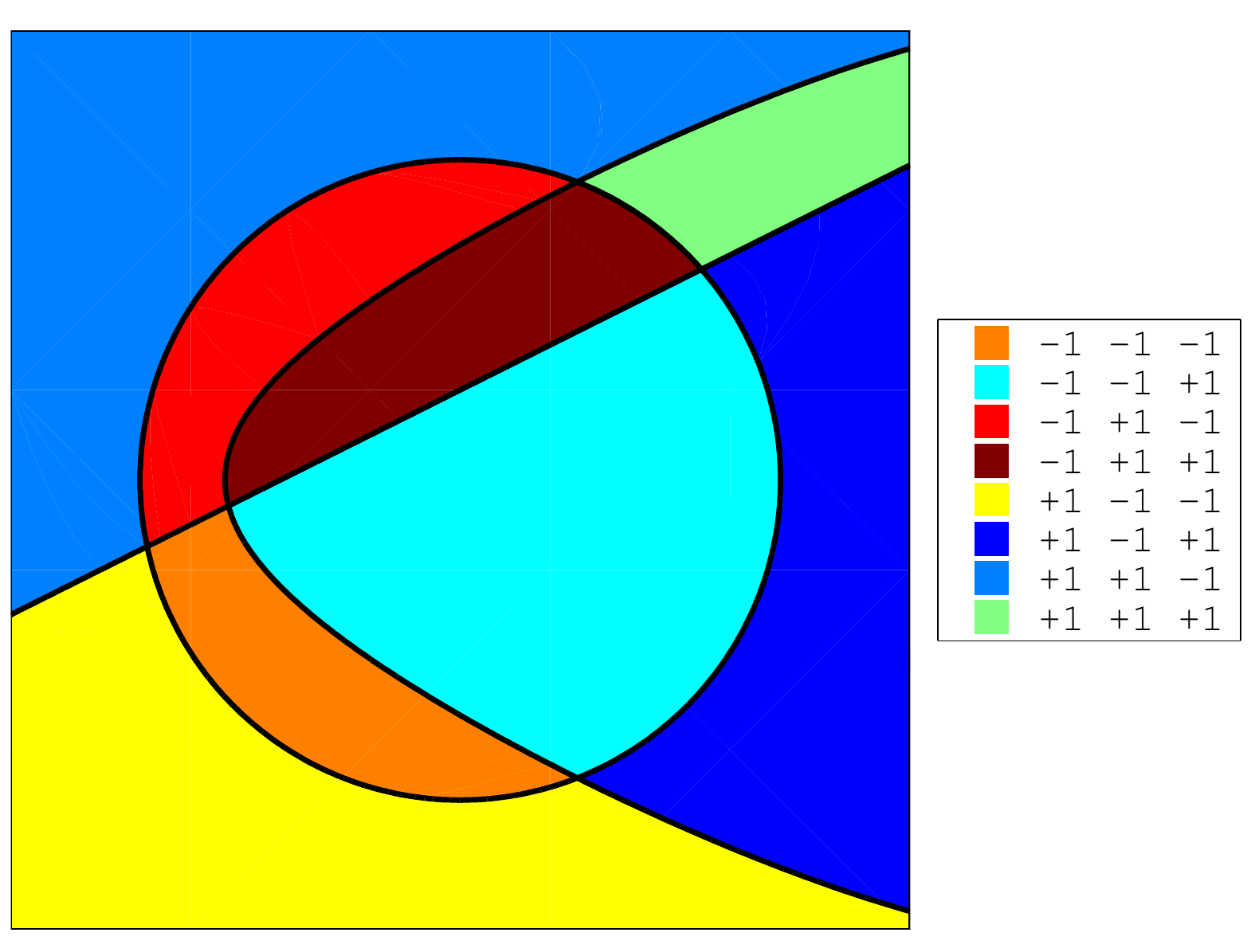}}\qquad\subfigure[conforming mesh]{\includegraphics[height=4cm]{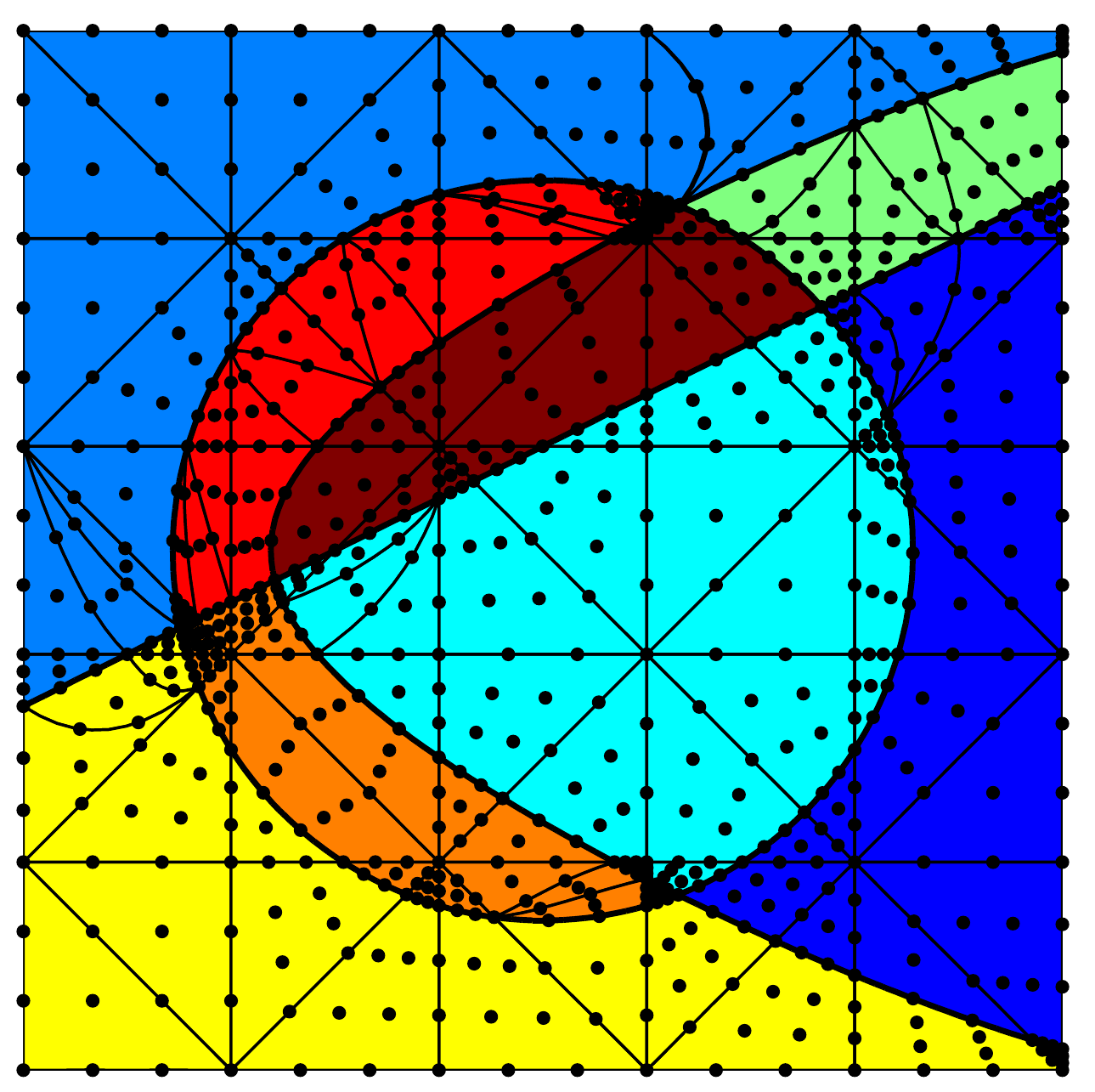}}

\caption{\label{fig:VisRegions}(a) to (c) show how the zero-level sets of
$k$ level-set functions are able to define $2^{k}$ sub-regions,
(d) shows an automatically generated, conforming mesh (composed by
cubic elements).}
\end{figure}

\begin{figure}
\centering

\subfigure[2D, level-sets and geometry]{\includegraphics[height=4cm]{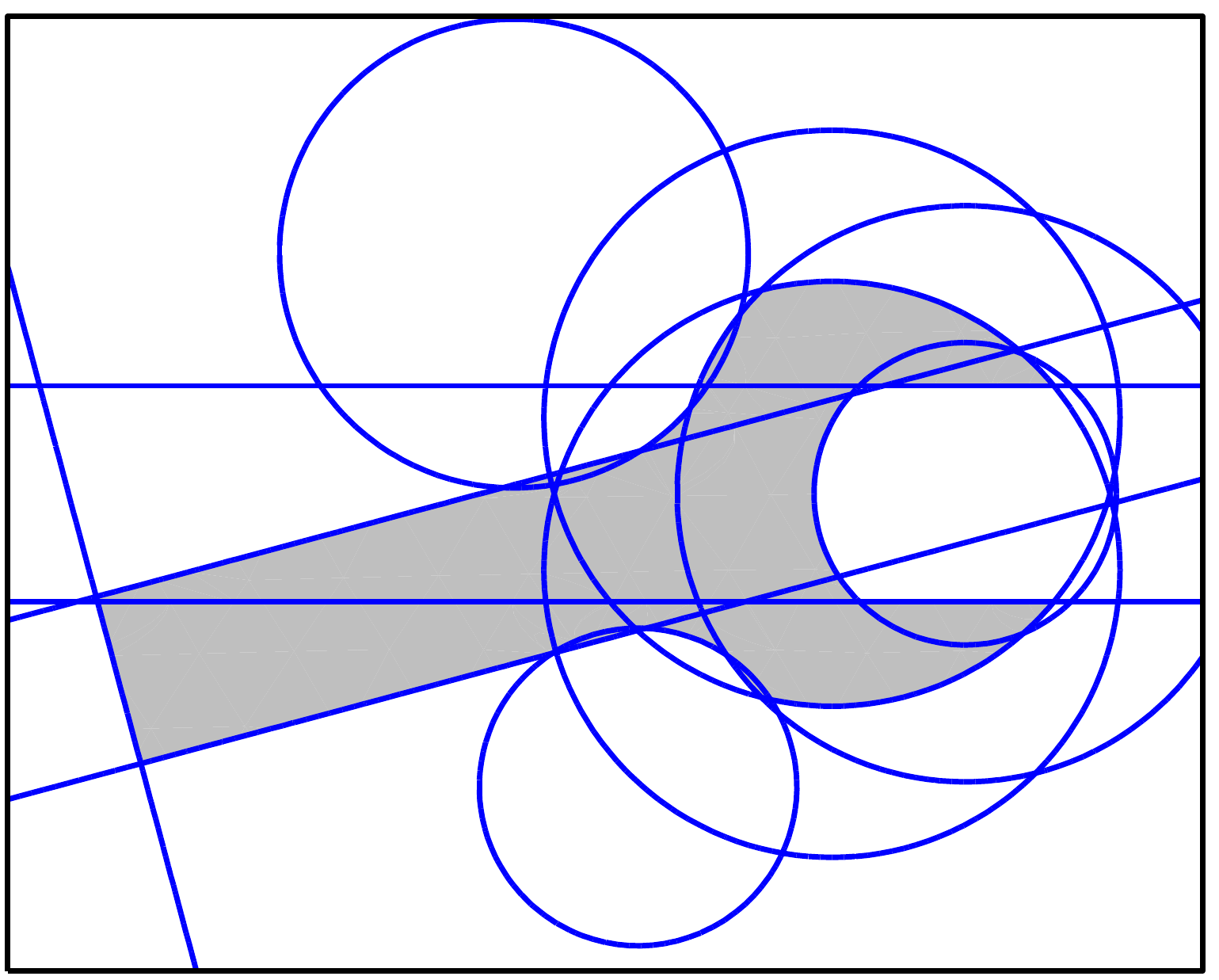}}\qquad\subfigure[3D, level-sets]{\includegraphics[height=4cm]{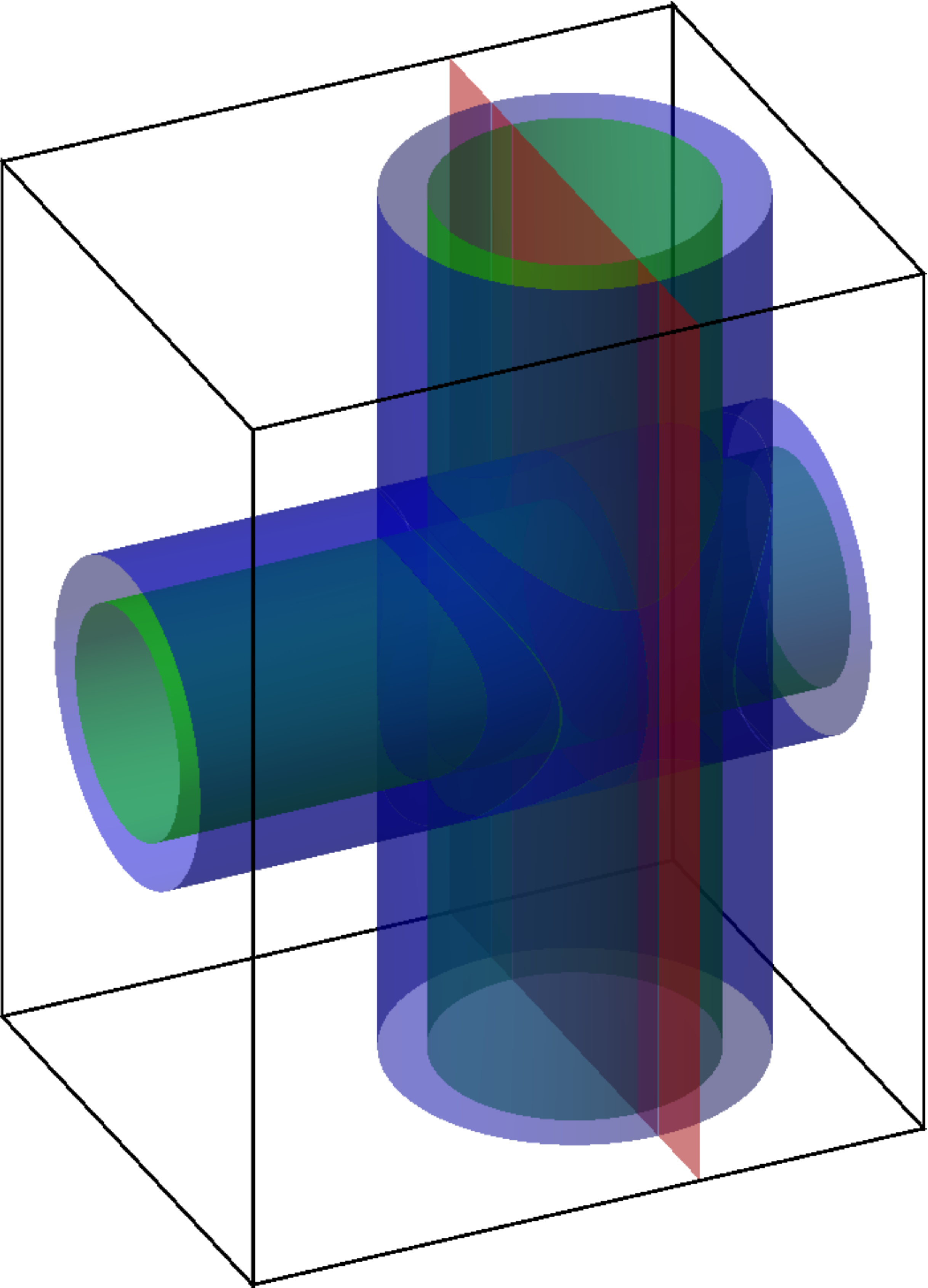}}\qquad\subfigure[3D, geometry]{\includegraphics[height=4cm]{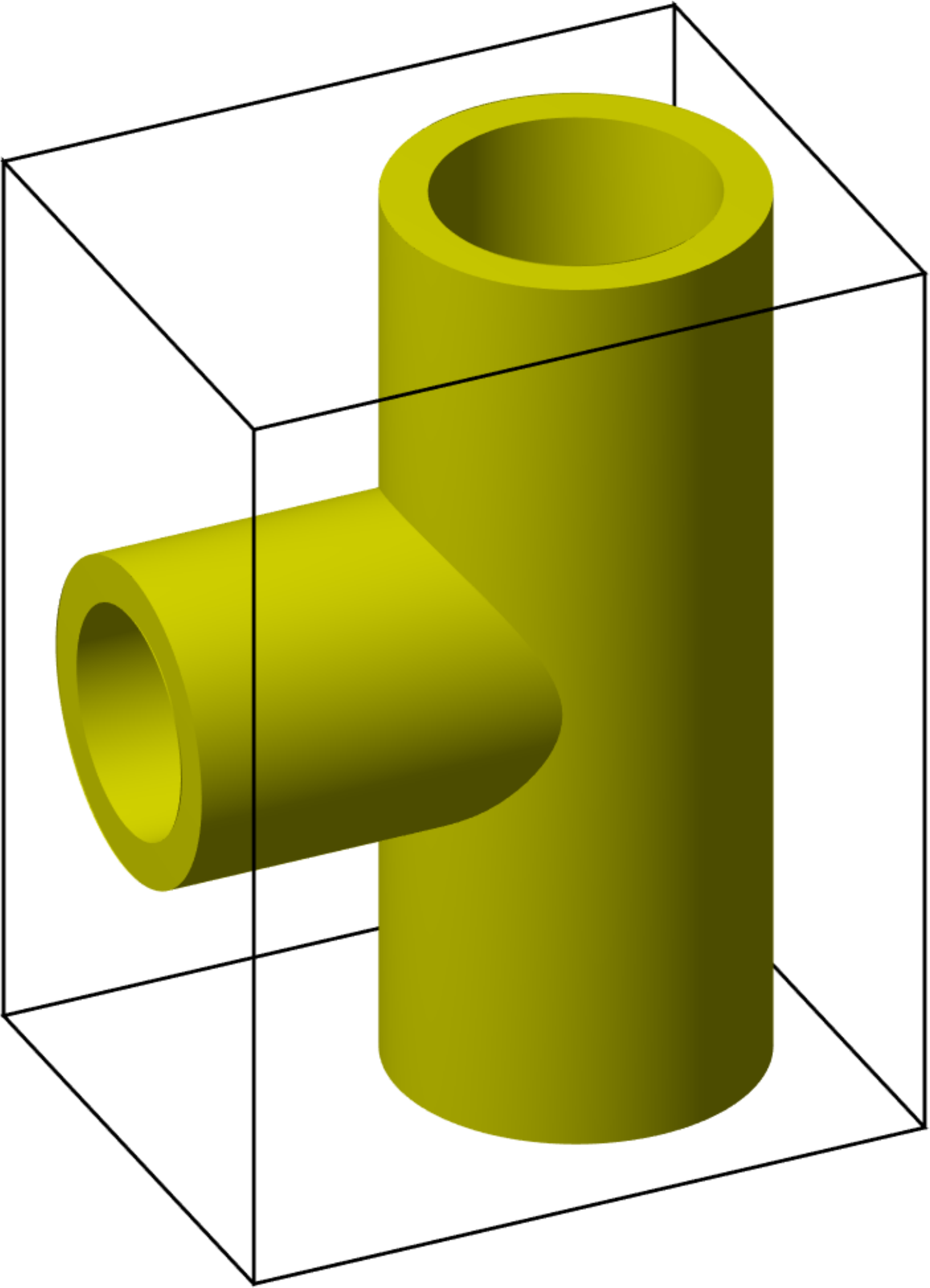}}

\caption{\label{fig:VisLevelSets}(a) Zero-level sets in 2D and the implied
geometry of a spanner, (b) and (c) show zero-level sets in 3D and
the implied geometry of a pipe junction.}
\end{figure}

Starting point is a domain of interest $\Omega$ in two or three dimensions
which is fully immersed in a background mesh composed by (possibly
unstructured) higher-order Lagrange elements. The boundary of the
domain and/or interfaces therein, for example between different materials,
are defined by the zero-contours of the level-set functions $\phi_{i}\left(\vek x\right)$.
The level-set functions are evaluated at the nodes of the background
mesh and, inbetween, interpolated by $\phi_{i}^{h}\left(\vek x\right)$
based on classical finite element shape functions. The signs of the
level-set functions define sub-regions in the background mesh and
it is easily verified that $k$ level-set functions may identify a
maximum of $2^{k}$ subregions, see Fig.~\ref{fig:VisRegions}(a)
to (c). It is important to note that several level-set functions naturally
imply corners and edges of the domain of interest which has already
been discussed in \cite{Fries_2016b}. Consequently, a purely implicit
description of complex geometries of practical interest is possible
using multiple level-set functions. See Fig.~\ref{fig:VisLevelSets}
for some examples where several zero-level sets in two and three dimensions
are shown and the implied geometries are highlighted.

\begin{figure}
\centering

\subfigure[Remeshing in 2D]{\includegraphics[width=11cm]{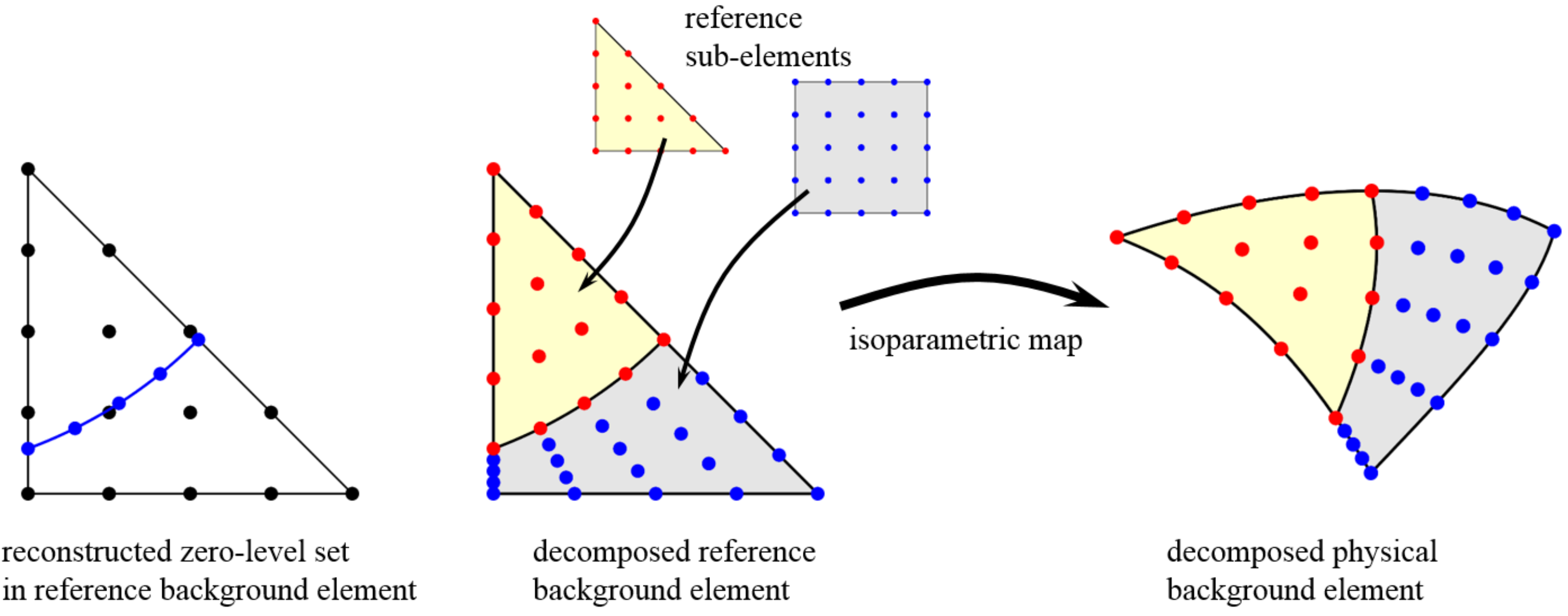}}\\\subfigure[Remeshing in 3D]{\includegraphics[width=11cm]{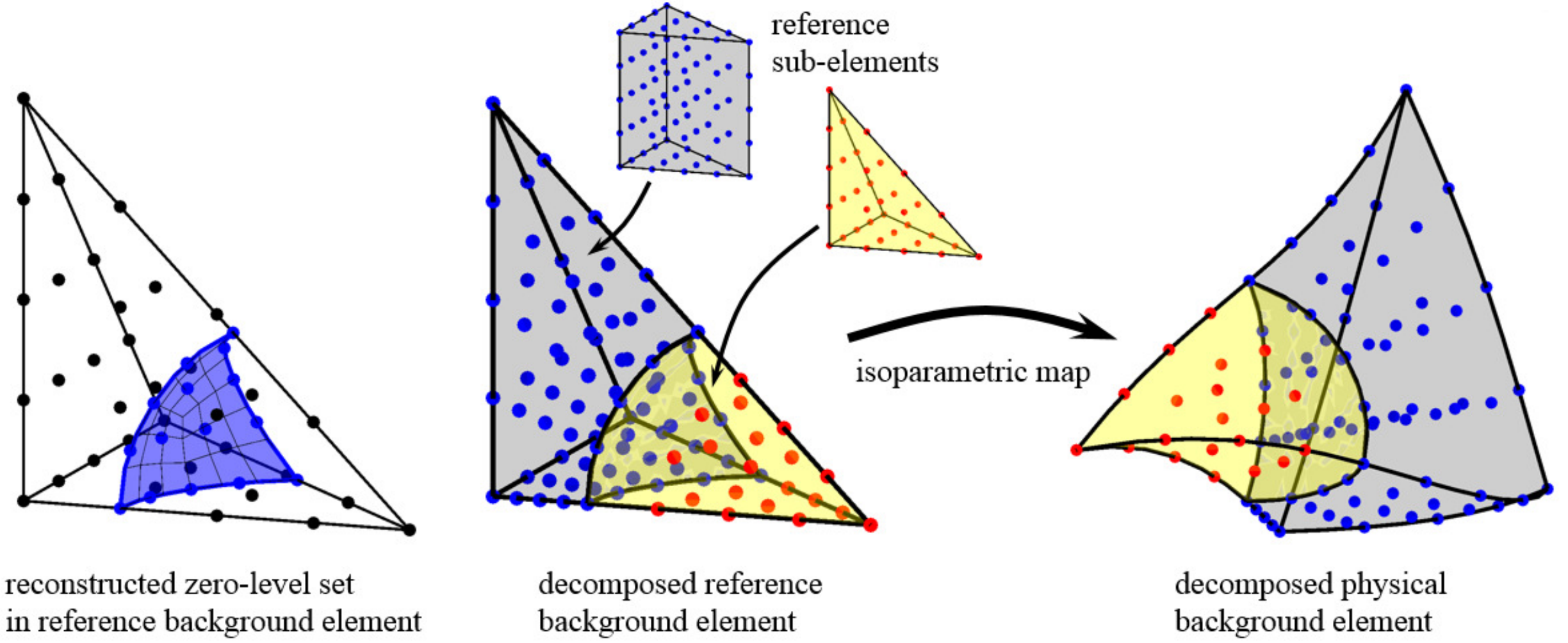}}

\caption{\label{fig:VisOverviewRemesh}The proposed automatic mesh generation
in (a) two and (b) three dimensions. The major steps are reconstruction
and decomposition in the reference element followed by the mapping
to the physical background element. Partly taken from \cite{Fries_2016b}.}
\end{figure}

The task is to automatically generate higher-order meshes which conform
to the boundaries and interfaces defined by the zero-level sets, see
e.g., Fig.~\ref{fig:VisRegions}(d). Based on the sign-combinations
of the involved level-set functions, void regions are easily identified
and/or material properties assigned. The mesh generation with respect
to all level-set functions is realized one after the other. For each
level-set function, the following steps are performed for every elements,
see Fig.~\ref{fig:VisOverviewRemesh} and references \cite{Fries_2015a,Fries_2016b}
for further details:
\begin{enumerate}
\item Detection whether the element is cut by the current level-set function
or not. This is based on a sample grid because nodal values are not
sufficient for this decision. For cut elements proceed with step 2,
otherwise with the next element.
\item Determine how the zero-level set cuts the element and classify the
topological cut situation provided that the level-set data is not
too complex, see below.
\item Reconstruction: In the reference element, identify the zero-level
set and define interface elements. Therefore, element nodes are identified
on the zero-level set along specified search paths for which a tailored
Newton-Raphson scheme is employed. The definition of such search paths
and the corresponding start values for the iteration are crucial for
the success.
\item Decomposition: Decompose the reference element based on the reconstructed
interface element wherefore customized mappings of sub-elements are
employed depending on the topological cut situation.
\item Map the decomposed sub-elements from the reference to the physical
background element.
\end{enumerate}
Note that these steps do not necessarily lead to a successul decomposition
of an element. For example, (a) the zero-level set may be too complex
and cuts the element several times so that no standard topological
cut situation is present, (b) the identified nodes on the zero-level
set are outside the element, or (c) the Jacobian of a decomposed sub-element
may be negative, hence, invalid. Therefore, it was suggested in \cite{Fries_2015a,Fries_2016b}
to use recursive refinements of the element until the reconstruction
and decomposition are successful. As a consequence, some of the resulting
sub-elements feature hanging nodes so that such meshes are ``irregular''.
This is not a problem in an integration and interpolation context
as in \cite{Fries_2015a,Fries_2016b}, however, in the context of
approximating BVPs, it is highly benefitial to have \emph{regular}
meshes without hanging nodes. Hence, herein we wish to avoid recursive
refinements in the sense of \cite{Fries_2015a,Fries_2016b} and suggest
to use \emph{adaptive mesh refinements} instead, which enables the
generation of regular, conforming\emph{ }meshes.

\section{Adaptive mesh refinements\label{X_Adaptivity}}

\begin{figure}
\centering

\subfigure[background mesh]{\includegraphics[height=4cm]{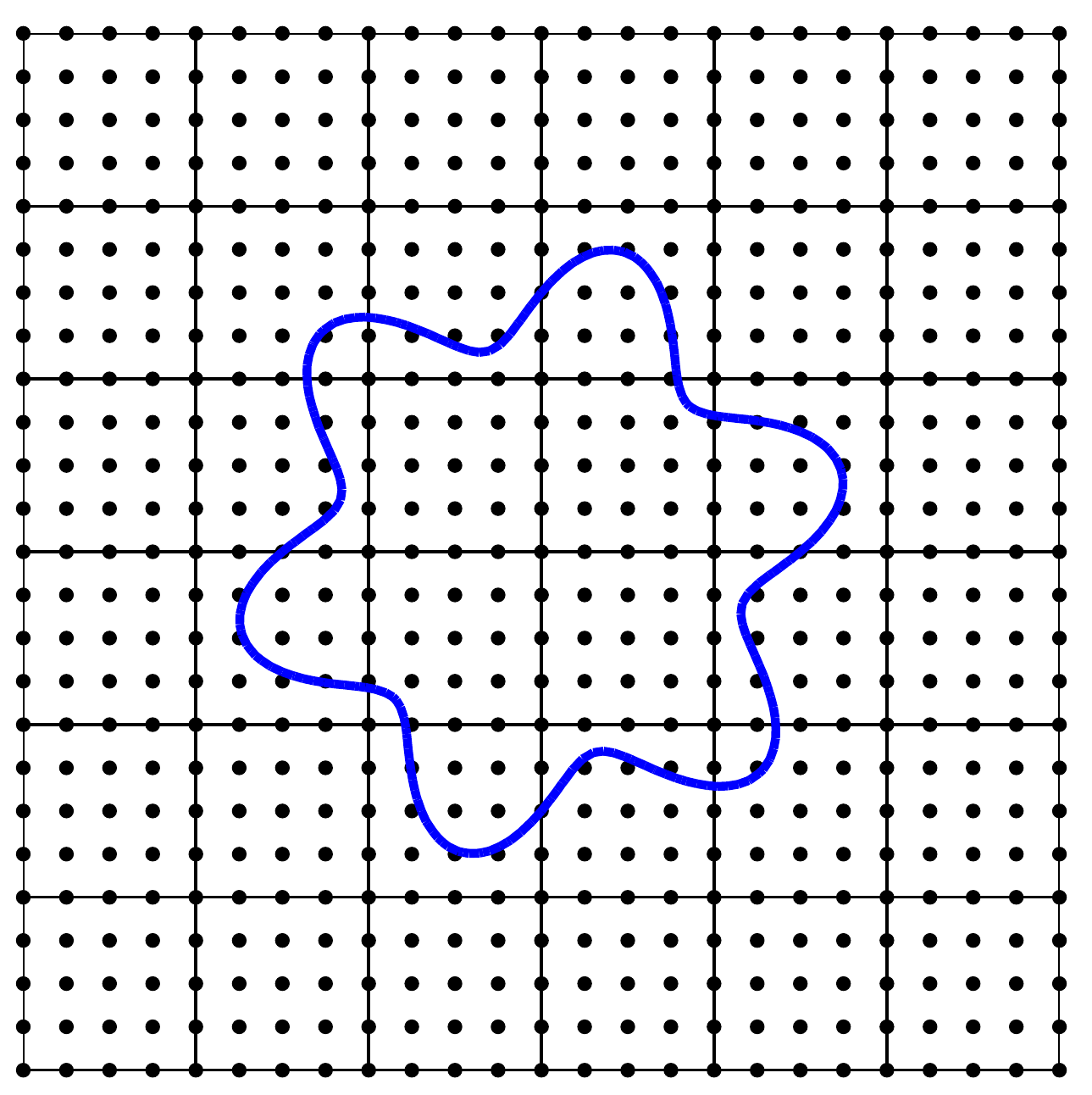}}\qquad\subfigure[decomposition]{\includegraphics[height=4cm]{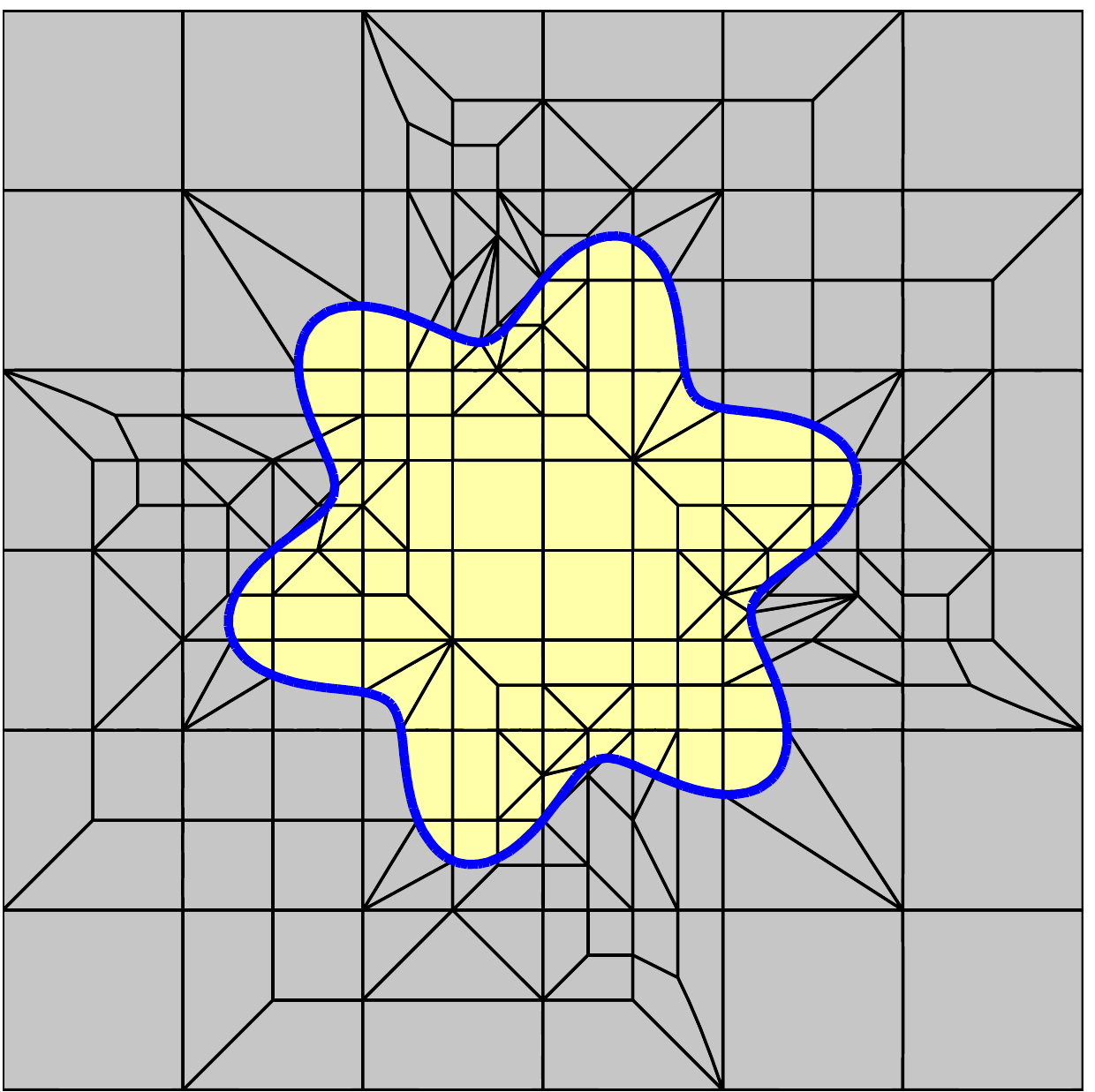}}

\subfigure[improve geometry]{\includegraphics[height=4cm]{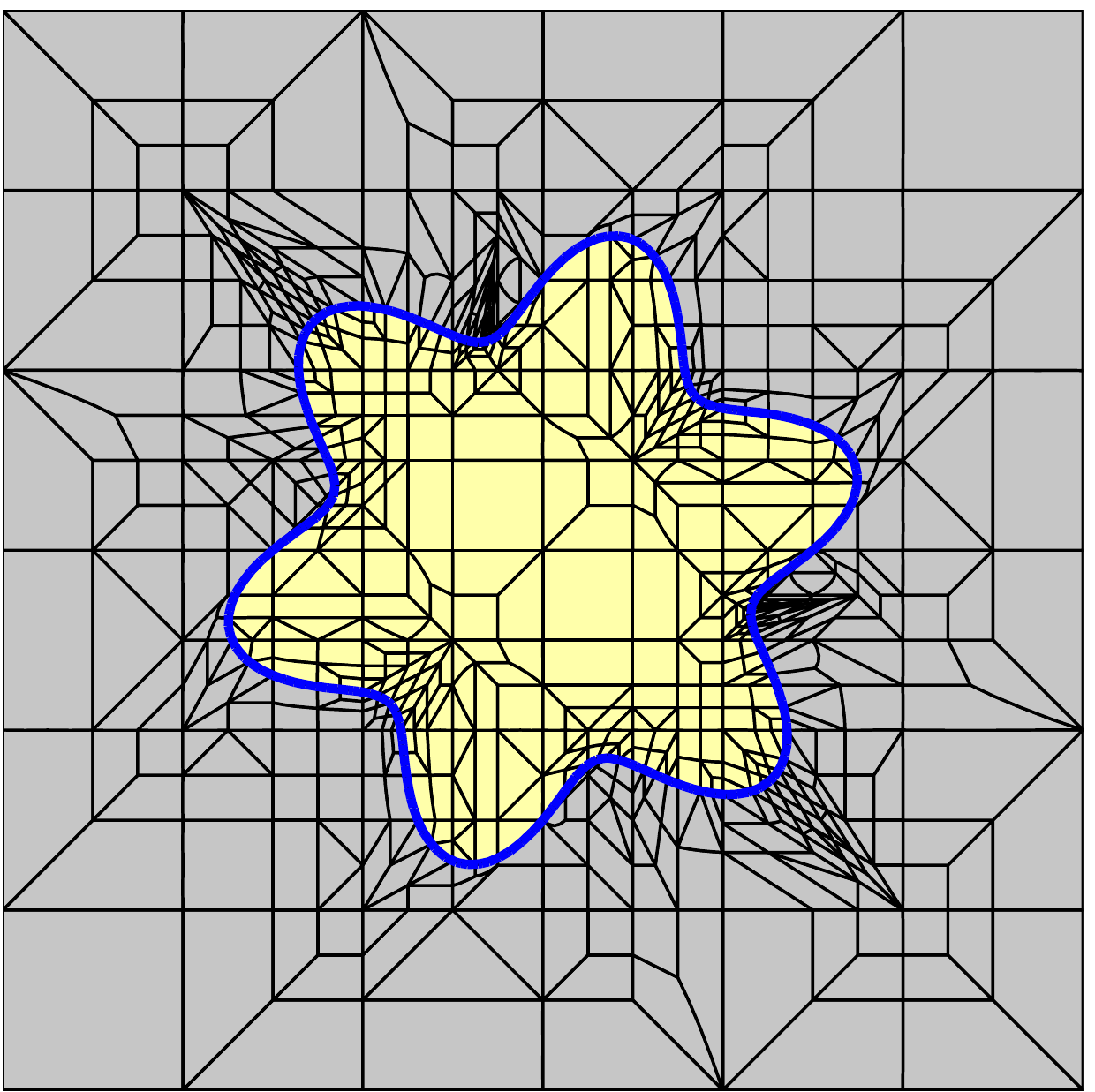}}\quad\subfigure[improve approx.]{\includegraphics[height=4cm]{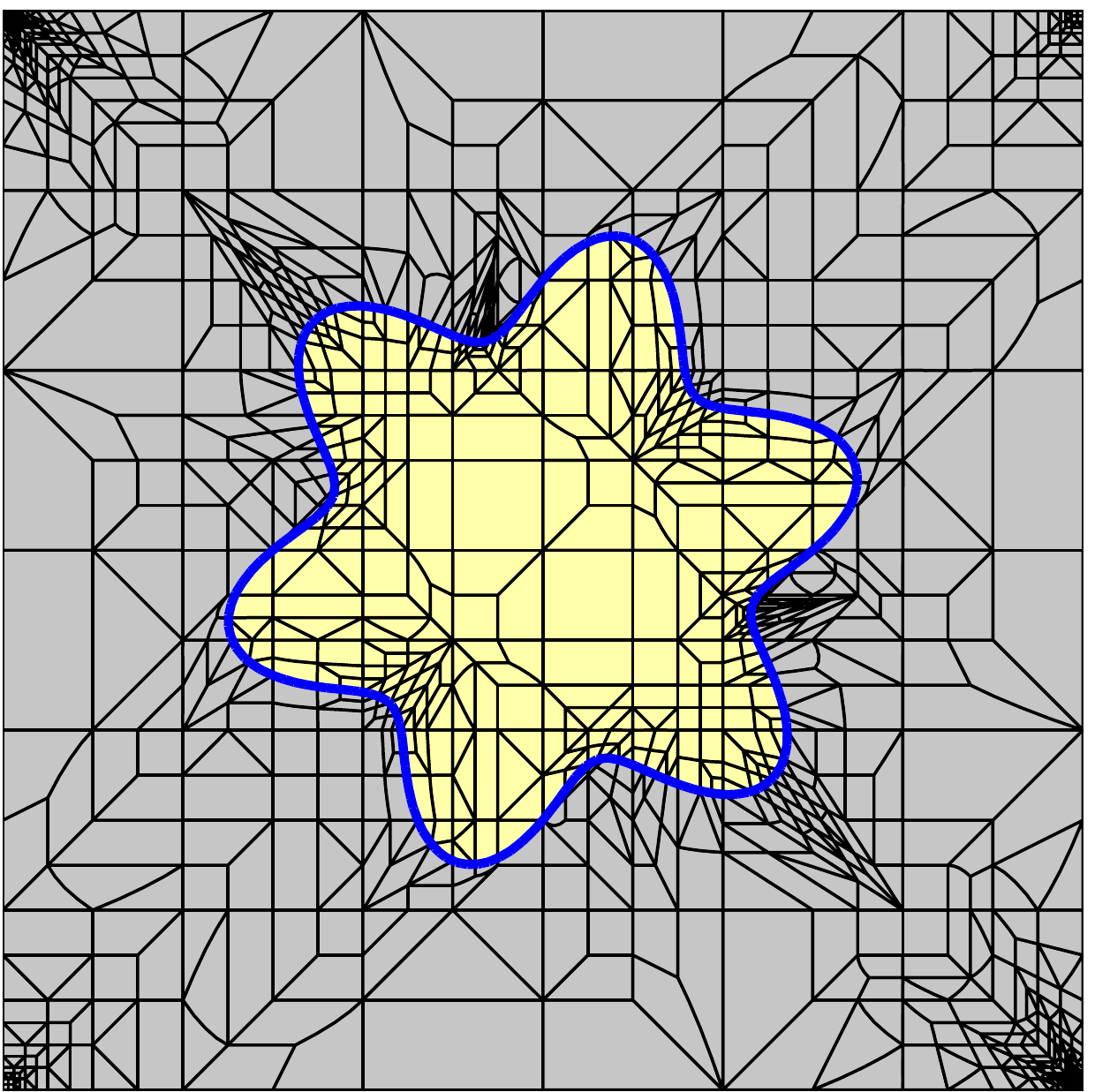}}\quad\subfigure[change order]{\includegraphics[height=4cm]{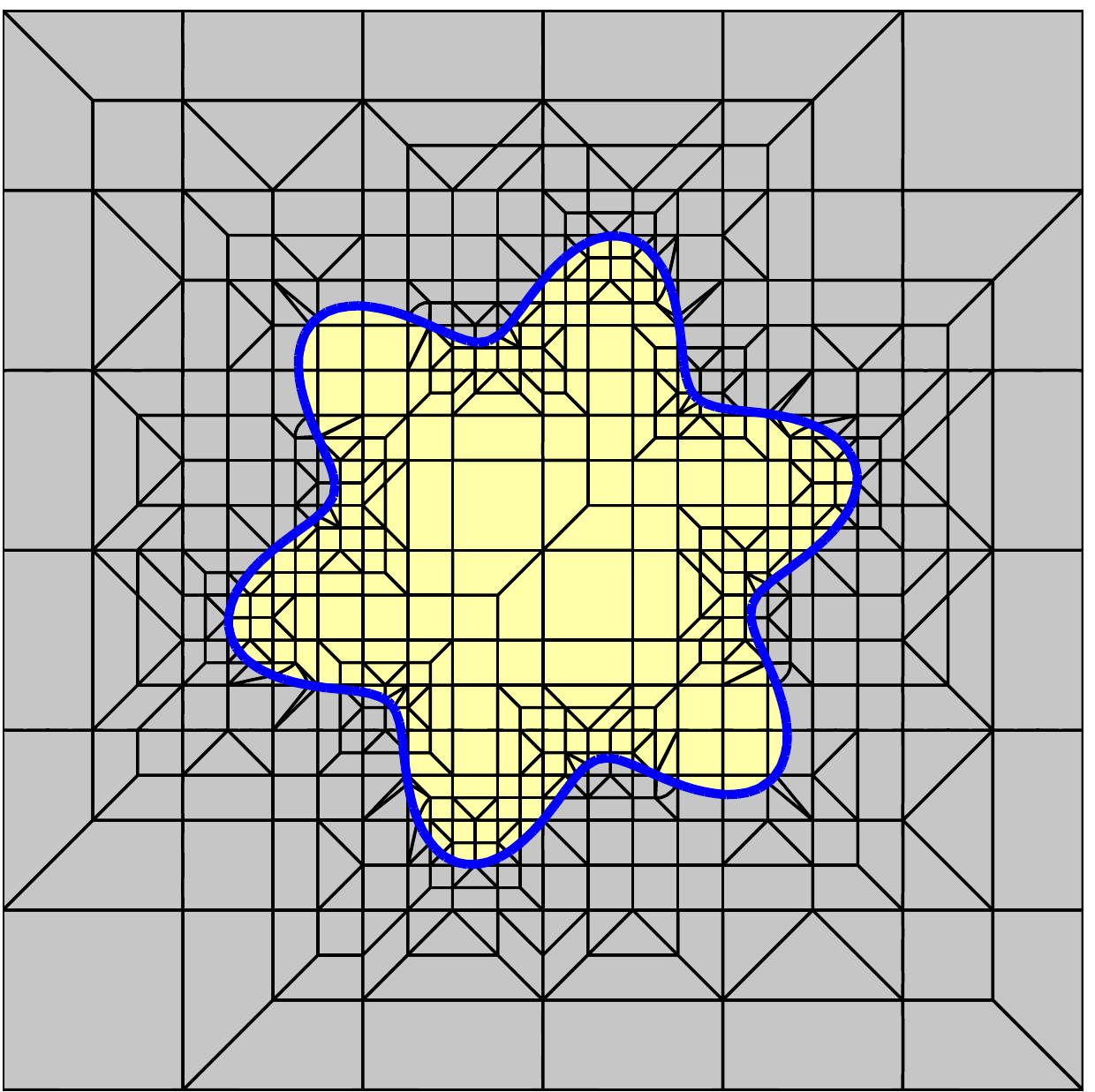}}

\caption{\label{fig:VisAdaptivity}(a) A background mesh and a zero-level set,
(b) shows adaptive refinements to cure elements where the decomposition
fails, (c) further adaptive refinements to improve the geometry description,
(d) further adaptive refinements to improve the approximation for
the example of a BVP whose solution features singularities in the
corners. In (e), the sequence of the refinement criteria is changed
to the recommended order, leading to a superior mesh.}
\end{figure}

The decomposition of all cut elements of a background mesh yields
a \emph{regular} mesh without hanging nodes provided that (i) the
background mesh itself is regular and (ii) the decompositions may
be successfully realized in \emph{all} cut elements. However, the
second criterion cannot, in general, be guaranteed and there is typically
a small number of background elements where the decomposition fails
even for smooth level-set data. Those elements have to be refined
and, in order to still meet criterion (i), also some neighboring elements
of the background mesh are affected. This is shown in Fig.~\ref{fig:VisAdaptivity}
for an extreme case where a very coarse background mesh is used for
a rather complex zero-level set. An (unusual) large number of elements
has to be refined because the decomposition fails for reasons mentioned
above, see Fig.~\ref{fig:VisAdaptivity}(b). Sometimes, even further
refinement steps are required for an original background element until,
finally, the decomposition of all refined sub-elements is valid. It
is also clearly seen, that the refined background mesh is regular,
for which neighboring elements may have to be refined as well (although
they are not even cut by the zero-level set). 

We find that at least for two further reasons, adaptive refinements
of the background mesh may be useful: to better capture the geometry
of the boundaries and interfaces and, in the context of approximating
BVPs, to improve the approximation. To improve the geometry description
it is useful to employ a curvature criterion in cut elements: The
curvature $\varkappa$ of the level-set function is evaluated in the
element and compared to the ``element length'' $h$. We use the
mean curvature, defined in two and three dimensions as
\begin{eqnarray*}
\varkappa_{\mathrm{2D}} & = & \dfrac{\phi_{,xx}\cdot\phi_{,y}^{2}-2\phi_{,x}\phi_{,y}\phi_{,xy}+\phi_{,yy}\cdot\phi_{,x}^{2}}{\left(\phi_{,x}^{2}+\phi_{,y}^{2}\right)^{\nicefrac{3}{2}}},\\
\varkappa_{\mathrm{3D}} & = & \nicefrac{1}{2}\cdot\dfrac{\begin{array}{c}
\left(\phi_{,yy}+\phi_{,zz}\right)\cdot\phi_{,x}^{2}+\left(\phi_{,xx}+\phi_{,zz}\right)\cdot\phi_{,y}^{2}+\left(\phi_{,xx}+\phi_{,yy}\right)\cdot\phi_{,z}^{2}\\
-2\phi_{,x}\phi_{,y}\phi_{,xy}-2\phi_{,x}\phi_{,z}\phi_{,xz}-2\phi_{,y}\phi_{,z}\phi_{,yz}
\end{array}}{\left(\phi_{,x}^{2}+\phi_{,y}^{2}+\phi_{,z}^{2}\right)^{\nicefrac{3}{2}}}.
\end{eqnarray*}

The element length may be the maximum Euclidean distance between every
pair of corner nodes of a physical element. One may then use the criterion
\begin{equation}
\dfrac{1}{\varkappa}\leq q\cdot h\qquad\text{with}\qquad0<q\in\mathbb{R}<2\label{eq:CurvatureCriterion}
\end{equation}
to mark elements for refinement. The value $q$ tunes the criterion
between $0$ for no curvature-driven refinement to $2$ (or larger)
for typically very curvature-sensitive refinements in cut elements.
An example for $q=0.6$ is shown in Fig.~\ref{fig:VisAdaptivity}(c).

Finally, adaptivity may also be useful when the automatically generated,
conforming meshes resulting from the previous steps are employed in
the context of approximating BVPs. This is the classical application
of adaptivity in the FEM, see e.g., \cite{Ainsworth_1997b,Demkowicz_2007a,Deuflhard_2012a,Solin_2003a}.
Although hanging nodes in the refined meshes are avoided herein, they
do not, in general, pose insurmountable problems in classical $hp$-FEM,
see e.g., \cite{Solin_2008a,Morton_1995a}. The refinement criteria
may be based on error indicators or heuristic criteria such as near
reentrant corners where singularities are expected. Algorithmically,
when adaptive refinements are already implemented for the reasons
mentioned above, it is only little effort to also enable adaptivity
to improve approximations. To continue the schematic example from
above, see Fig.~\ref{fig:VisAdaptivity}(d) where a refinement has
been realized in the corners of the square domain, for example, because
singularities are expected there for point supports in a solid mechanics
context.

It is noted that adaptivity with respect to elements where the decomposition
fails and those where the curvature criterion from Eq.~(\ref{eq:CurvatureCriterion})
fails, may easily be combined in one element loop. It is recommended
to first check the curvature criterion and only try the decomposition
when the curvature of the level-set function is sufficiently small
with respect to the element size. It is then typical that the decomposition
fails in less than $1\%$ of the elements. Note that the order of
the adaptive refinements is not commutative. Fig.~\ref{fig:VisAdaptivity}(e)
shows the resulting mesh when the curvature criterion is enforced
first followed by the decomposition; only very few elements have to
be further refined then. It is obvious that this leads to a superior
mesh than in Fig.~\ref{fig:VisAdaptivity}(c) where refinements have
first been made to enable the decomposition and, thereafter, to enforce
the curvature criterion.

\section{Mesh generation\label{X_MeshGeneration}}

The procedure from above yields a set of higher-order elements conforming
to the inner-element boundaries and/or interfaces, yet without information
on the inter-element connectivity. It is important to ensure that
across element boundaries, the generated element nodes are exactly
at the same positions. For example, in three dimensions it may be
useful to first generate element nodes on the element faces (achieved
in 2D reference elements and mapped to the physical face element),
generating a wireframe model of the zero-level sets. Next, the inner
element nodes are generated based on 3D reference elements mapped
to the physical background elements. It is then simple to generate
the connectivity of the nodes needed for the complete definition of
a finite element mesh. It is noted that the resulting meshes are \emph{mixed},
e.g., in two dimensions, they are composed by triangular \emph{and}
quadrilateral Lagrange elements. Of course this could be avoided by
converting elements to one type only.

With each element, we store the information of the signs of all level-set
functions needed for the definition of boundaries and interfaces.
Based on this information, one may easily identify elements that are
outside the domain of interest or associate material properties in
individual sub-regions of the domain.

Valid meshes for the approximation of BVPs must feature \emph{shape
regular} elements. However, this cannot generally be guaranteed for
arbitary level-set data on a given background mesh. Therefore, we
suggest to move nodes of the background mesh to ensure the shape regularity.
The aim is to bound the areas/volumes of the elements from below.
This is ensured by moving corners nodes of the elements away from
the zero-level sets. Only the nodes in a close band around the zero
isosurface are moved. The procedure was described by the author in
detail in \cite{Fries_2017a} and is only outlined here. It is applied
before the decomposition is started (however, after a potential adaptive
refinement of the background mesh due to the curvature criterion from
above). It is also noted that \cite{Loehnert_2014a,Kramer_2017a}
suggest node manipulations in related contexts, however, the concrete
algorithm from \cite{Fries_2017a} and herein is quite different.

\begin{figure}
\centering

\subfigure[example 1]{\includegraphics[width=5cm]{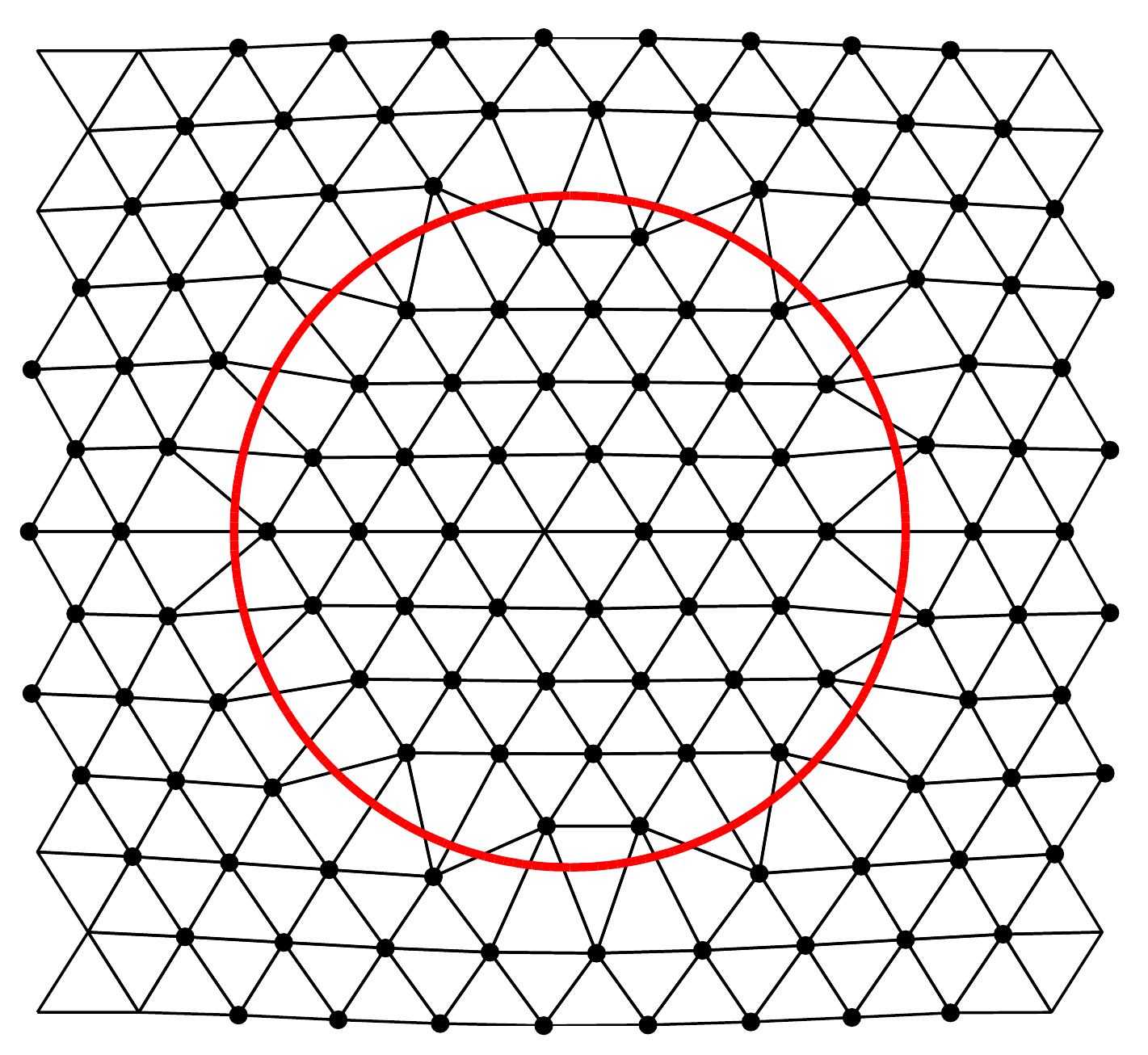}}\qquad\subfigure[example 2]{\includegraphics[width=5cm]{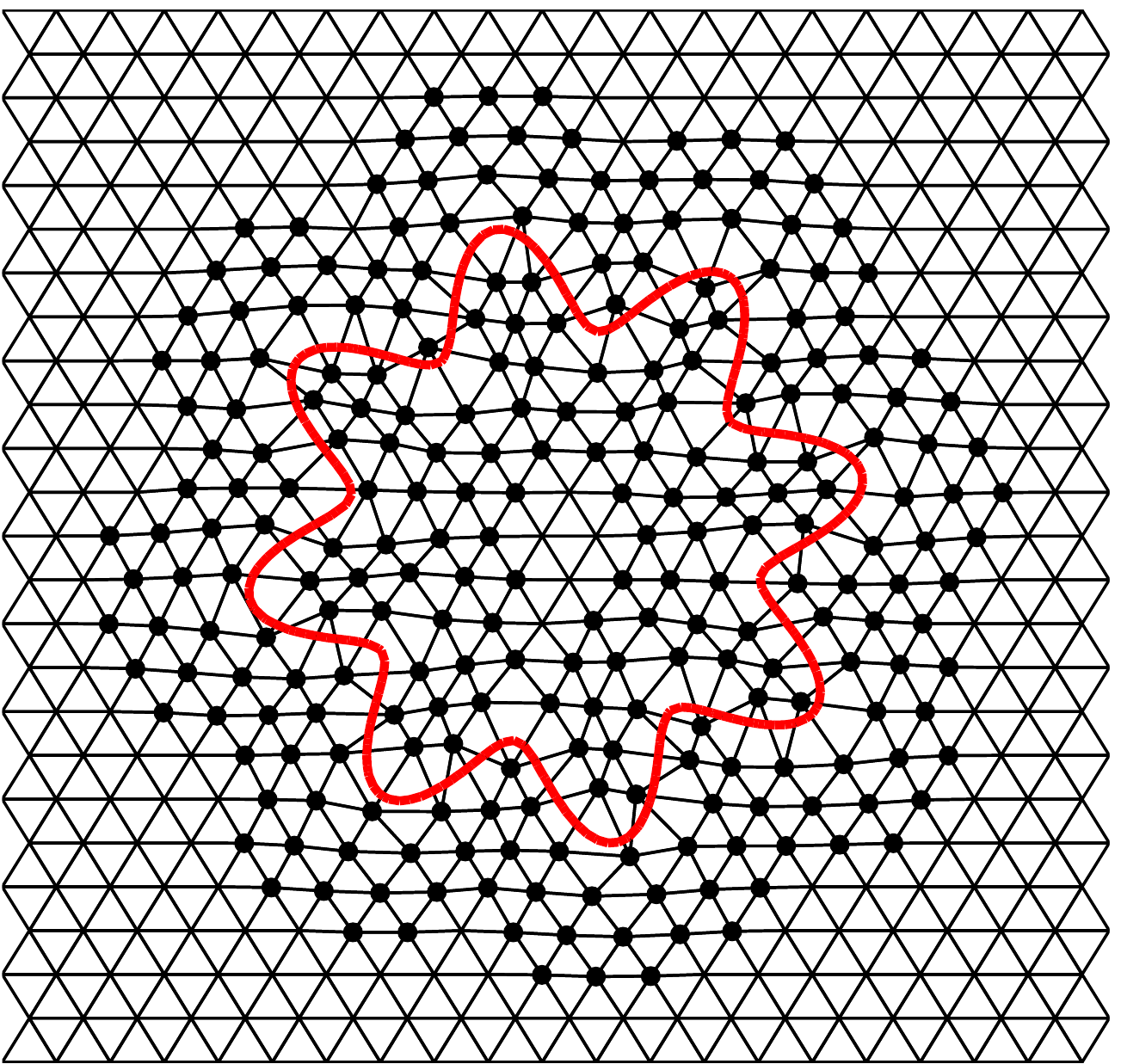}}

\caption{\label{fig:PlotMeshManip2d}Examples for node movements in 2D based
on different level-set functions and background meshes, the red lines
are the zero-level sets. Figures taken from \cite{Fries_2017a}.}
\end{figure}

The procedure of the node manipulations is split into the following
steps which are realized successively for all level-set functions
involved: (i) The distance of the nodes to the zero-level set is approximated
using a Newton-Raphson-type approach (this step is not needed when
the level-set functions feature signed-distance property). (ii) The
direction to the corresponding node on the zero-level set is measured.
This is not necessarily the exact \emph{shortest} distance, however,
it will be a good approximation for nodes which are close to the zero-level
set. (iii) If the distance is below a given threshold depending on
the element length $h$, the node is moved away from the zero-level
set. The moving distance depends on the distance itself and is ramped
linearly within the narrow band around the zero-level set. The procedure
(i) to (iii) is repeated resulting in a fix-point iteration. For further
details, see \cite{Fries_2017a}. Examples of manipulated background
meshes are seen in Fig.~\ref{fig:PlotMeshManip2d}. We note that
the node manipulations must be sufficiently small to maintain the
validity of the background mesh which is not a problem in general.
In particular the concept of ``universal meshes'' \cite{Rangarajan_2012a}
allows for a large range of manipulations of individual nodes without
leading to invalid elements (with negative Jacobians).

We summarize the differences between the proposed strategy to automatically
generate valid finite element meshes for approximations compared to
those decompositions needed only for integration purposes, e.g., in
the context of the XFEM and FDMs as proposed in \cite{Fries_2015a,Fries_2016b}.
Here, regular adaptive refinements are suggested in contrast to recursive
refinements that generate hanging nodes. A node manipulation scheme
is used to ensure the shape regularity which was not a concern in
the integration context. The generated elements must be $C_{0}$-continuous
whereas they may be discontinuous for the numerical integration. It
is noted that the related method suggested in \cite{Omerovic_2016a}
is (i) only two-dimensional, (ii) avoids refinements for the price
of more topologically different decompositions of cut elements into
sub-elements, and (iii) does not describe the issue of node manipulations.
Core features of the proposed method herein are adaptivity, node manipulations,
a simple decomposition into a minimal number of sub-elements, and
a consistent, similar treatment in two and three dimensions.

\section{Governing Equations of Linear Elasticity \label{X_Governing-Equations-of-Lin-Elast}}

The proposed higher-order CDFEM is applicable for the approximation
of general BVPs when inner-element boundaries and interfaces are present
and manufactured meshes are to be avoided. Herein, as a representative
field of application, we focus on structural mechanics and linear
elasticity. The corresponding governing equations are presented within
this section.

\begin{figure}[htbp]
\begin{centering}
\includegraphics[height=4cm]{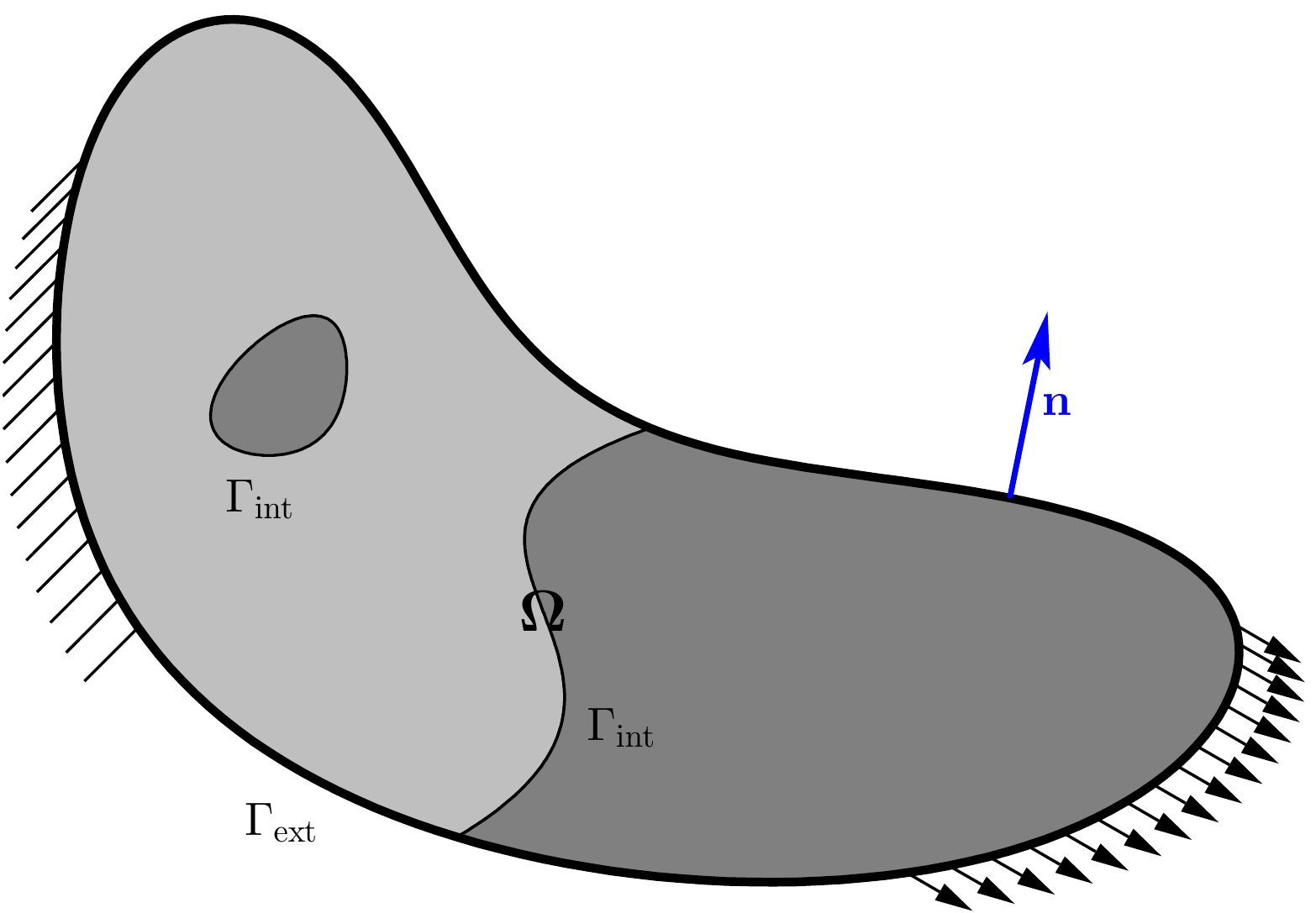}
\par\end{centering}

\caption{\label{MechDomain}Structural domain with external interfaces (boundaries)
$\Gamma_{\mathrm{ext}}$ and internal interfaces $\Gamma_{\mathrm{int}}$
implied by zero-level sets.}
\end{figure}

A domain of interest $\Omega$ in two or three dimensions is considered
which is completely immersed in a background domain $\Omega_{\mathrm{BG}}$.
The boundary of $\Omega$ may be called \emph{external} interface
$\Gamma_{\mathrm{ext}}$ and interfaces within $\Omega$, e.g., between
different materials, are called \emph{internal}, $\Gamma_{\mathrm{int}}$.
Displacements are continuous accross $\Gamma_{\mathrm{int}}$, however,
stresses and strains are discontinuous there. External and internal
interfaces are implied by zero-level sets as described above. See
Fig.~\ref{MechDomain} for a sketch of the situation in two dimensions.

The boundary $\Gamma_{\mathrm{ext}}$ is decomposed into the complementary
sets $\Gamma_{\vek u}$ and $\Gamma_{\vek t}$. Displacements $\hat{\vek u}$
are prescribed along the Dirichlet boundary $\Gamma_{\vek u}$, and
tractions $\hat{\vek t}$ along the Neumann boundary $\Gamma_{\vek t}$.
The strong form for an elastic solid undergoing small displacements
and strains under static conditions, is \cite{Belytschko_2000b,Zienkiewicz_2000d}
\begin{equation}
\nabla\cdot\boldsymbol{\sigma}=\vek f,\quad\mathrm{on}\,\,\Omega\subseteq\mathbb{R}^{2},\label{eq:LinElastStrongForm}
\end{equation}
where $\vek f$ describe volume forces, and $\boldsymbol{\sigma}$
is the following stress tensor
\begin{equation}
\boldsymbol{\sigma}=\mat C:\boldsymbol{\varepsilon}=\lambda\left(\textrm{tr}\,\boldsymbol{\varepsilon}\right)\mat I+2\mu\boldsymbol{\varepsilon},\label{eq:CauchyStress}
\end{equation}
with $\lambda$ and $\mu$ being the Lam\'e constants which are easily
related to Young's modulus $E$ and Poisson's ratio $\nu$. In two
dimensions, we shall always consider \emph{plane strain} herein. Then,
\[
\mu=\dfrac{E}{2\left(1+\nu\right)},\qquad\lambda=\dfrac{E\nu}{\left(1+\nu\right)\cdot\left(1-2\nu\right)}
\]
holds for the two and three-dimensional case. The linearized strain
tensor $\boldsymbol{\varepsilon}$ is
\begin{equation}
\boldsymbol{\varepsilon}=\frac{1}{2}\left(\nabla\vek u+\left(\nabla\vek u\right)^{\mathrm{T}}\right).
\end{equation}

For the approximation of the displacements $\vek u$, the following
test and trial function spaces $\mathcal{S}_{\vek u}^{h}$ and $\mathcal{V}_{\vek u}^{h}$
are introduced as
\begin{eqnarray}
\mathcal{S}_{\vek u}^{h} & = & \left\{ \left.\vek u^{h}\right|\vek u^{h}\in\left(\mathcal{H}^{1h}\right)^{d},\:\vek u^{h}=\hat{\vek u}^{h}\:\textrm{on}\:\Gamma_{\vek u}\right\} ,\label{eq:DisplacementTrialSpace}\\
\mathcal{V}_{\vek u}^{h} & = & \left\{ \left.\vek w^{h}\right|\vek w^{h}\in\left(\mathcal{H}^{1h}\right)^{d},\:\vek w^{h}=\vek0\:\textrm{on}\:\Gamma_{\vek u}\right\} ,\label{eq:DisplacementTestSpace}
\end{eqnarray}
where $\mathcal{H}^{1h}\subset\mathcal{H}^{1}$ is a finite dimensional
Hilbert space consisting of the shape functions. The space $\mathcal{H}^{1}$
is the set of functions which are, together with their first derivatives,
square-integrable in $\Omega$. The discretized weak form may be formulated
in the following Bubnov-Galerkin setting \cite{Belytschko_2000b,Zienkiewicz_2000d}:
Find $\vek u^{h}\in\mathcal{S}_{\vek u}^{h}$ such that
\begin{equation}
\int_{\Omega}\boldsymbol{\sigma}\left(\vek u^{h}\right):\boldsymbol{\varepsilon}\left(\vek w^{h}\right)d\Omega=\int_{\Omega}\vek w^{h}\cdot\vek f^{h}d\Omega+\int_{\Gamma_{\vek t}}\vek w^{h}\cdot\hat{\vek t}^{h}d\Gamma\quad\,\,\forall\vek w^{h}\in\mathcal{V}_{\vek u}^{h}.\label{eq:Bubnov-Galerkin-WeakFormStructure}
\end{equation}
which is the (discrete) principle of virtual work. Obviously, the
correct material parameters have to be assigned during the integration
of Eq.~(\ref{eq:Bubnov-Galerkin-WeakFormStructure}) which is based
on the signs of the involved level-set functions in each element.

\section{Numerical results\label{X_NumericalResults}}

Different test cases in two and three dimensions are considered next.
On the one hand, test cases with known analytical solutions are considered
to investigate the achieved convergence rates. On the other hand,
more technical applications aim to show the potential of the proposed
method in practice. Solutions are then compared to ``overkill solutions''
obtained on extremely fine higher-order meshes. Special attention
is given to situations where singularities are present in the solutions,
e.g., at reentrant corners of the domain. There, optimal convergence
rates can no longer be expected, however, it is found that adaptive
refinements still enable highly accurate approximations.

In the following, the errors are measured in the $L_{2}$-norm of
the displacements when analytic solutions are available. Otherwise,
it is useful to study the convergence of scalar quantities such as
the stored elastic energy or selected displacements. The condition
numbers are computed using Matlab's condest-function. The values are
normalized by dividing through the smallest condition number obtained
in a certain convergence study.

\subsection{Square shell with circular hole\label{XX_SquareWithHole}}

A square shell with dimensions $\left[-1,1\right]\times\left[-1,1\right]$
is considered with a circular void region of radius $R=0.7123$. Plane
strain conditions are assumed with Young's modulus $E=1000$ and Poisson's
ratio $\nu=0.3$. The exact solution is given in the appendix \ref{XX_ExactSolSquareWithHole}
and is also found in \cite{Szabo_1991a,Liu_2002a,Cheng_2009a}. The
corresponding displacements are prescribed along the outer boundary
of the domain, the inner boundary to the void is traction-free.

\begin{figure}
\centering

\subfigure[background mesh, quad]{\includegraphics[height=4cm]{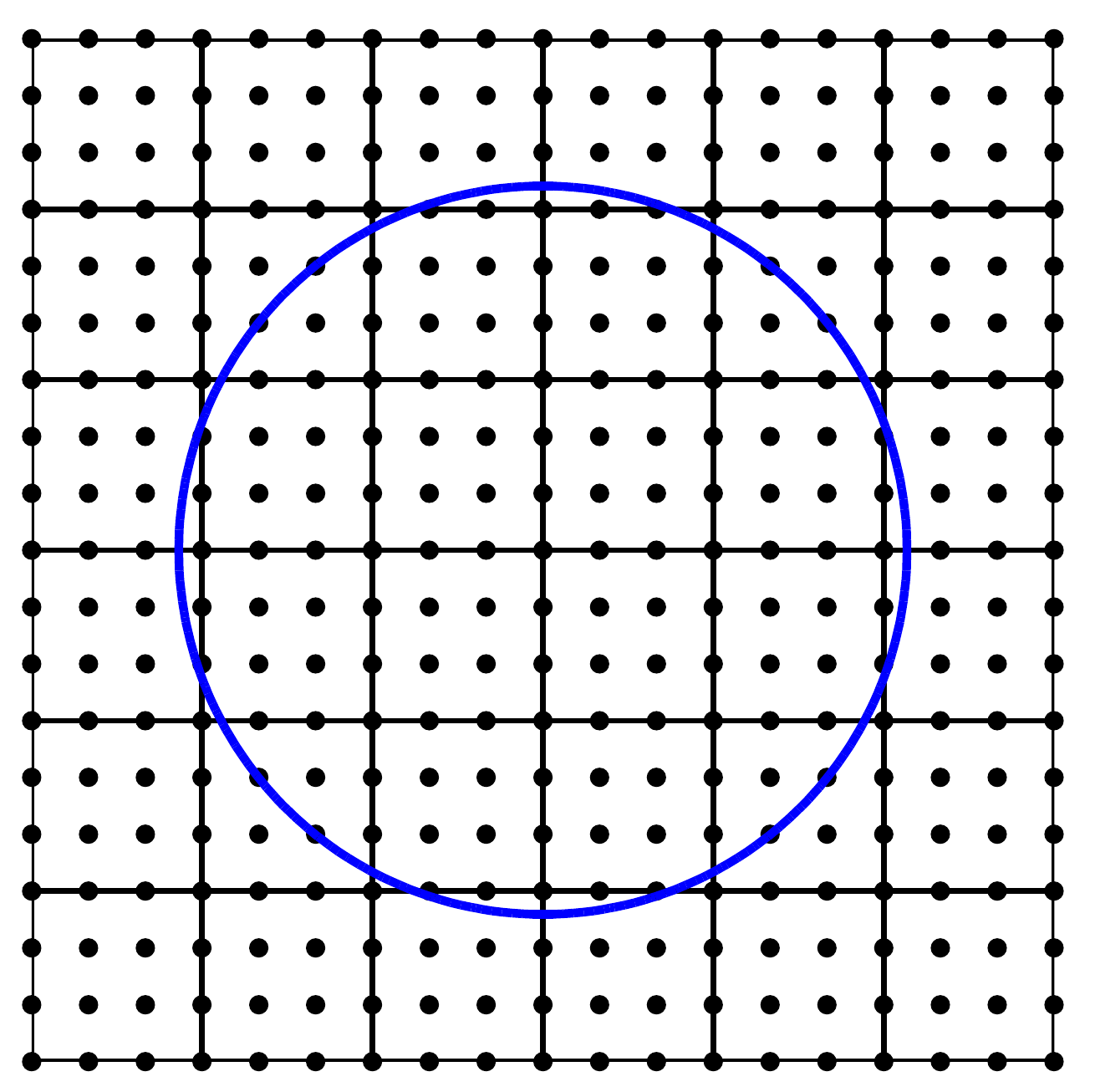}}\qquad\subfigure[background mesh, tri]{\includegraphics[height=4cm]{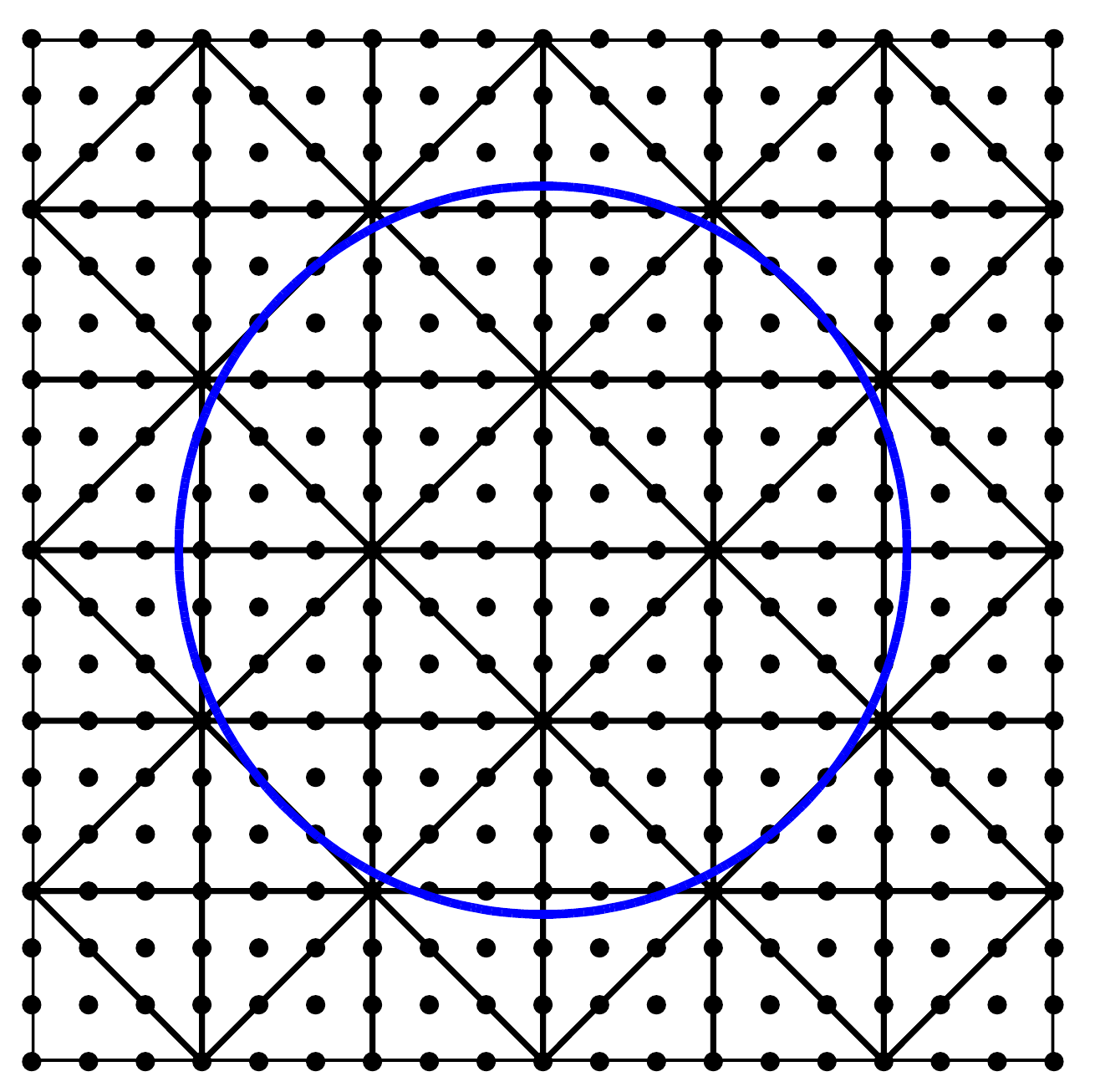}}

\subfigure[generated mesh]{\includegraphics[height=4cm]{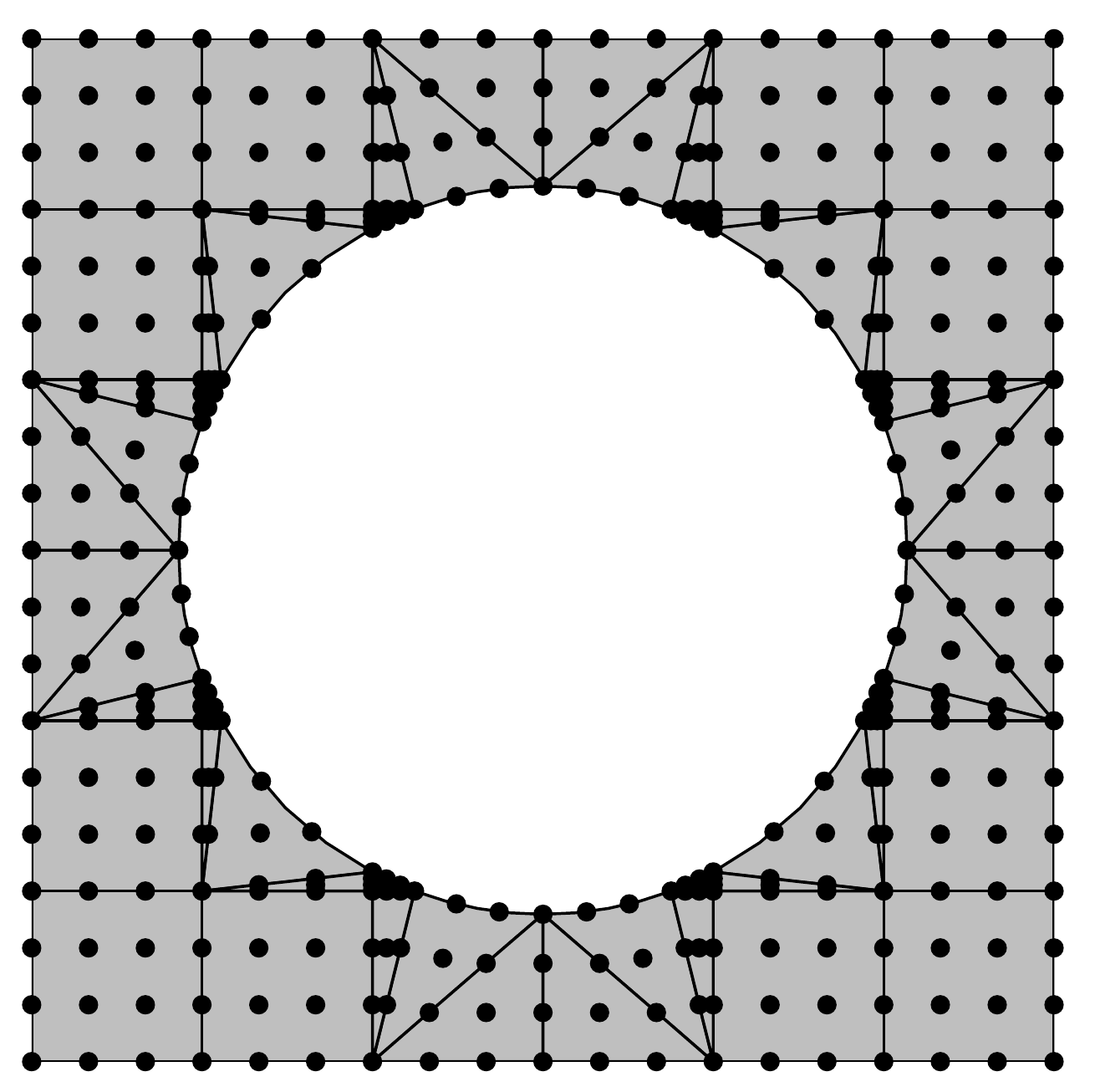}}\qquad\subfigure[exact solution]{\includegraphics[height=4cm]{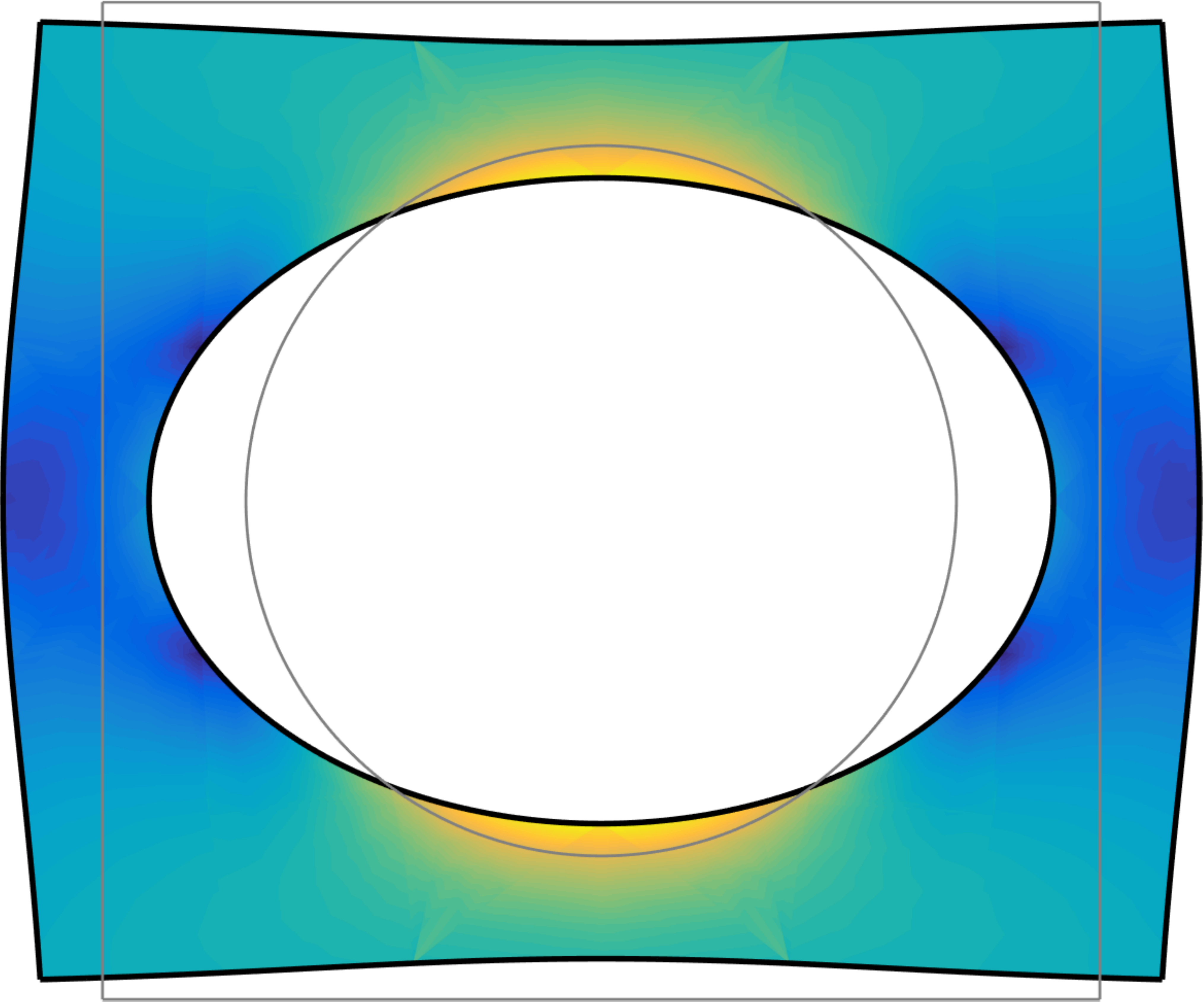}}

\caption{(a) quadrilateral and (b) triangular background meshes with circular
zero-level set, (c) generated higher-order mesh with hole, (d) deformed
configuration and von Mises stress.}

\label{fig:TestCase2dA_Sketch} 
\end{figure}

Background meshes in $\left[-1,1\right]\times\left[-1,1\right]$ with
quadrilateral and triangular elements of different orders are considered,
see Figs.~\ref{fig:TestCase2dA_Sketch}(a) and (b). The inner boundary
is defined implicitly by the level-set function
\begin{equation}
\phi\left(\vek x\right)=\sqrt{x^{2}+y^{2}}-R\label{eq:LevelSetCircle}
\end{equation}
which is evaluated at the nodes of the background mesh, so that, in
fact, only the interpolation $\phi^{h}\left(\vek x\right)$ is used.
An example of an automatically generated mesh based on a background
mesh with $6\times6$ cubic elements is seen in Fig.~\ref{fig:TestCase2dA_Sketch}(c).
The exact solution is plotted in Fig.~\ref{fig:TestCase2dA_Sketch}(d)
in terms of the deformed configuration and the corresponding von Mises
stress.

\begin{figure}
\centering

\subfigure[integration error, quad]{\includegraphics[height=4cm]{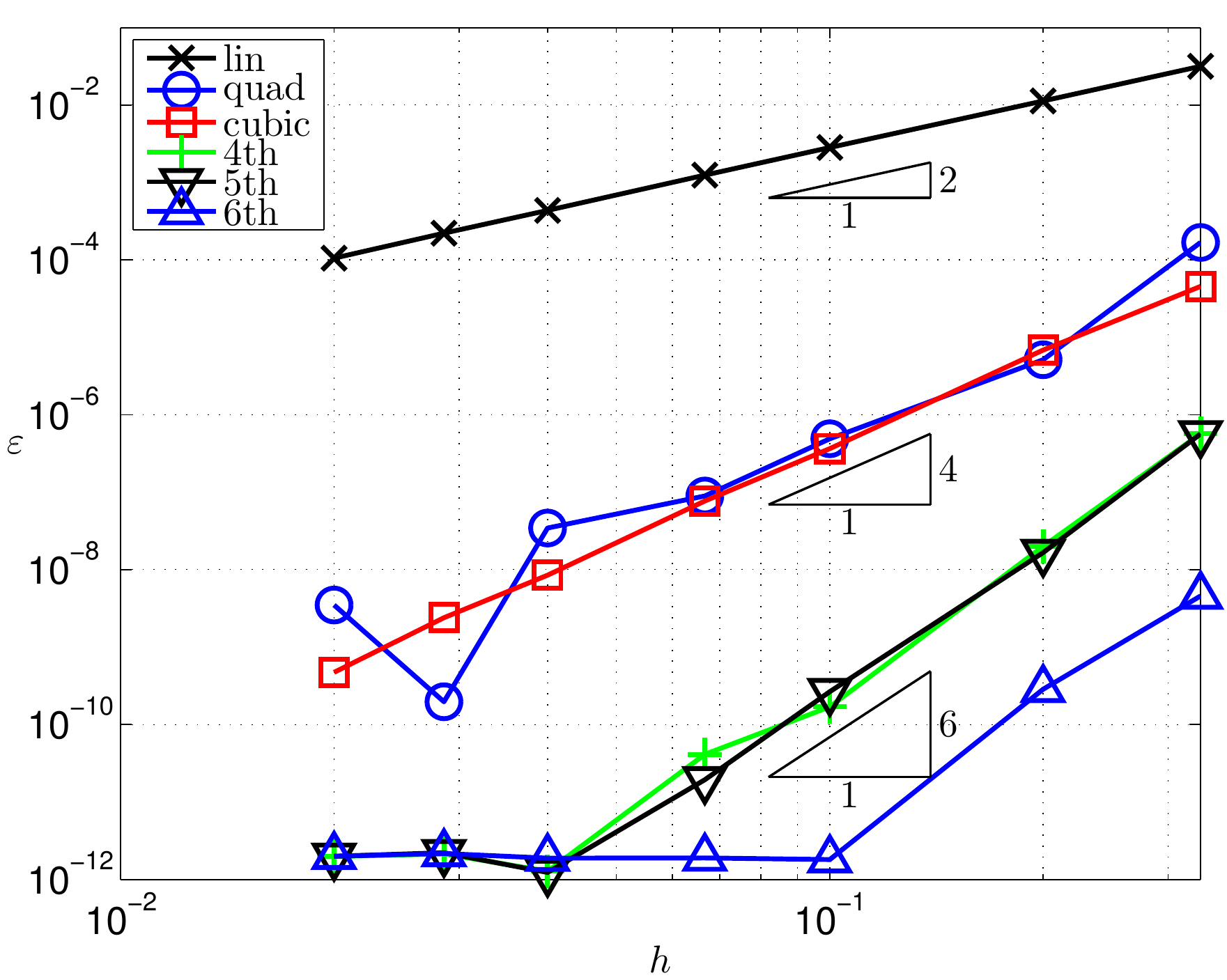}}\quad\subfigure[interpolation error, quad]{\includegraphics[height=4cm]{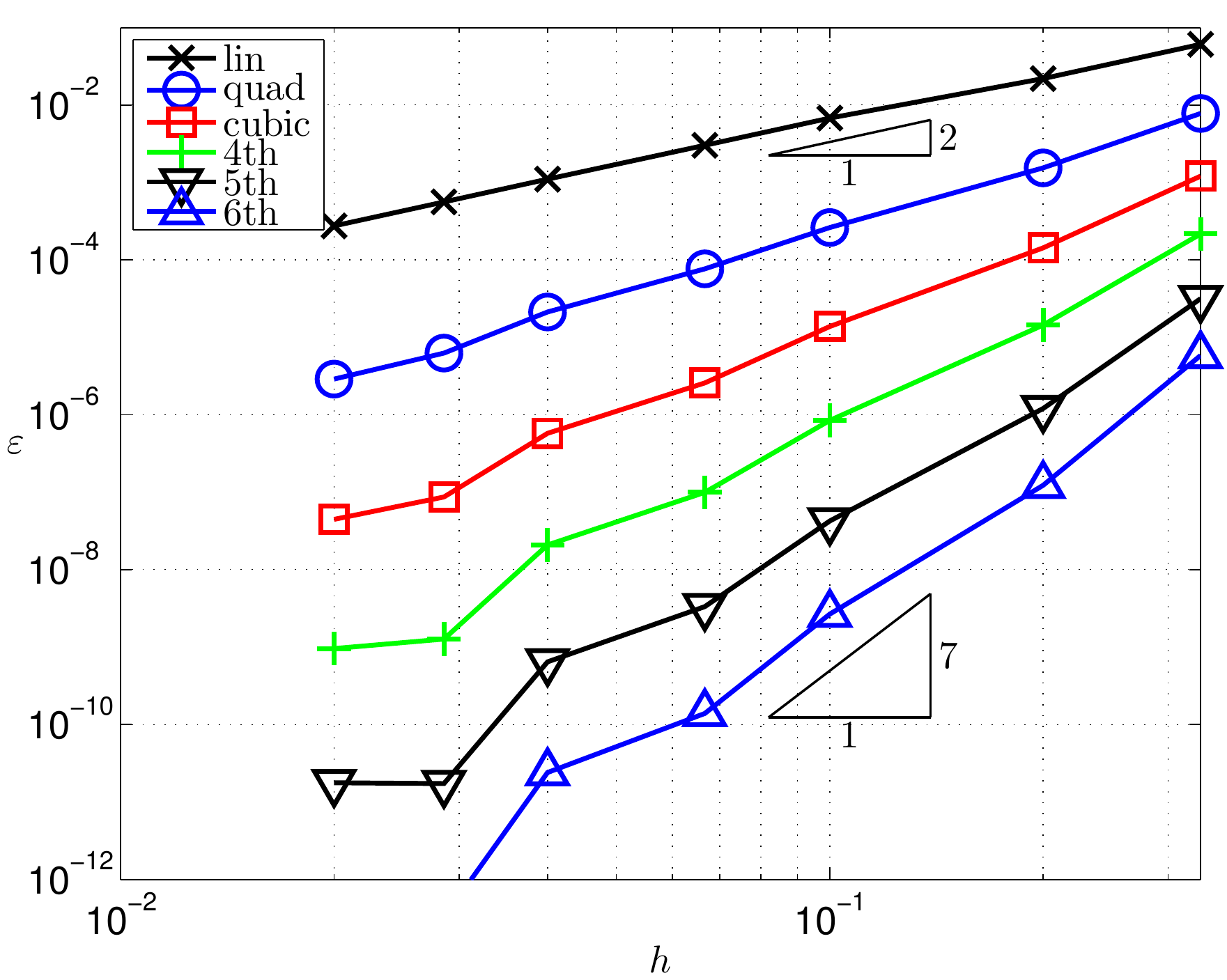}}\quad\subfigure[approximation error, quad]{\includegraphics[height=4cm]{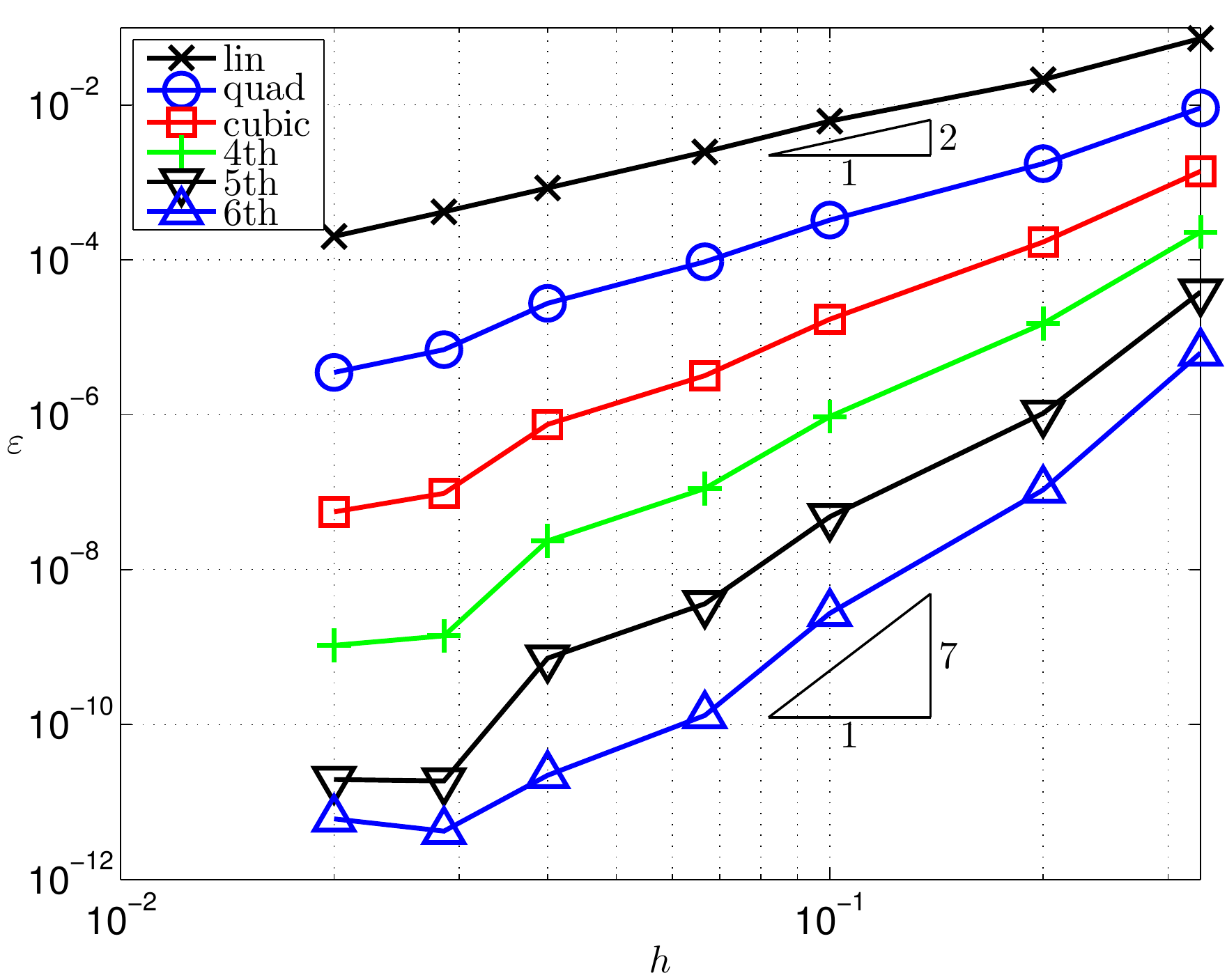}}

\subfigure[approximation error, tri]{\includegraphics[height=4cm]{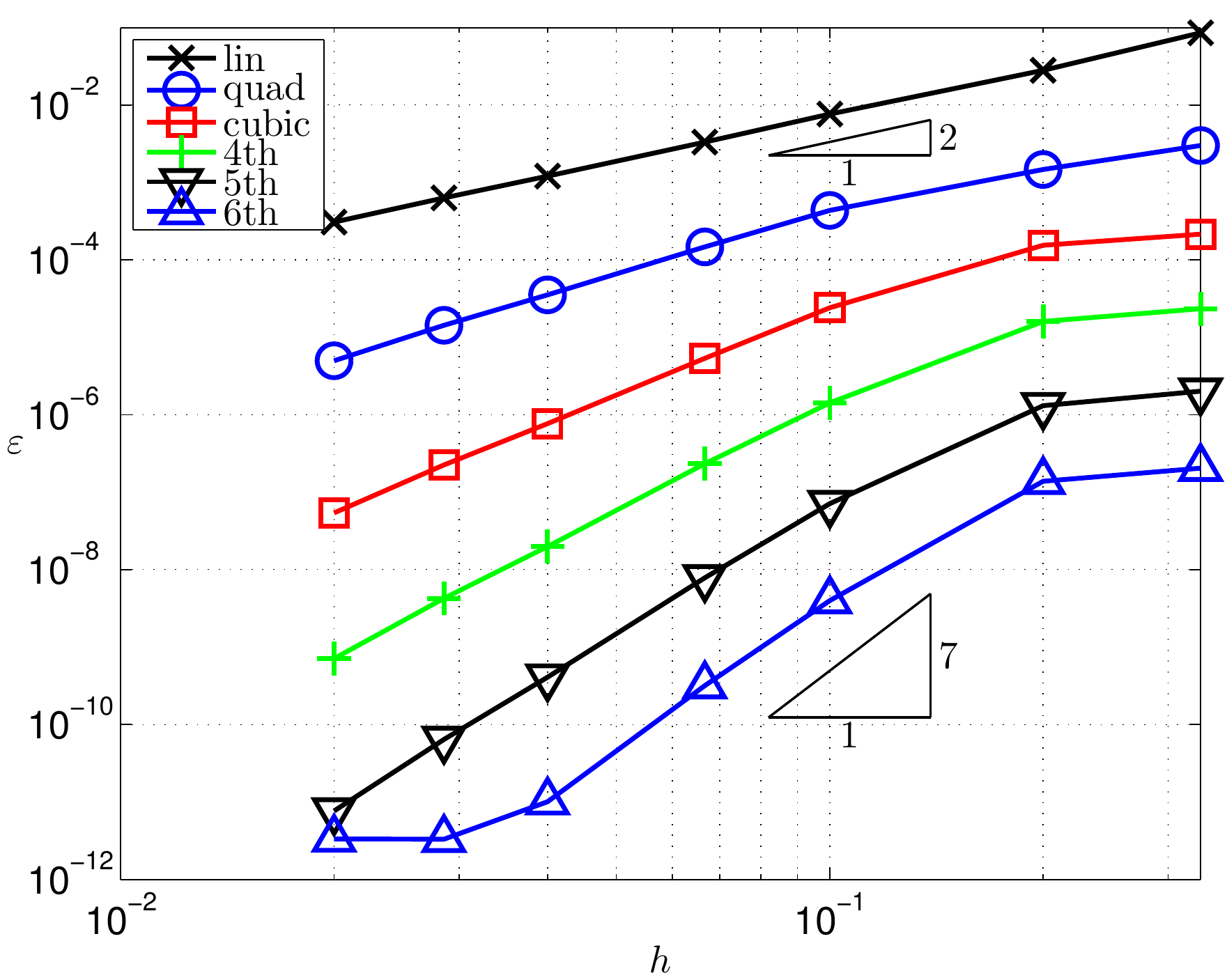}}\quad\subfigure[condition number, quad]{\includegraphics[height=4cm]{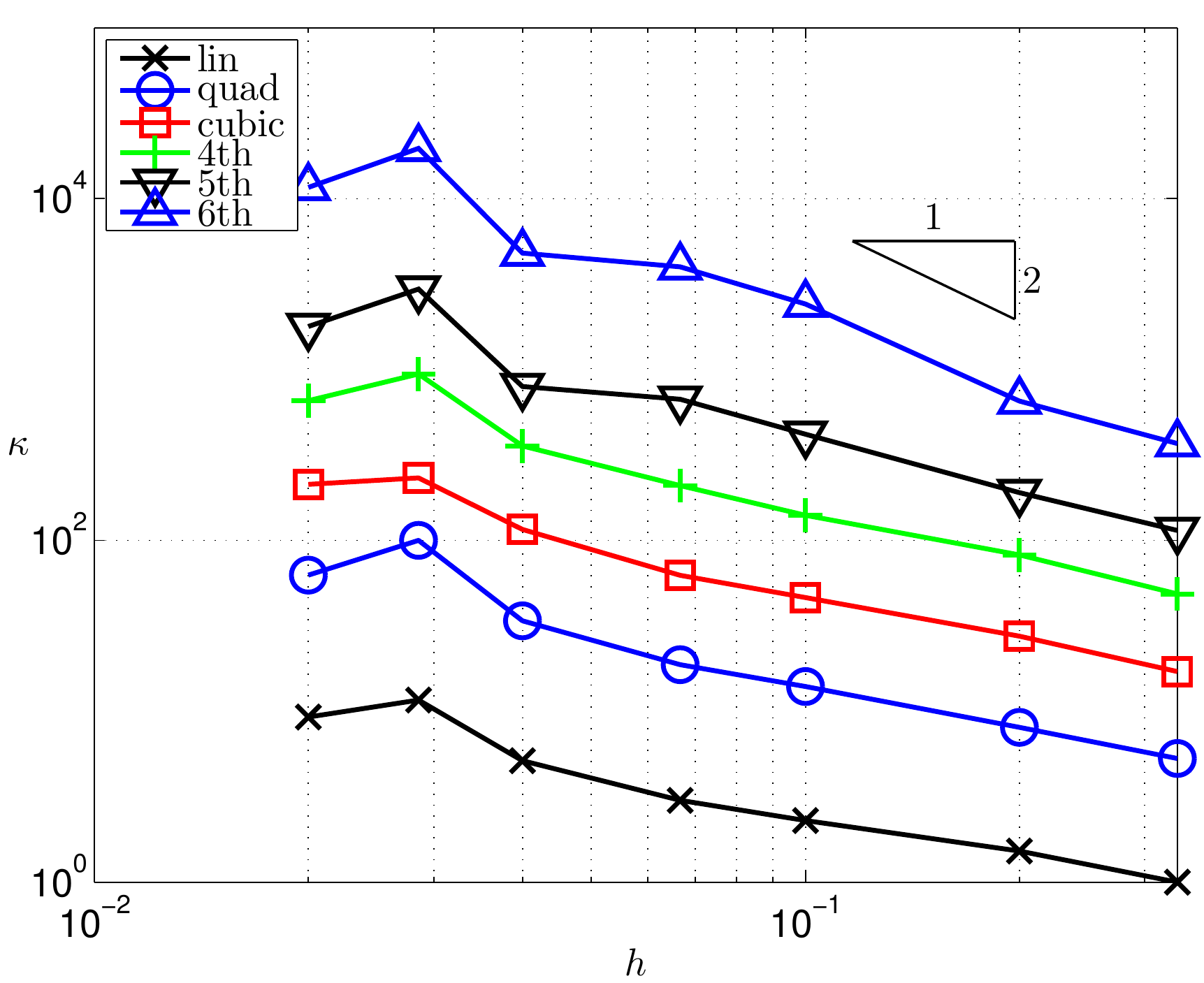}}\quad\subfigure[condition number, tri]{\includegraphics[height=4cm]{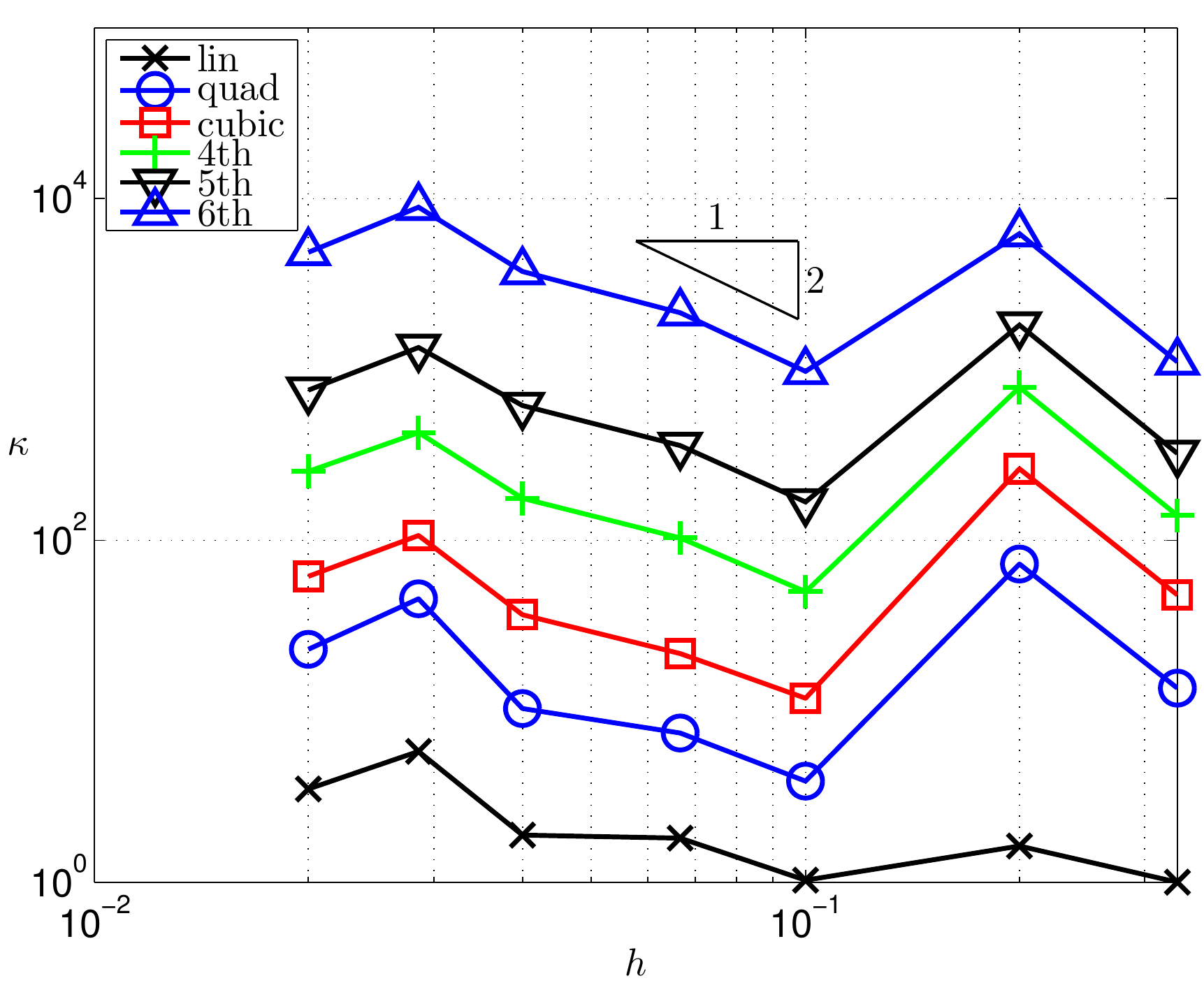}}

\caption{(a) integration, (b) interpolation, and approximation error for (c)
quadrilateral or (d) triangular background meshes, (e) and (f) show
the corresponding condition numbers of the system matrices, respectively.}

\label{fig:TestCase2dA_Res} 
\end{figure}

For the convergence study, the number of elements, $n_{d}$, per dimensions
of the background mesh is systematically increased and $n_{d}=\left\{ 6,10,20,30,50,70,100\right\} $
elements are used with varying orders between $1$ and $6$. Results
are shown in Fig.~\ref{fig:TestCase2dA_Res}. In \cite{Fries_2015a,Fries_2016a},
the focus is on the integration properties of the automatically generated
meshes. Therefore, for example, the area of the mesh may be computed
and compared to the exact area. For this example, this is shown in
Fig.~\ref{fig:TestCase2dA_Res}(a) and optimal convergence rates
are achieved. In \cite{Fries_2016b}, also the interpolation error
is studied, i.e., the ability to reproduce functions on the generated
meshes. The error between a given example function and its interpolation
is shown in Fig.~\ref{fig:TestCase2dA_Res}(b) and is, again, optimal.
It is clear that the ability to integrate and interpolate optimally
is a necessary requirement for an optimal convergence in an approximation
context as well. Convergence results of the approximated displacements
in the $L_{2}$-norm are shown in Figs.~\ref{fig:TestCase2dA_Res}(c)
and (d) based on background meshes composed by quadrilateral and triangular
elements, respectively. The corresponding condition numbers $\kappa$
of the resulting system matrices for the different meshes are seen
in Figs.~\ref{fig:TestCase2dA_Res}(e) and (f).

\subsection{Square shell with circular inclusion\label{XX_SquareWithCircInclusion}}

\begin{figure}
\centering

\subfigure[generated mesh, $n_{d}=6$]{\includegraphics[height=4cm]{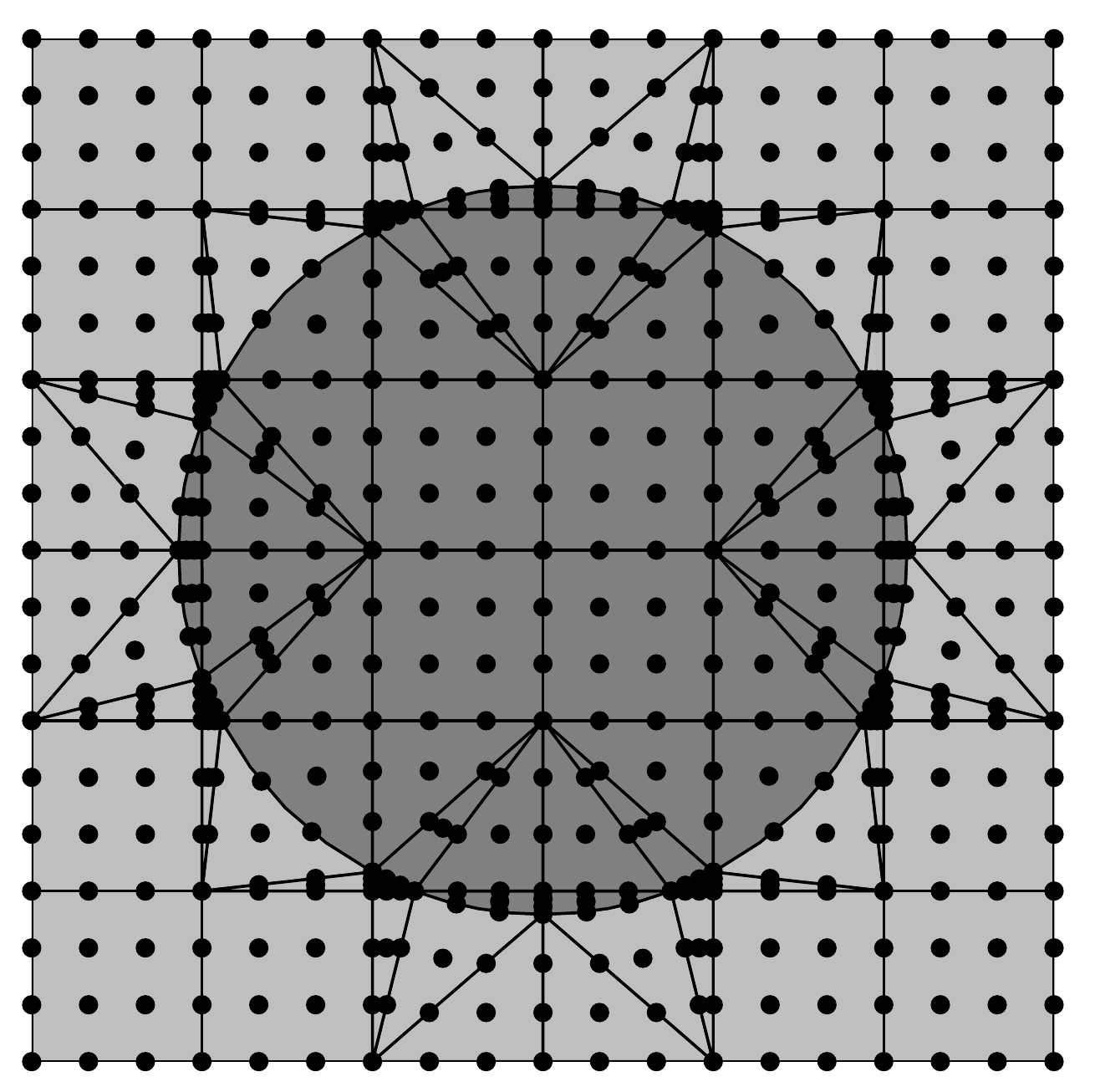}}\qquad\subfigure[generated mesh, $n_{d}=10$]{\includegraphics[height=4cm]{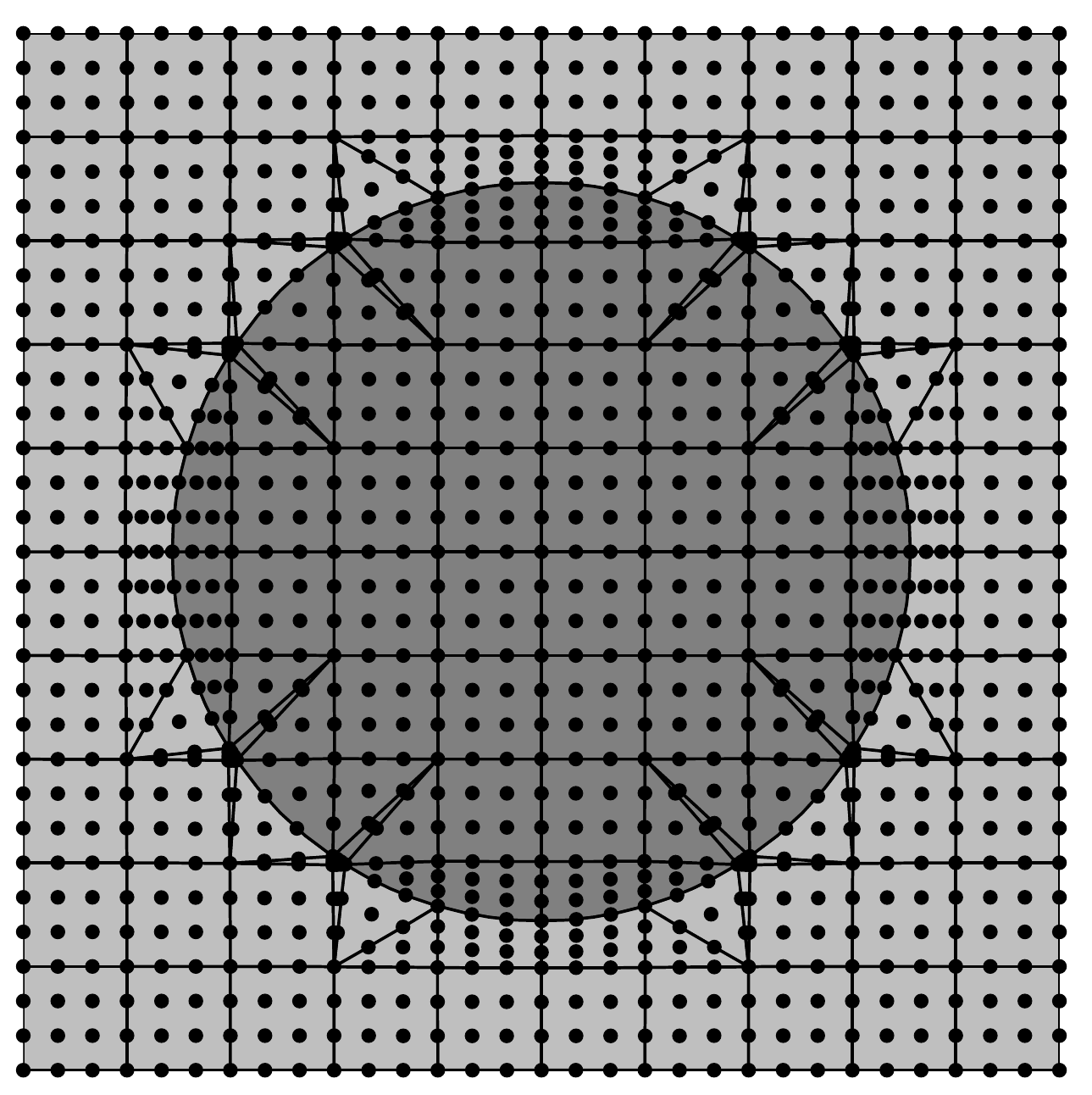}}

\subfigure[generated mesh, $n_{d}=20$]{\includegraphics[height=4cm]{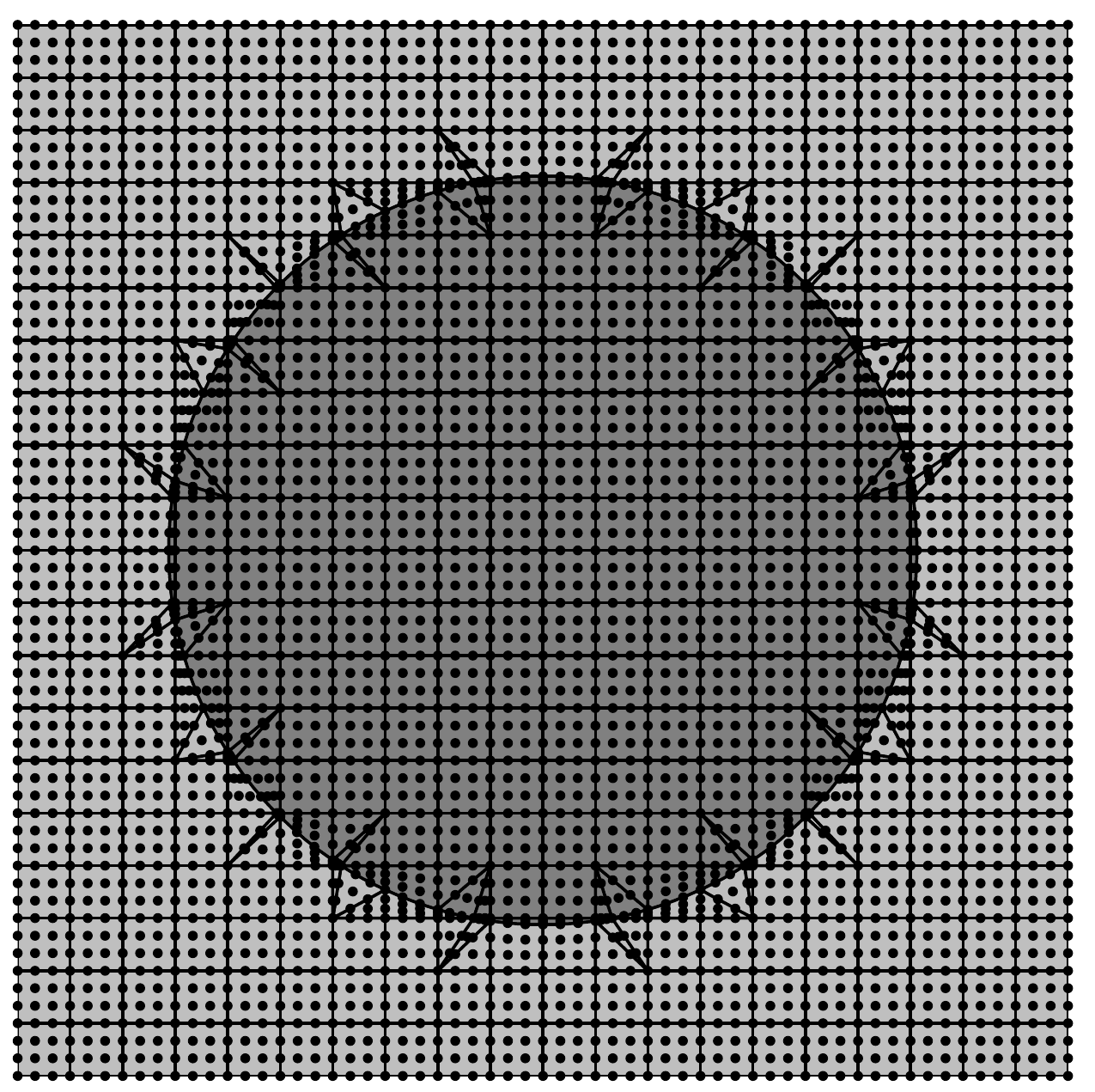}}\qquad\subfigure[exact solution]{\includegraphics[height=4cm]{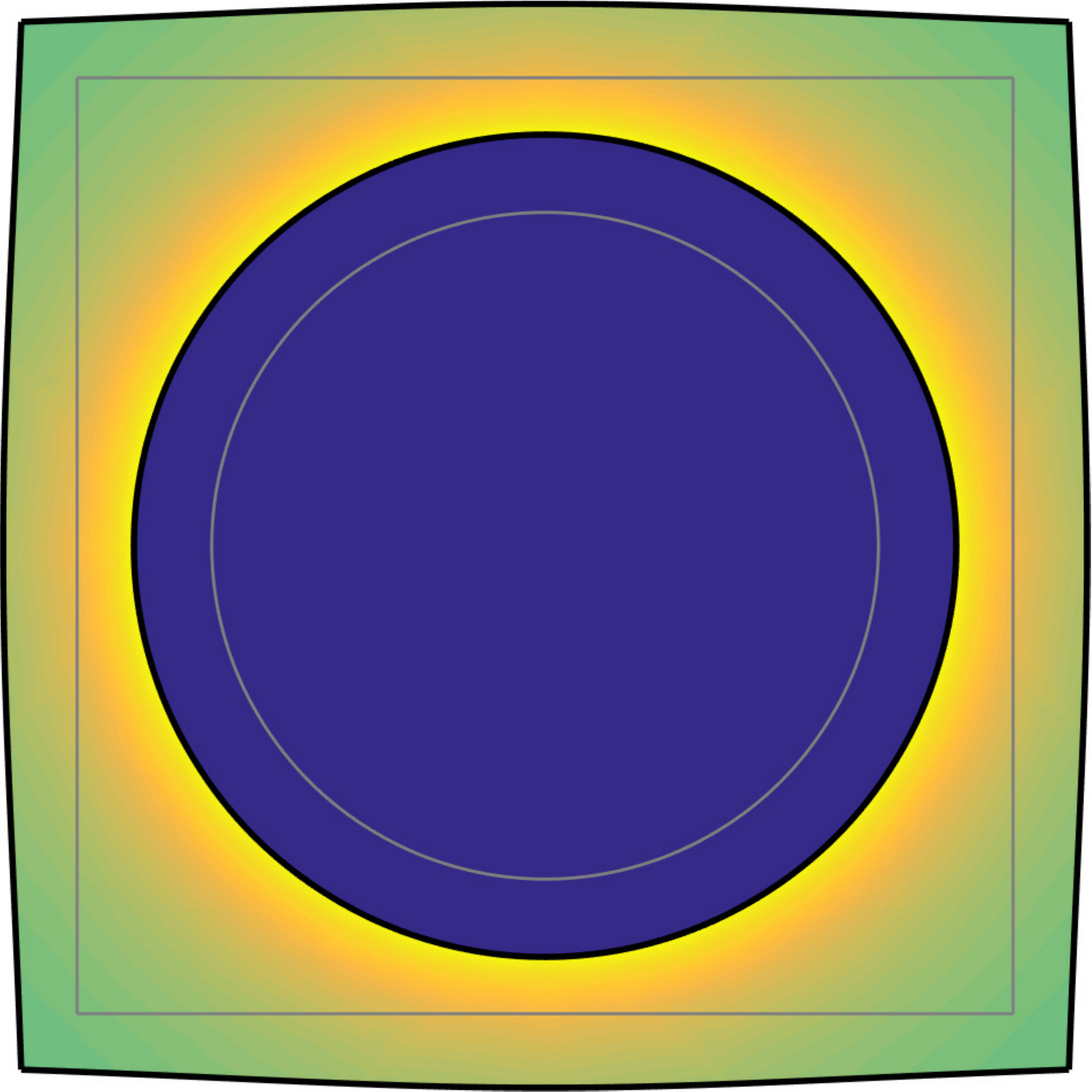}}

\caption{(a) to (c) show generated meshes based on different resolutions of
the background mesh (composed by cubic quadrilateral elements), (d)
deformed configuration and von Mises stress.}

\label{fig:TestCase2dB_Sketch} 
\end{figure}

A shell with the same geometry from above is considered, however,
the void region is now filled with a different material. That is,
the domain is again $\left[-1,1\right]\times\left[-1,1\right]$ and
the interface is defined by the zero-level set of Eq.~(\ref{eq:LevelSetCircle}).
In the outer region, Young's modulus is $E_{1}=10$ and Poisson's
ratio $\nu_{1}=0.3$. Inside the circular inclusion, there is $E_{2}=1$
and $\nu_{2}=0.25$. An exact solution for this problem is found in
\cite{Sukumar_2001c,Fries_2006b} and is given in the Appendix \ref{XX_ExactSolSquareWithCircInclusion}.
The deformed configuration with von Mises stress is seen in Fig.~\ref{fig:TestCase2dB_Sketch}(d).
The corresponding displacements are prescribed at the outer boundary.
For the convergence studies, the same background meshes with different
numbers of elements per dimensions $n_{d}$ and element orders than
in Section \ref{XX_SquareWithHole} are used. Examples for the generated
conforming meshes based on background meshes with $n_{d}=\left\{ 6,10,20\right\} $
and cubic elements are seen in Fig.~\ref{fig:TestCase2dB_Sketch}(a)
to (c), respectively.

\begin{figure}
\centering

\subfigure[approximation error, quad]{\includegraphics[height=4cm]{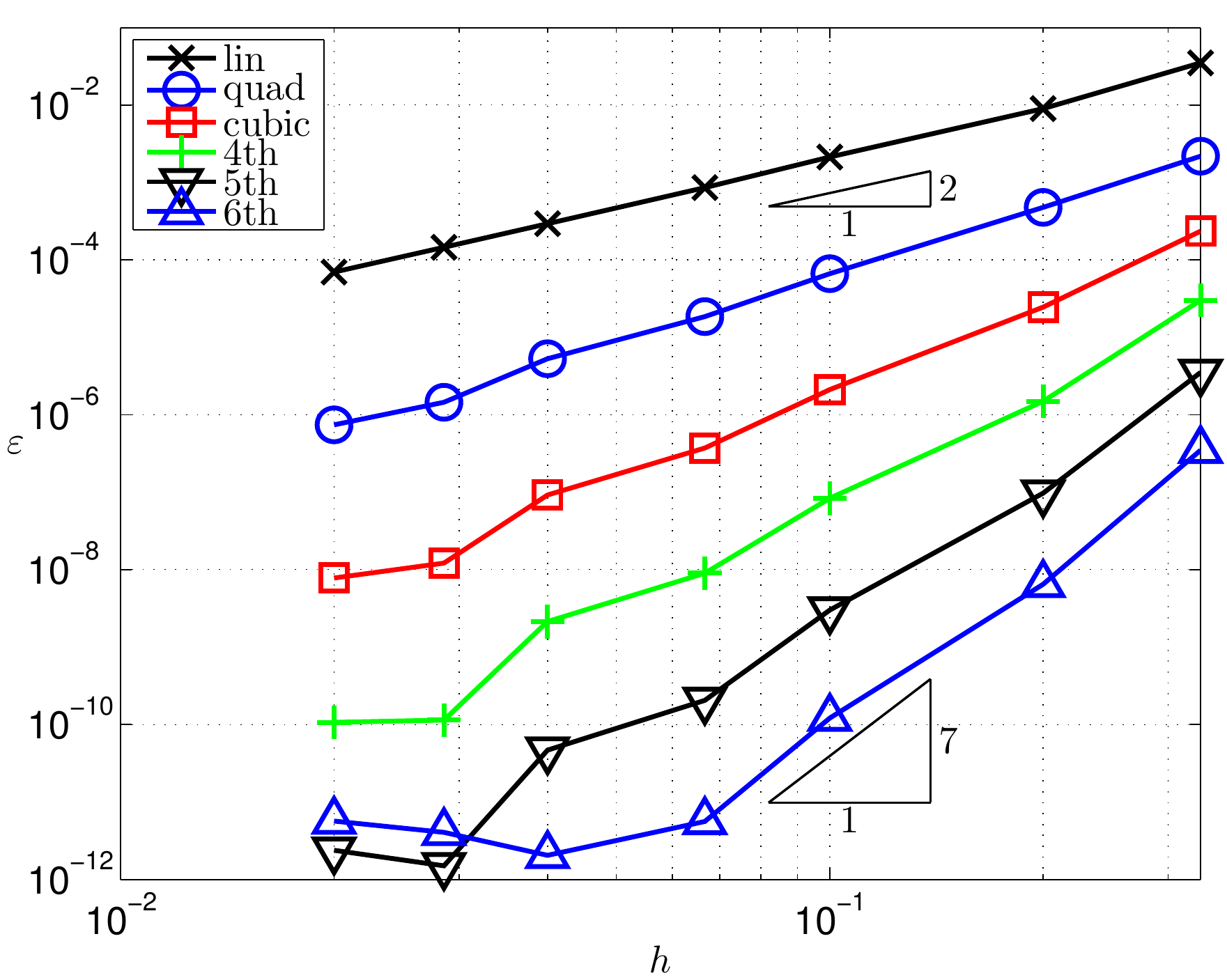}}\quad\subfigure[approximation error, tri]{\includegraphics[height=4cm]{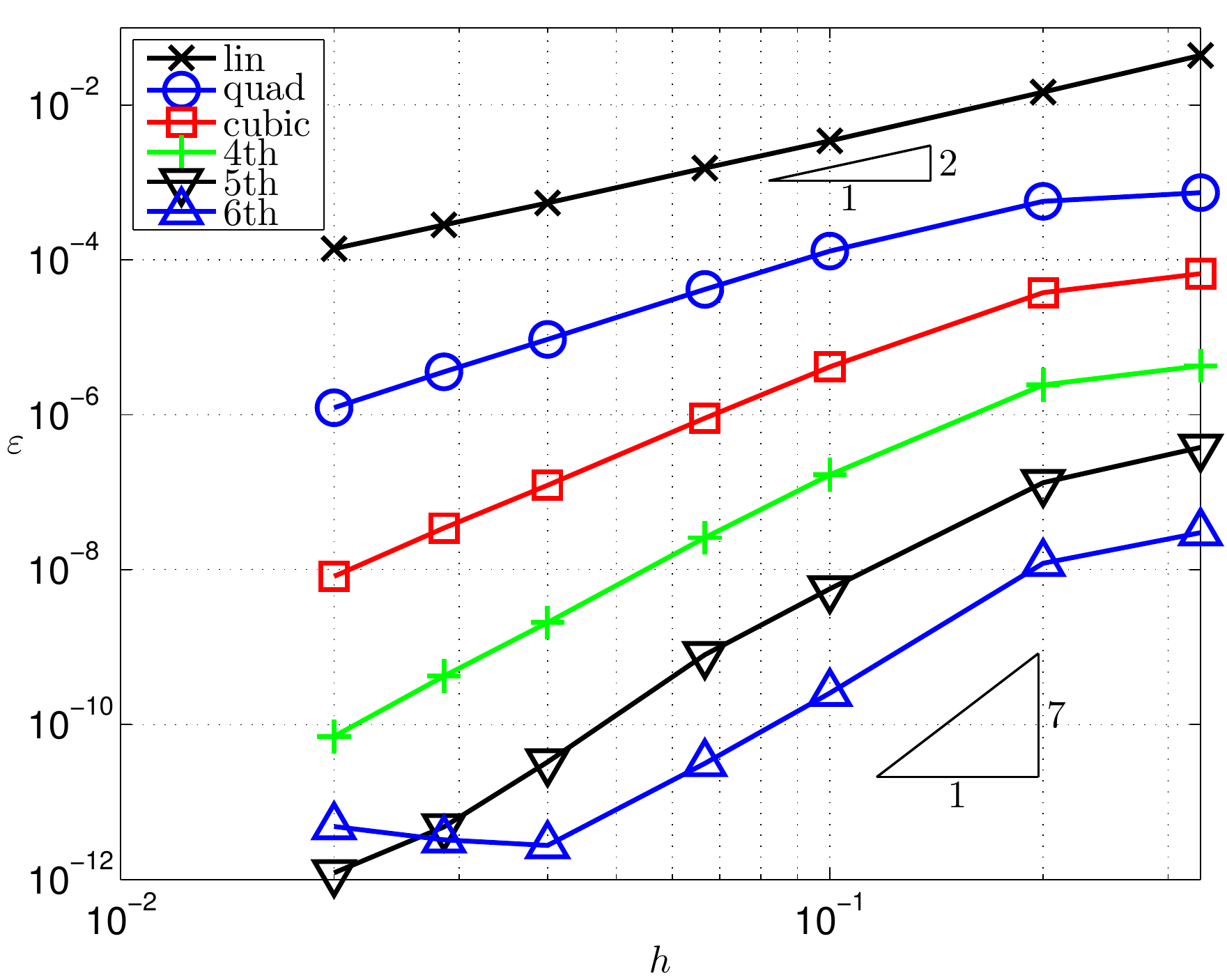}}

\subfigure[condition number, quad]{\includegraphics[height=4cm]{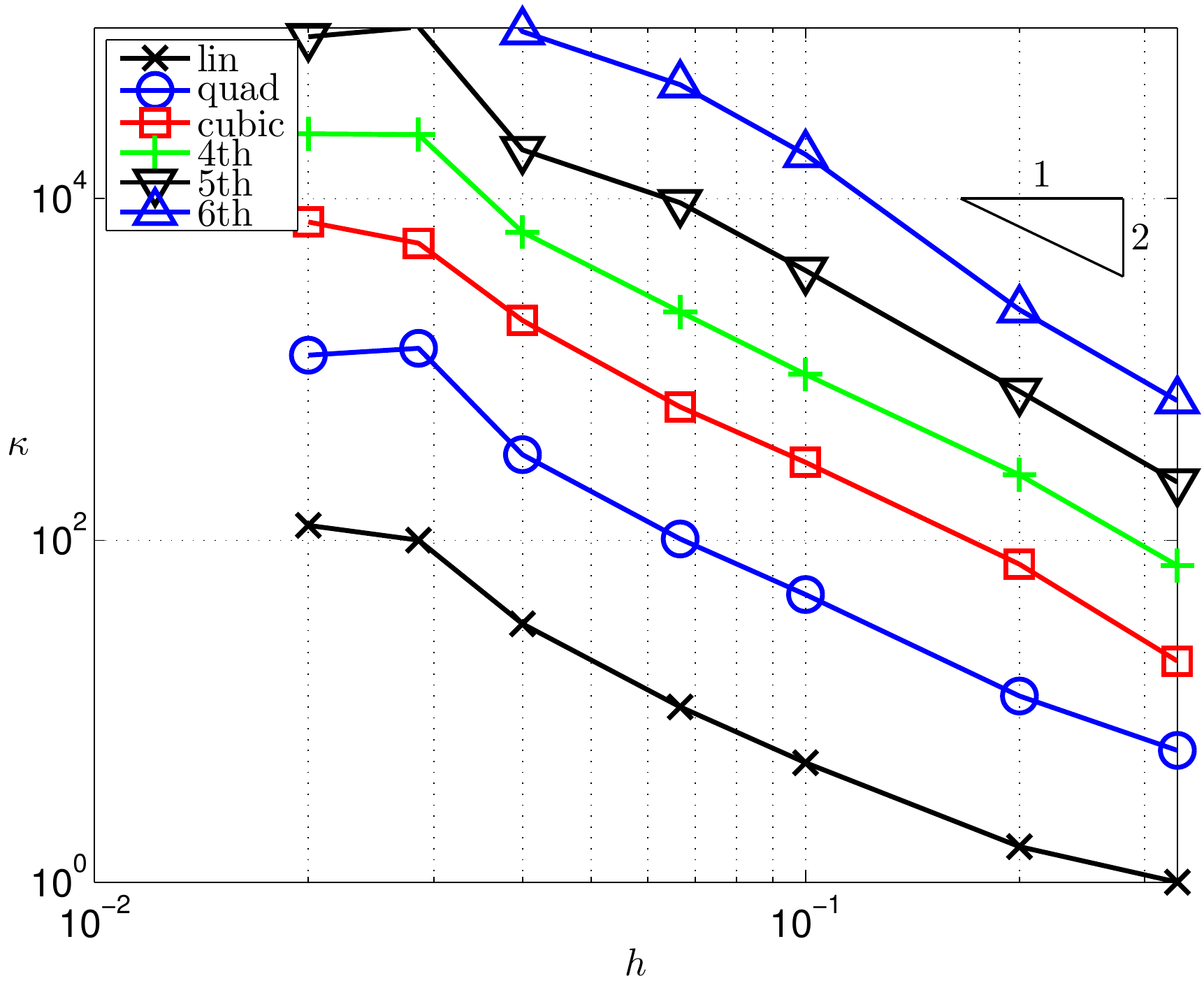}}\quad\subfigure[condition number, tri]{\includegraphics[height=4cm]{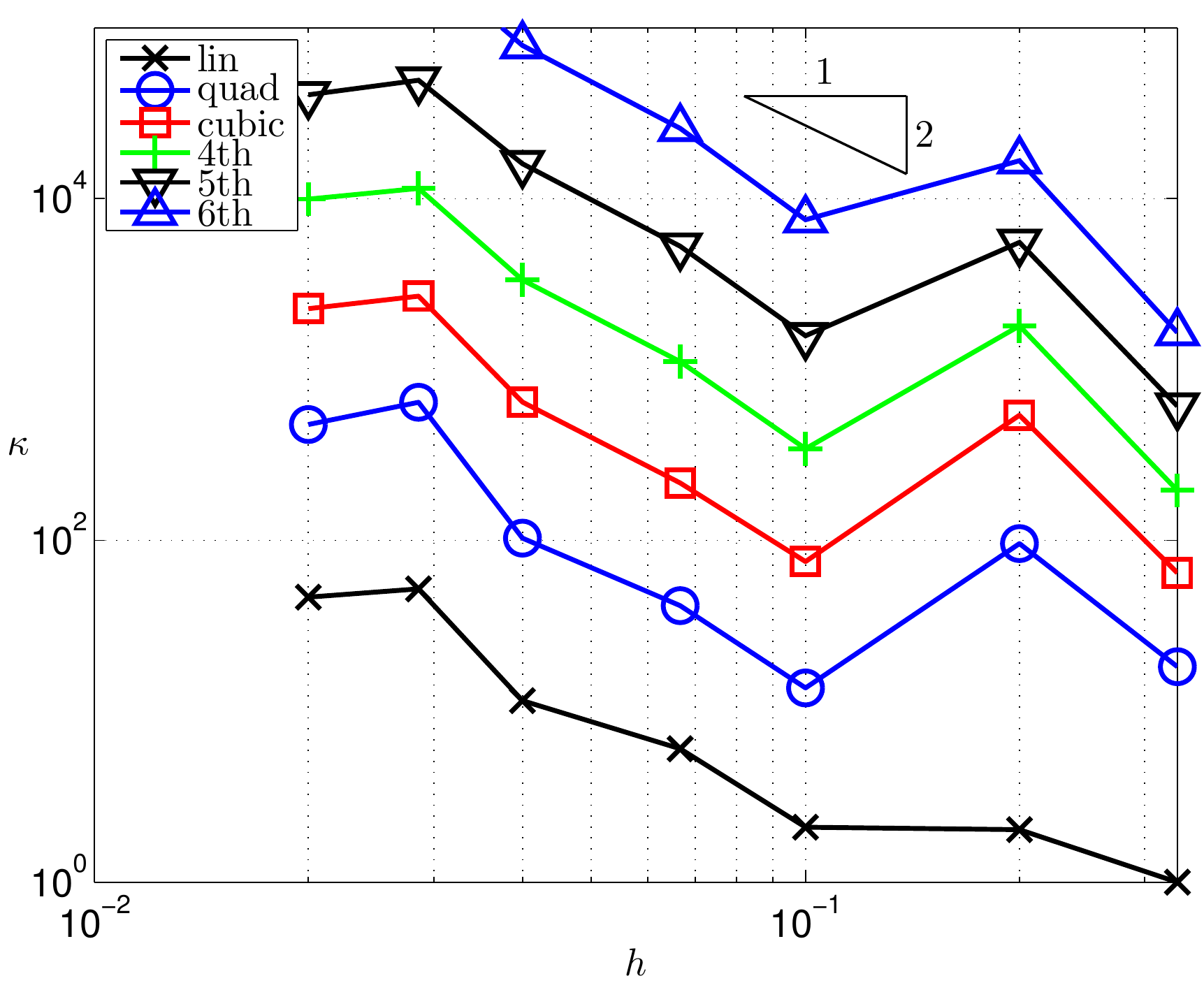}}

\caption{The approximation error for (a) quadrilateral or (b) triangular background
meshes, (c) and (d) show the corresponding condition numbers of the
system matrices, respectively.}

\label{fig:TestCase2dB_Res} 
\end{figure}

Convergence results are presented in Fig.~\ref{fig:TestCase2dB_Res}.
Integration and interpolation errors are no longer considered and
the focus is only the approximation error of the displacements in
the $L_{2}$-norm. Figs.~\ref{fig:TestCase2dB_Res}(a) and (b) show
optimal convergence rates for background meshes composed by quadrilateral
and triangular elements up to order $6$, respectively. The corresponding
condition numbers are shown in Figs.~\ref{fig:TestCase2dB_Res}(c)
and (d). It is seen that they are well-bounded, however, not as smooth
as for manufactured meshes. Nevertheless, they behave with $O\left(h^{2}\right)$
as expected.

\subsection{Cantilever beam\label{XX_CantileverBeam}}

The next test cases in two dimensions serve the purpose to demonstrate
the proposed higher-order CDFEM for more technical rather than academic
setups. A cantilever beam with length $L=5.0\mathrm{m}$ and a variable
thickness between $h=0.2\mathrm{m}$ and $0.4\mathrm{m}$ is considered
first. The material is composed of steel with $E=2.1\cdot10^{8}\,\unitfrac{kN}{m^{2}}$
and $\nu=0.3$. The beam is loaded by gravity acting as a body force
of $f_{y}=-78.5\,\unitfrac{kN}{m^{3}}$ and a vertical traction on
the right side. This traction is distributed in a quadratic profile
being zero at the upper and lower right side and reaching a maximum
of $\sigma_{y}=-100\unitfrac{kN}{m}$ inbetween, leading to a force
resultant of $F_{y}=-26.\bar{6}\,\mathrm{kN}$.

\begin{figure}
\centering

\subfigure[background mesh]{\includegraphics[width=0.8\textwidth]{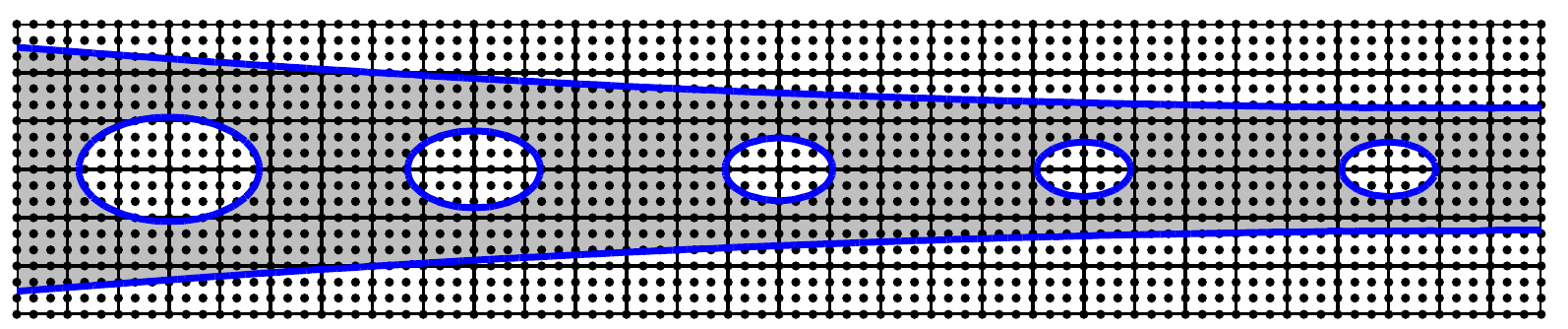}}

\subfigure[generated mesh, coarse]{\includegraphics[width=0.8\textwidth]{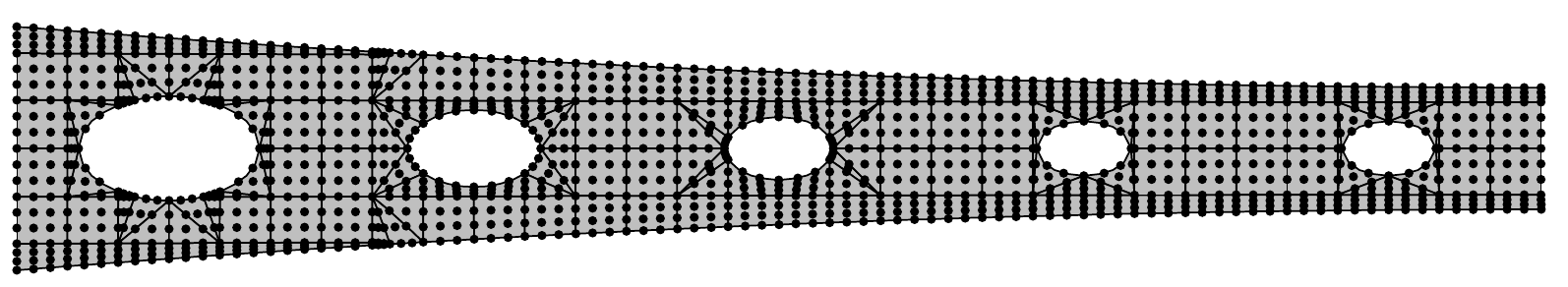}}

\subfigure[generated mesh, finer]{\includegraphics[width=0.8\textwidth]{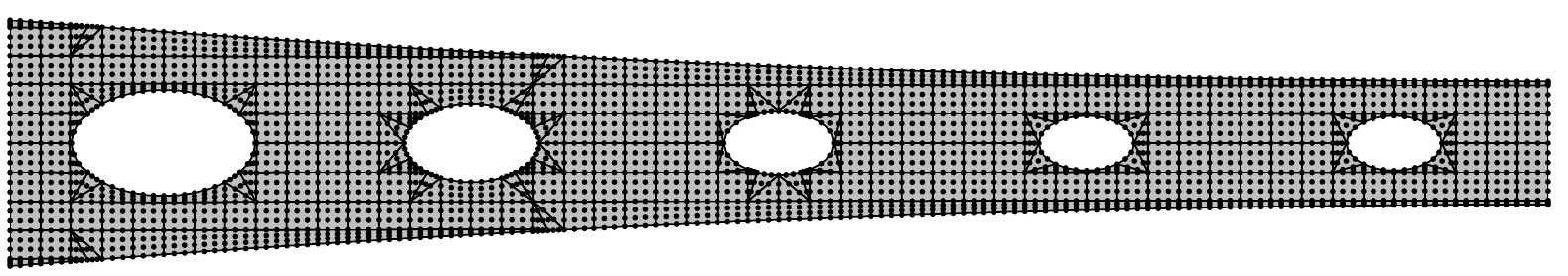}}

\caption{(a) Background mesh and zero-level sets of the $7$ level-set functions,
(b) and (c) show examples of generated conforming higher-order meshes.}

\label{fig:TestCase2dC_Sketch} 
\end{figure}

The beam features $5$ elliptical void regions. We place the $xy$-coordinate
system at the left side in the middle axis of the beam. $7$ level-set
functions are used to define the geometry and the background meshes
conform to the left and right side from the beginning. See Fig.~\ref{fig:TestCase2dC_Sketch}(a)
for the zero-level sets and an example background mesh. The $5$ elliptical
void regions are defined from left to right as
\[
\phi_{i}\left(\vek x\right)=a\cdot\left(x-x_{i}\right)^{2}+b\cdot y^{2}-R_{i}
\]
with $a=\nicefrac{1}{4}$, $b=\nicefrac{3}{4}$, $x_{i}=\left\{ 1,3,5,7,9\right\} \cdot\nicefrac{1}{2}$,
and $R_{i}=\left\{ 22,12,8,6,6\right\} \cdot10^{-3}$. The $2$ level-set
functions which define the upper and lower side of the beam are given
as
\[
\phi_{6}\left(\vek x\right)=g\left(x\right)-y,\,\phi_{7}\left(\vek x\right)=-g\left(x\right)-y\quad\text{with}\quad g\left(x\right)=\dfrac{x^{2}}{125}-\dfrac{2\cdot x}{25}+\dfrac{2}{5}.
\]
None of these level-set functions has signed-distance property. Examples
for automatically generated, conforming higher-order meshes for this
test case are seen in Fig.~\ref{fig:TestCase2dC_Sketch}(b) and (c).

We consider two different support cases: In case 1, the left side
of the beam is fully fixed. That is, zero-displacements are prescribed
for the horizontal and vertical displacements at all nodes on the
left. The deformed configuration is seen in Fig.~\ref{fig:TestCase2dC_Sketch}(a)
and the resulting von Mises stress in 3D view are shown in Fig.~\ref{fig:TestCase2dC_Sketch}(b).
It is seen that the stresses are \emph{singular} at the upper and
lower left corners, where the boundary conditions change from supported
to free. This is well known in structural mechanics and will effect
the convergence properties in a higher-order FEM as confirmed below.
Support case 2 fixes all nodes on the boundary to the left elliptical
void region. The resulting deformed configuration and von Mises stress
are seen in Figs.~\ref{fig:TestCase2dC_Sketch}(c) and (d). As seen,
there are no singularities in the stresses for this support case and
optimal convergence rates in the analysis are possible.

\begin{figure}
\centering

\subfigure[case 1, 2D view]{\includegraphics[width=7cm]{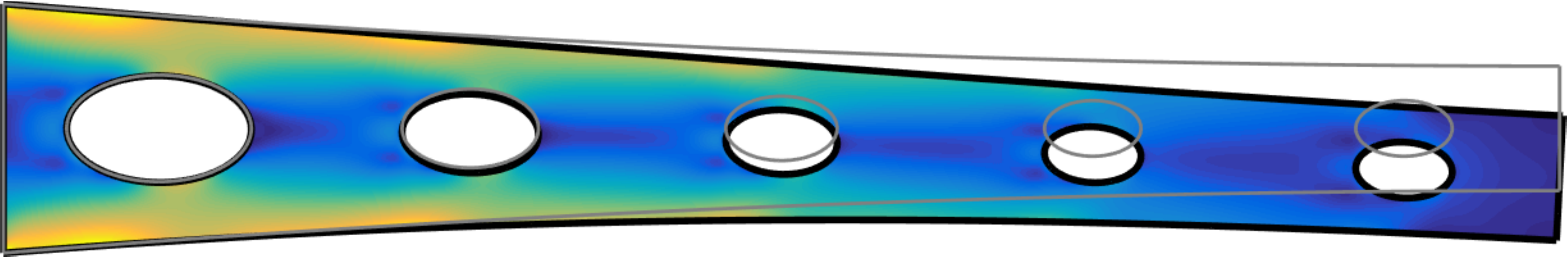}}\qquad\subfigure[case 1, 3D view]{\includegraphics[width=6cm]{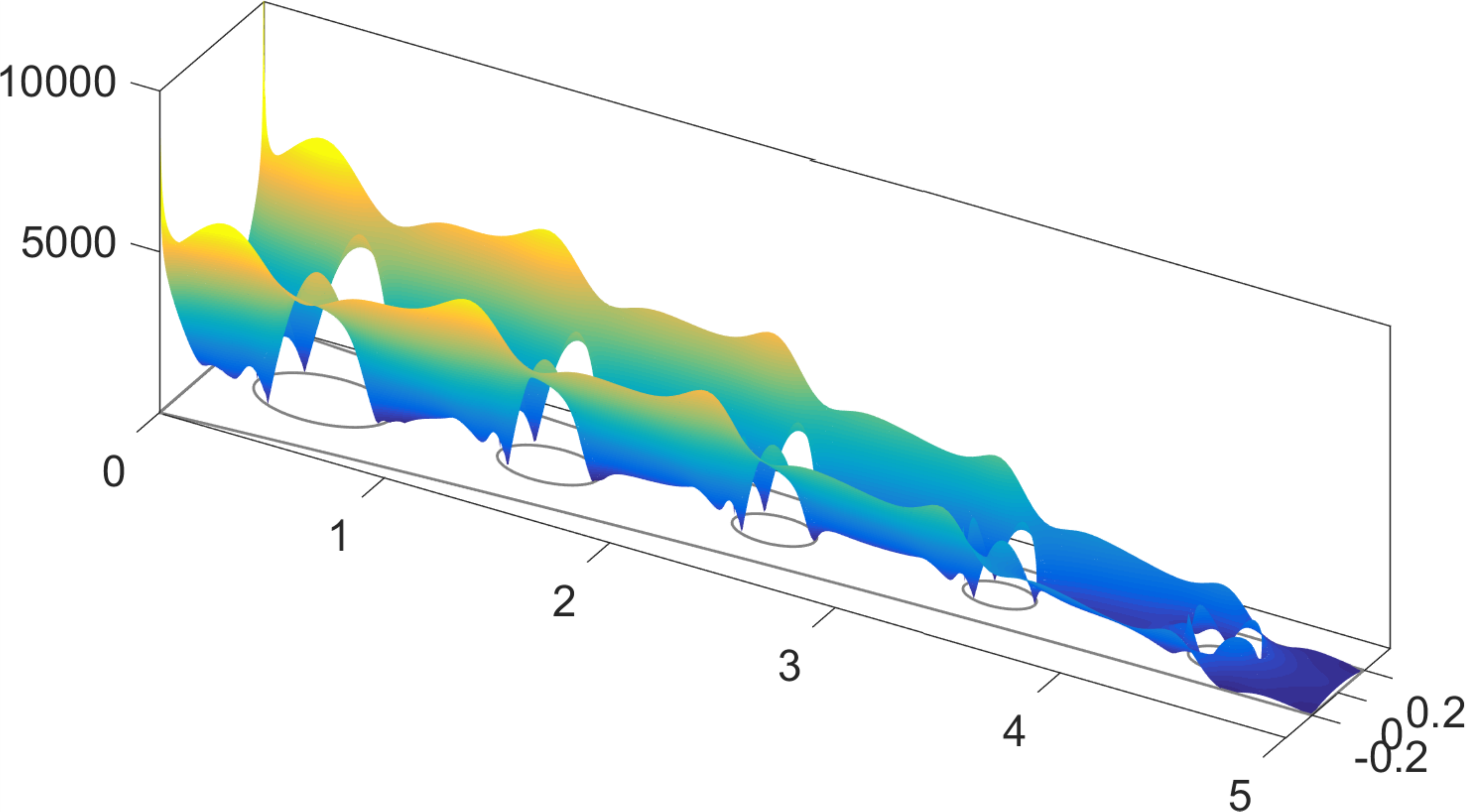}}

\subfigure[case 2, 2D view]{\includegraphics[width=7cm]{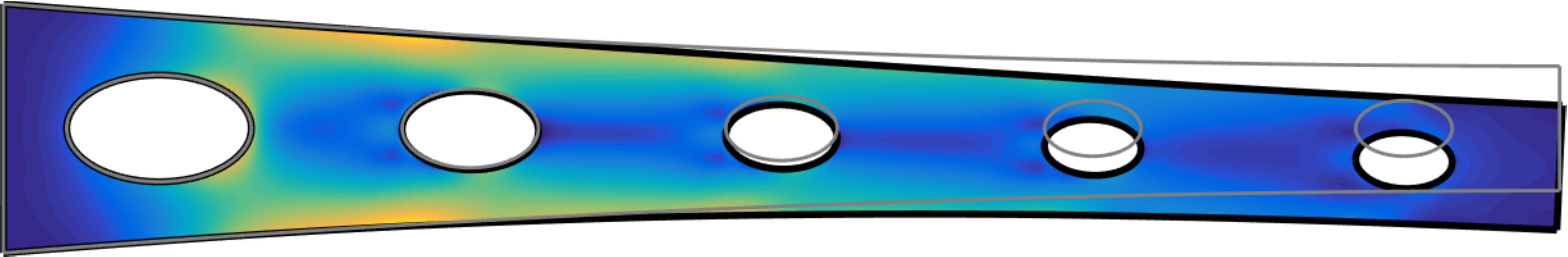}}\qquad\subfigure[case 2, 3D view]{\includegraphics[width=6cm]{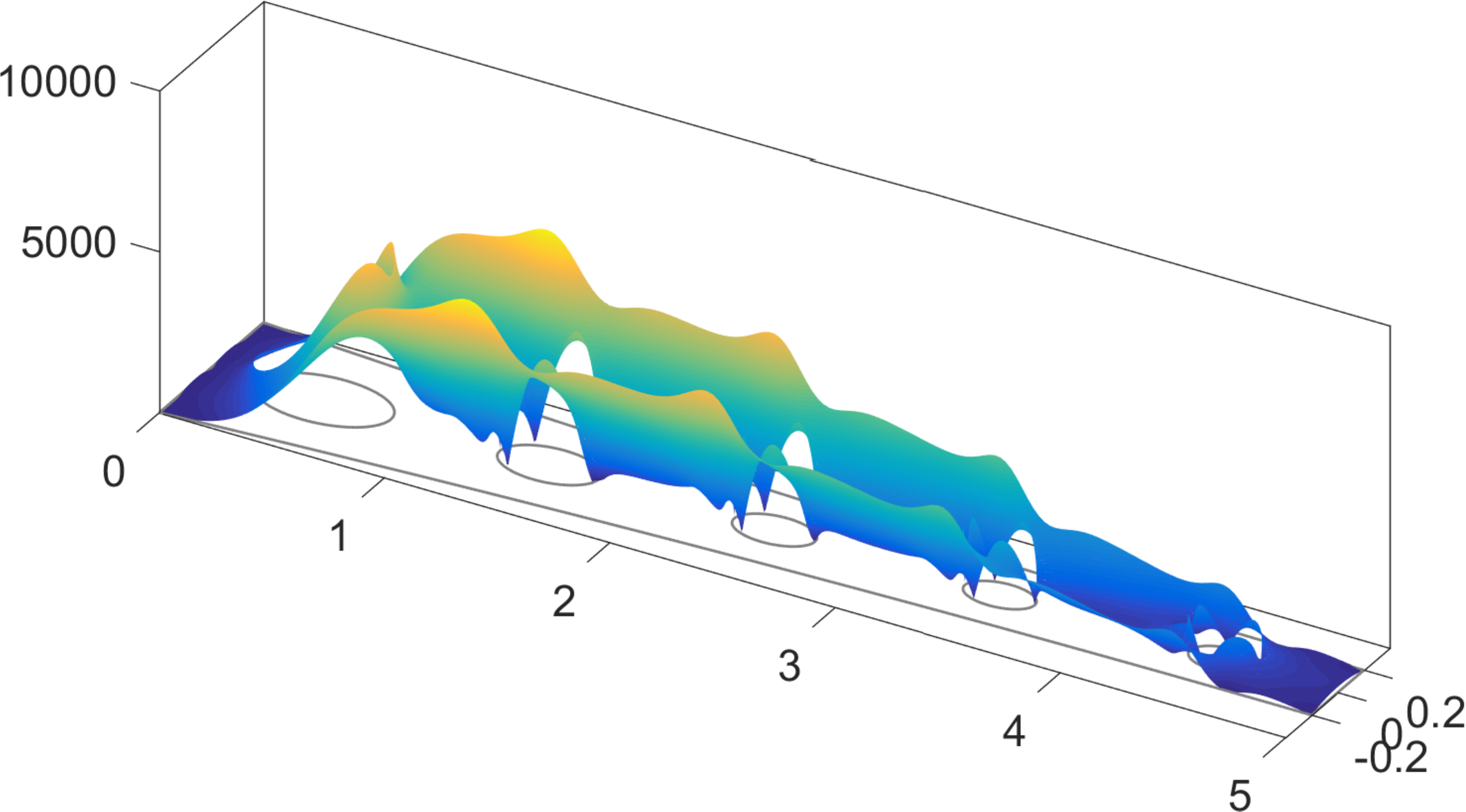}}

\caption{Deformed configurations scaled by a factor of $200$ and von Mises
stress for the two different support cases, (a) and (b) refer to a
fixed support on the left side (support case 1), (c) and (d) to a
fixed support in the left ellipsoid (case 2).}

\label{fig:TestCase2dC_Sketch2} 
\end{figure}

There are no analytical solutions available for the two different
support scenarios of this test case, however, the stored energy has
been computed by an overkill solution using standard $p$-FEM. For
support case 1, the elastic energy is $\mathfrak{e}=0.03246547385\pm10^{-8}\mathrm{kNm}$
and for support case 2, $\mathfrak{e}=0.02361112384\pm10^{-10}\mathrm{kNm}$.
The different uncertainties reflect the fact, that the singularities
hinder an optimal convergence of the $p$-FEM for support case 1 when
computing the overkill solution.

\begin{figure}
\centering

\subfigure[approx.~error, case 1]{\includegraphics[height=4cm]{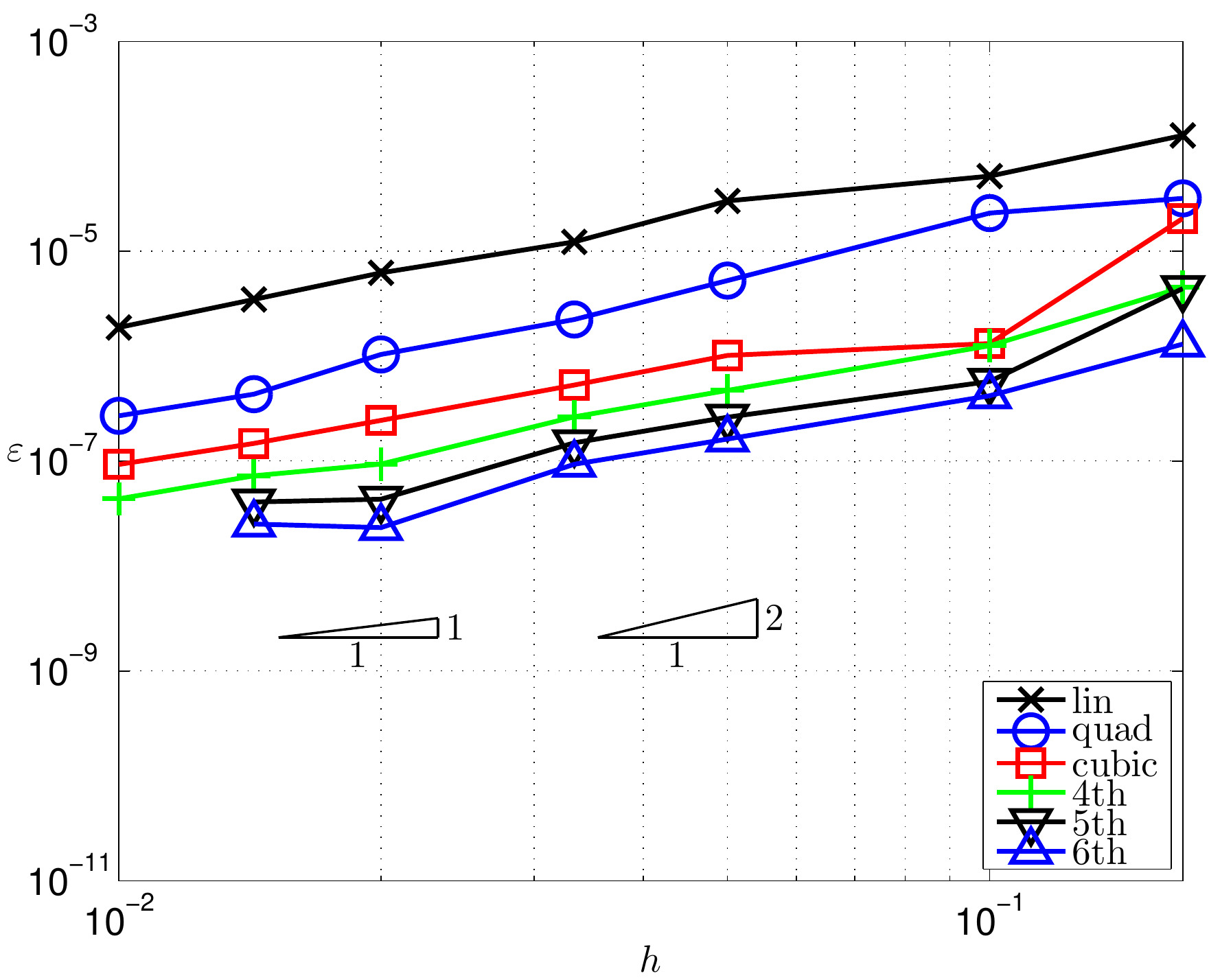}}\quad\subfigure[approx.~error, case 2]{\includegraphics[height=4cm]{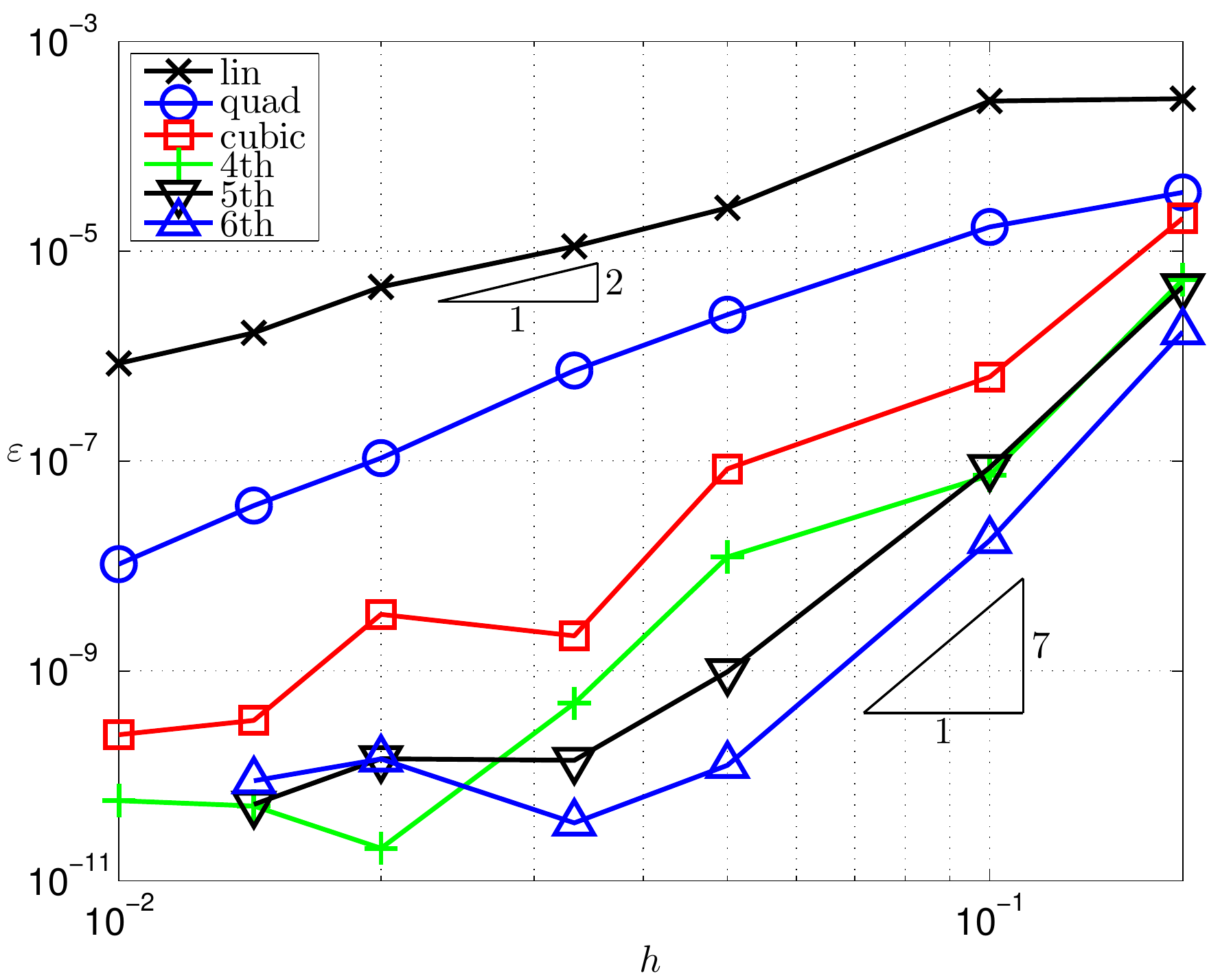}}\quad\subfigure[condition number]{\includegraphics[height=4cm]{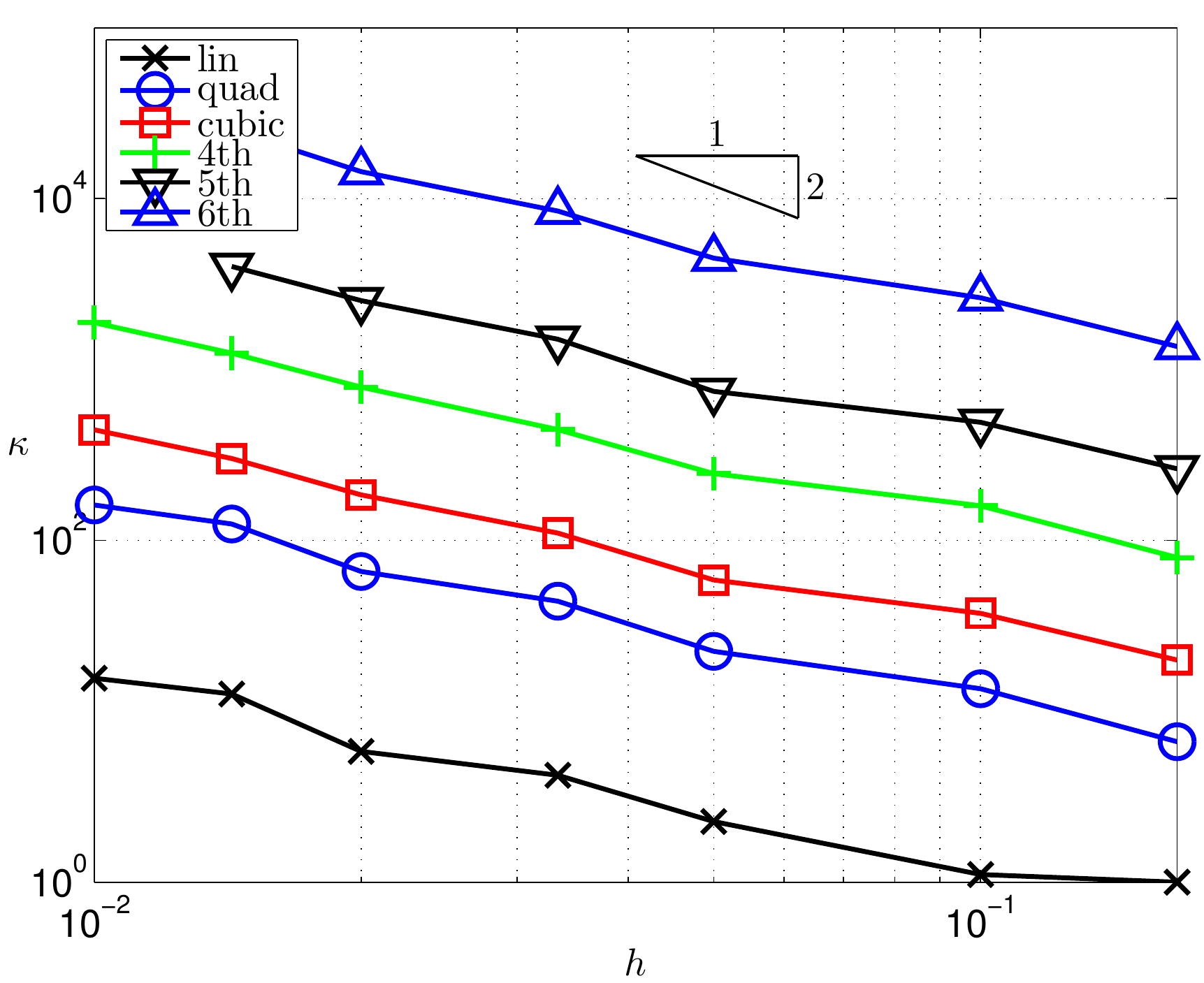}}

\caption{The approximation error for (a) support case 1 and (b) case 2, (c)
shows the corresponding condition numbers of the system matrices which
are almost identical for the two cases.}

\label{fig:TestCase2dC_Res} 
\end{figure}

Convergence results are shown in Fig.~\ref{fig:TestCase2dC_Res}.
For the first support case where singularities are present, it is
clearly seen from Fig.~\ref{fig:TestCase2dC_Res}(a) that only first
order convergence rates in the energy are achieved. It is nevertheless
noted that from linear to cubic elements the results improve by about
one order of magnitude. When comparing the results for even higher
orders the improvements is more than two orders of magnitude. That
is, although no higher-order convergence rates are achieved, there
is still a significant improvement of the results. The energy error
in Fig.~\ref{fig:TestCase2dC_Res}(b) refers to the second support
case and clearly converges with optimal rates as no singularities
are present. The condition numbers for the different resolutions of
the background meshes and various element orders are seen in Fig.~\ref{fig:TestCase2dC_Res}(c).
They are almost identical for the two support cases and behave as
expected.

\begin{figure}
\centering

\subfigure[approx.~error]{\includegraphics[height=4cm]{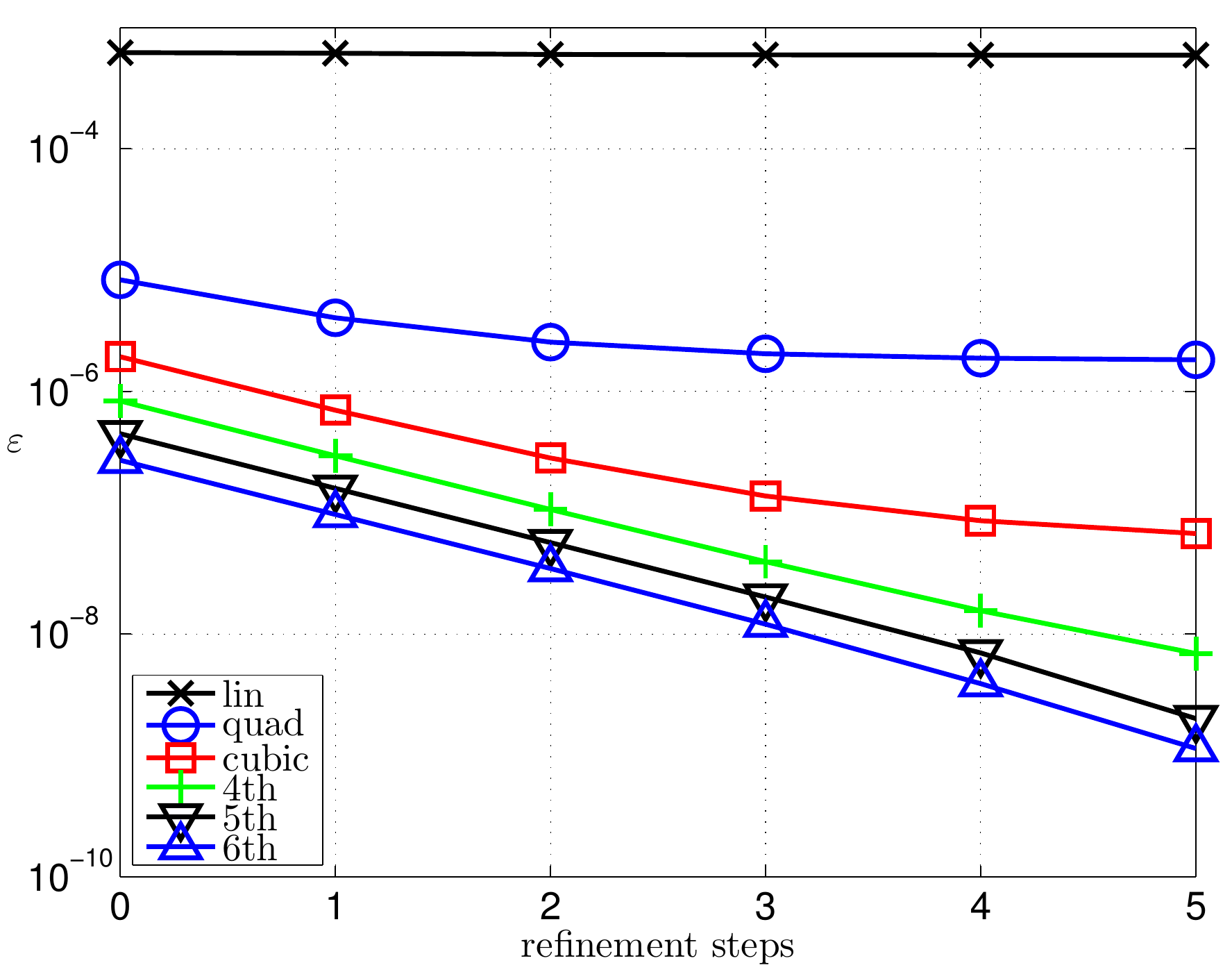}}\qquad\subfigure[condition number]{\includegraphics[height=4cm]{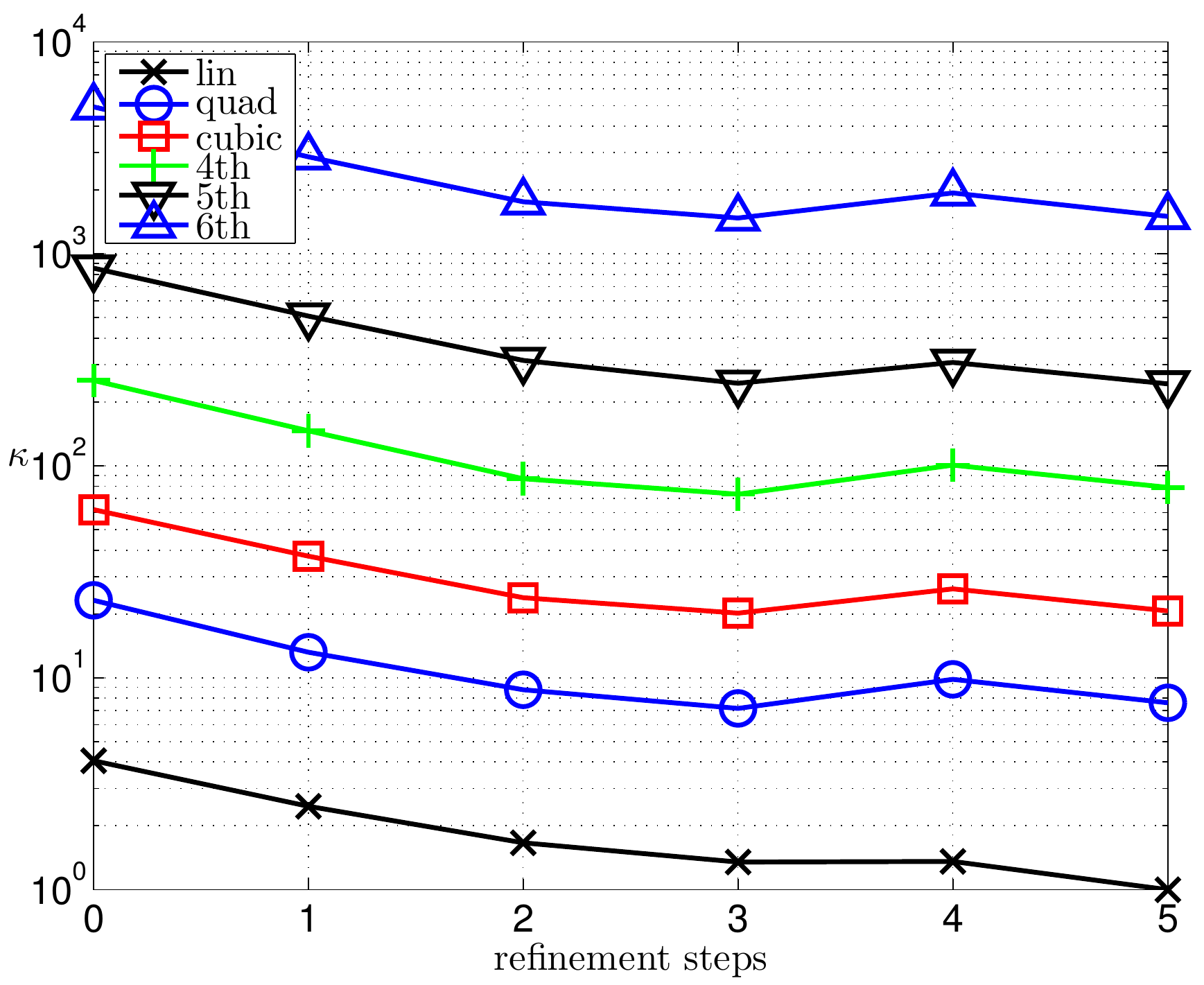}}

\caption{(a) The approximation error for a different number of adaptive refinements
at the singularities for support case 1 and (b) the corresponding
condition numbers of the system matrices.}

\label{fig:TestCase2dC_ResAdapt} 
\end{figure}

The focus is now only on support case 1 with the singularities at
the upper and lower left corner of the beam. The results in Fig.~\ref{fig:TestCase2dC_Res}
have been achieved \emph{without} adaptive refinements at these singularities.
Next, it is investigated how adaptive refinements improve the approximation
error. The following results, presented in Fig.~\ref{fig:TestCase2dC_ResAdapt},
are achieved for a background mesh with a \emph{fixed} resolution
but varying orders of the elements and a different number of refinement
steps at the corners (up to $5$). The corresponding unrefined mesh
taken as the starting mesh for each computation is seen in \ref{fig:TestCase2dC_Sketch}(c).
Results are visualized in Fig.~\ref{fig:TestCase2dC_ResAdapt} where
the horizontal axis shows the number of refinement steps. It is seen
that for \emph{linear} meshes, the adaptive refinements at the singularities
does not change the results noticeably because the error in the bulk
mesh is still dominant. However, for increasing orders of the elements,
the local refinements improve the results to a great extent. So it
is obvious that the singularities hinder optimal convergence rates
for the higher-order elements as expected and that adaptivity is highly
useful. For the $6th$-order elements, the error is improved by $3$
orders of magnitude after $5$ refinements steps at the singularities!

\subsection{Spanner}

\begin{figure}
\centering

\subfigure[background mesh]{\includegraphics[width=0.3\textwidth]{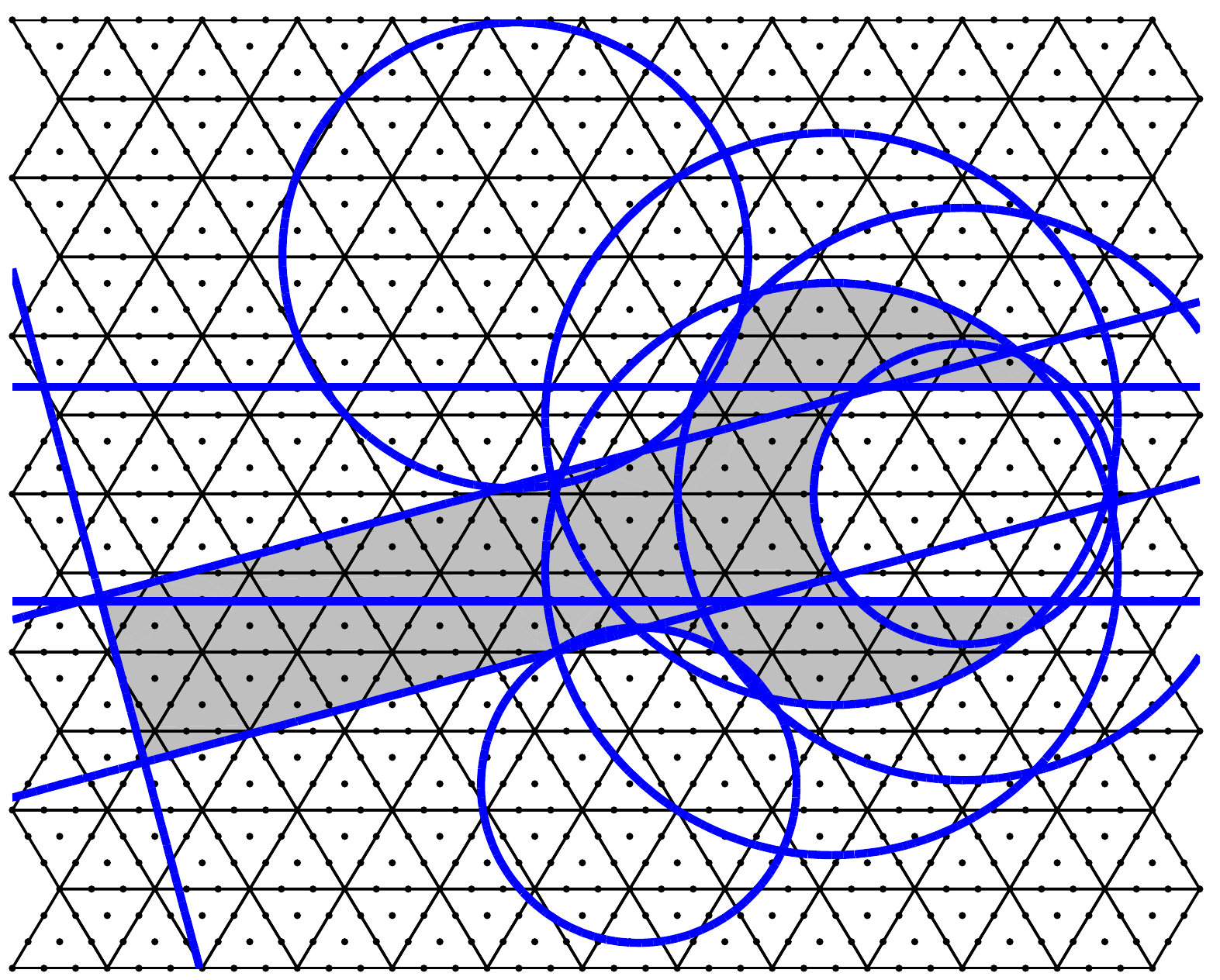}}\quad\subfigure[exact sol., 2D view]{\includegraphics[width=0.3\textwidth]{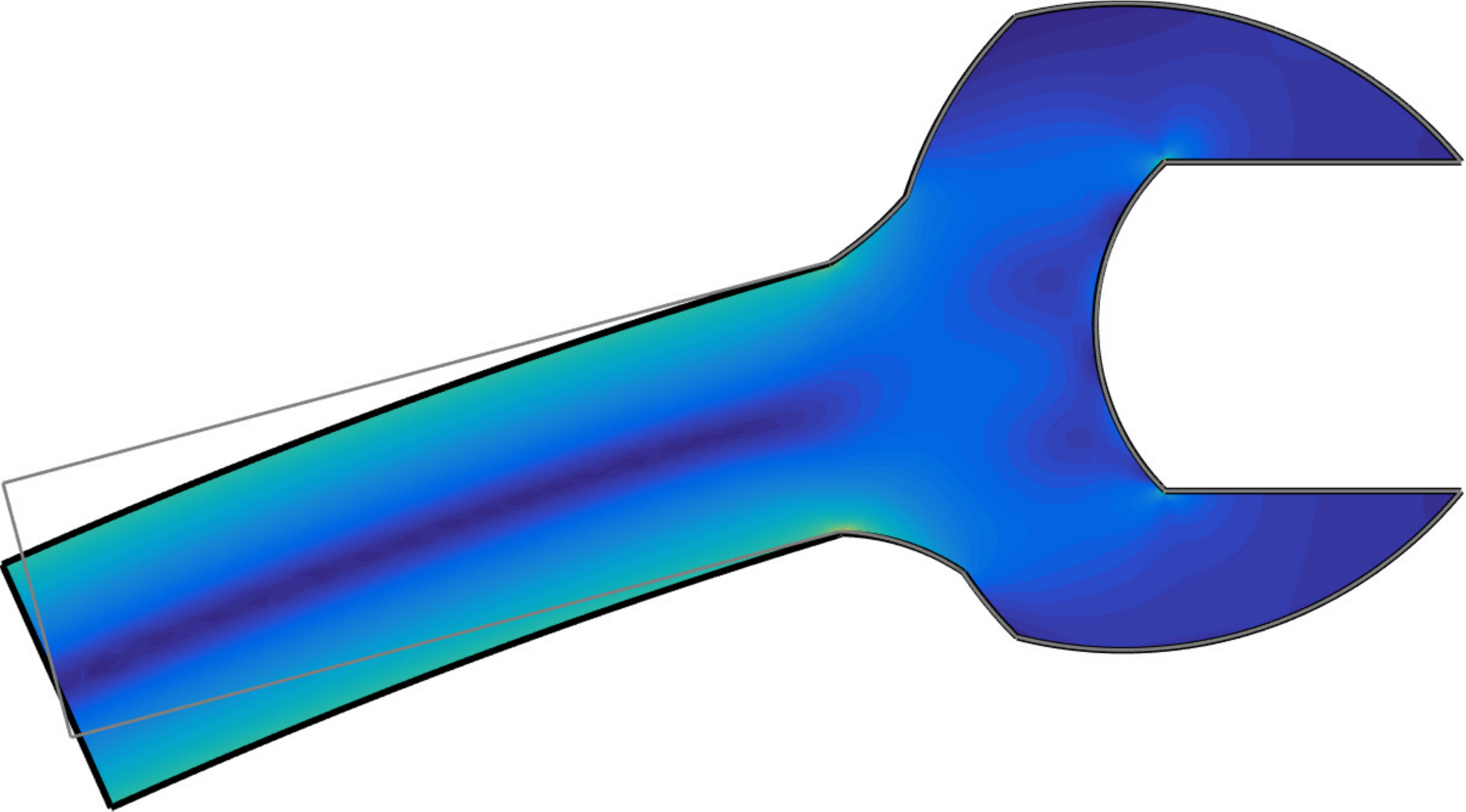}}\quad\subfigure[exact sol., 3D view]{\includegraphics[width=0.3\textwidth]{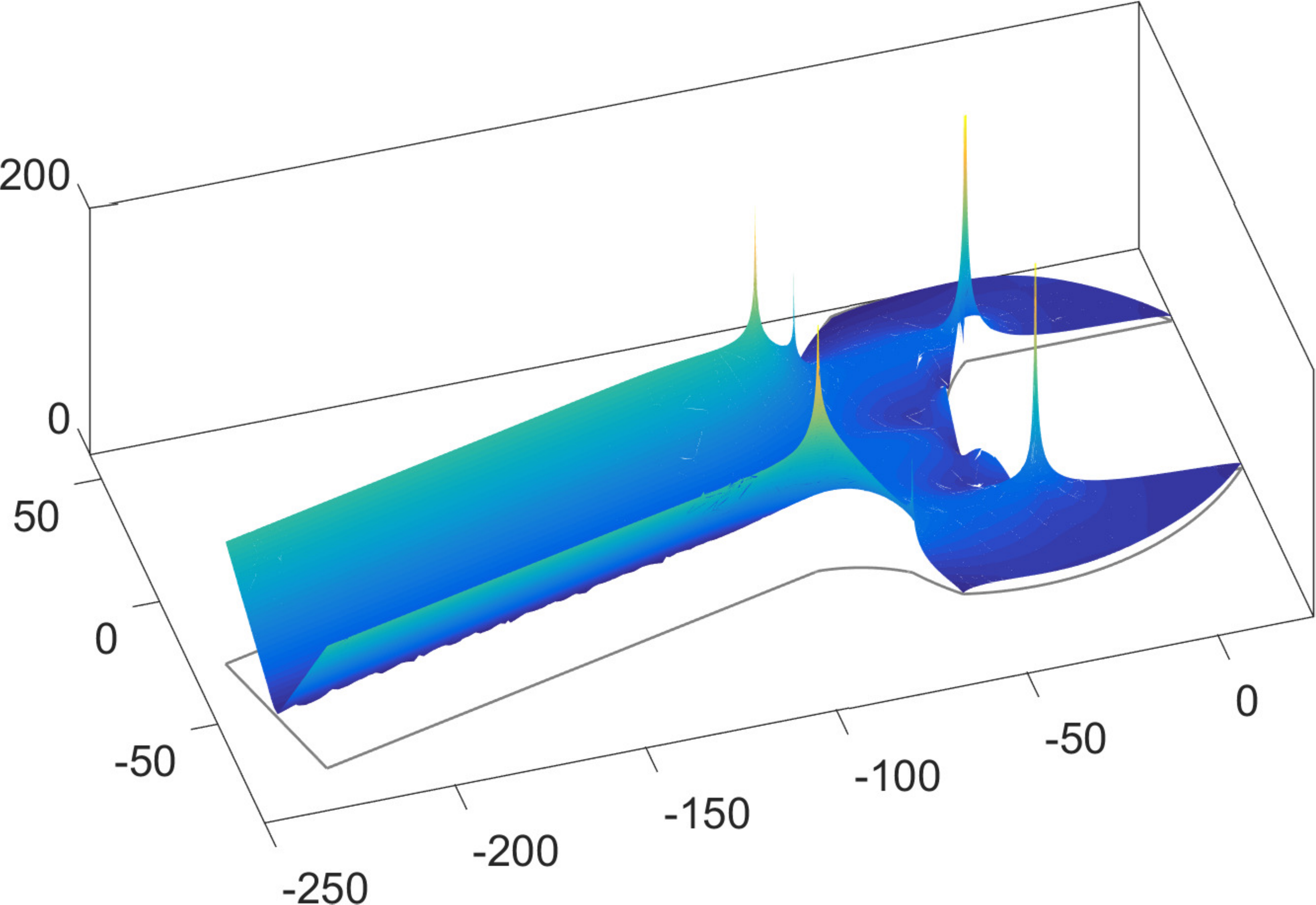}}

\subfigure[generated mesh, coarse]{\includegraphics[width=0.3\textwidth]{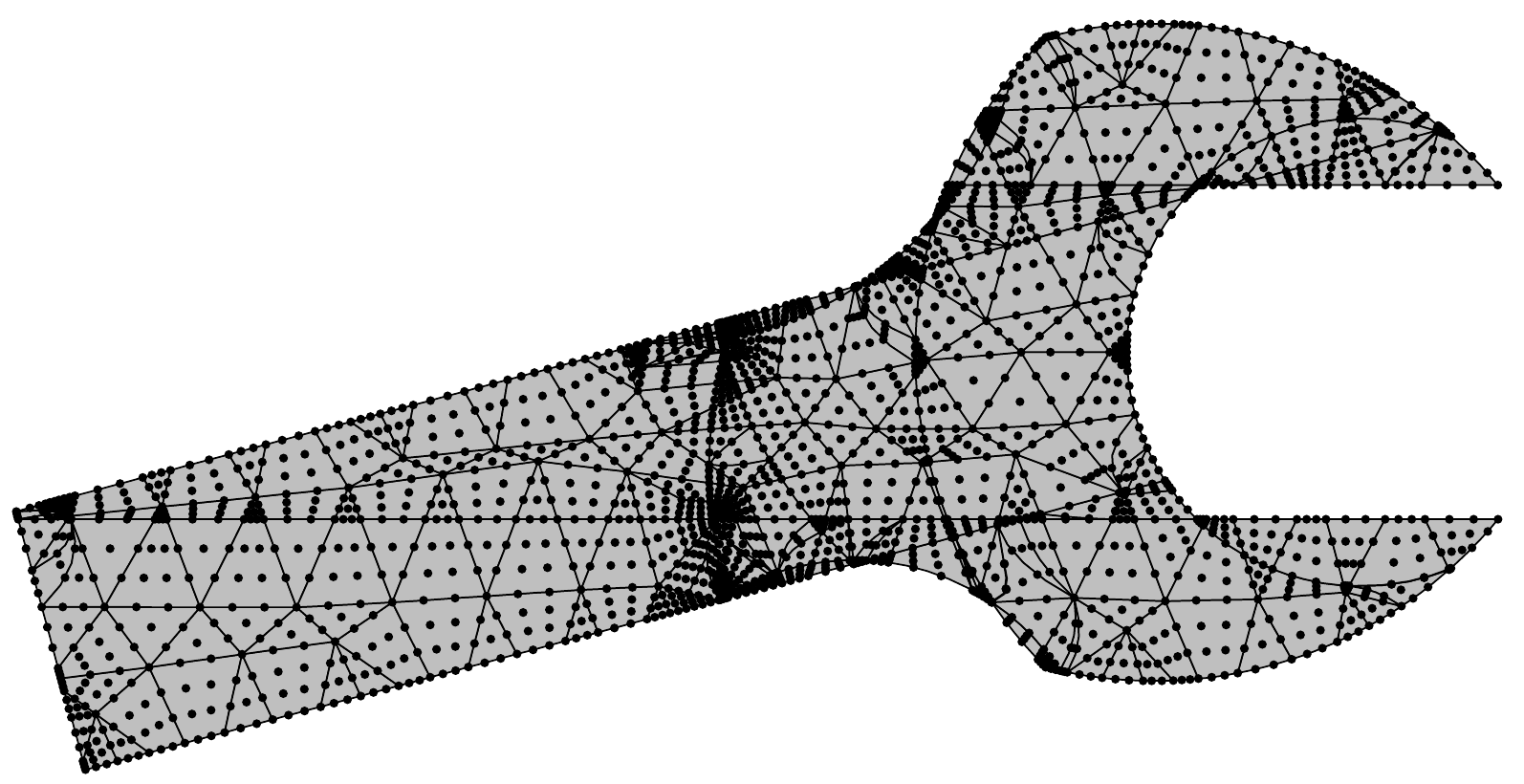}}\quad\subfigure[generated mesh, finer]{\includegraphics[width=0.3\textwidth]{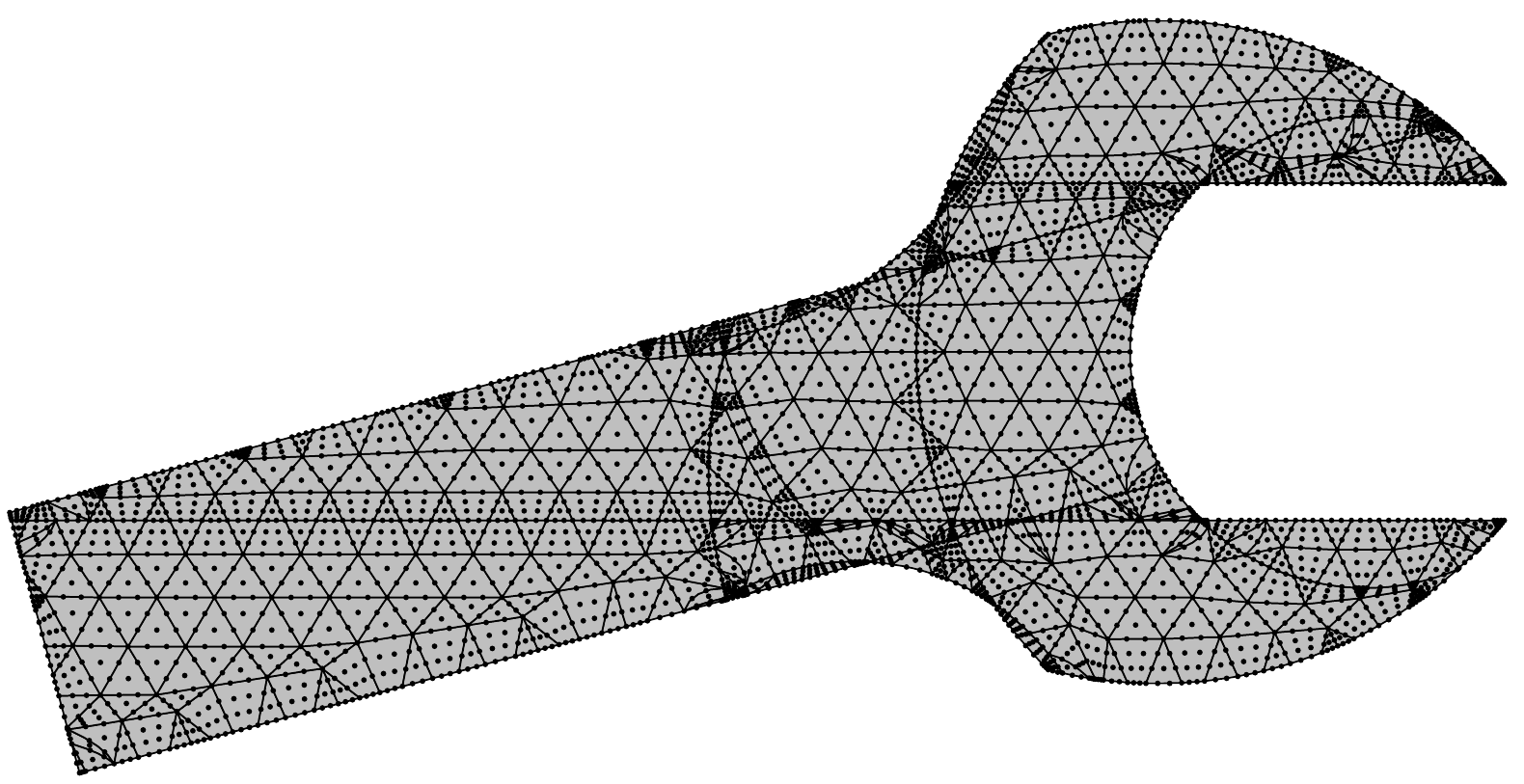}}\quad\subfigure[generated mesh, adapt.]{\includegraphics[width=0.3\textwidth]{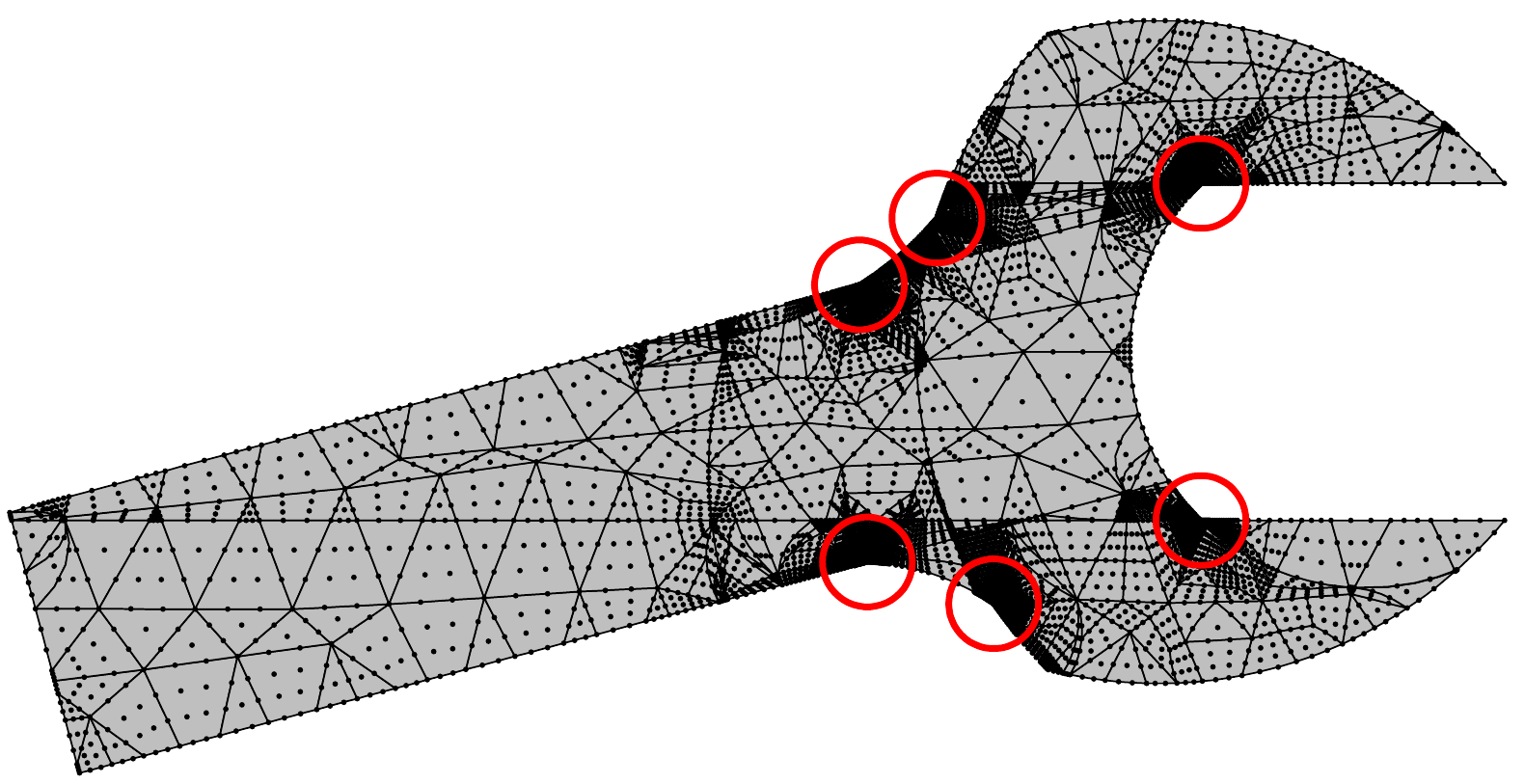}}

\caption{(a) Universal background mesh and zero-level sets of the $11$ level-set
functions, (b) the deformed configuration, (c) von Mises stress in
3D view showing the singularities, (d) to (f) show examples of generated
conforming higher-order meshes, (f) is adaptively refined at the singularities.}

\label{fig:TestCase2dD_Sketch} 
\end{figure}

Next, the geometry shown in Fig.~\ref{fig:TestCase2dD_Sketch} is
considered and refers to a spanner being very similar to normed spanner
geometries defined in DIN 895. The geometry is embedded into a universal
mesh composed by triangular elements of different orders. $11$ level-set
functions are employed to define the geometry. There are $5$ straight
lines defined by linear level-set functions
\begin{eqnarray*}
\phi_{i}\left(\vek x\right) & = & \vek n_{i}\cdot\left(\vek x-\vek x_{i}^{\star}\right),\\
 & = & n_{x,i}\cdot\left(x-x_{i}^{\star}\right)+n_{y,i}\cdot\left(y-y_{i}^{\star}\right),
\end{eqnarray*}
which imply zero-level sets going through points $\vek x_{i}^{\star}$
with normal vectors $\vek n_{i}$. See Table \ref{tab:SpannerSraightLines}
for the concrete values of $\vek x_{i}^{\star}$ and $\vek n_{i}$.
All measurements for this test case are given in $\mathrm{mm}$. Furthermore,
there are a number of arc segments in the boundary definition of the
spanner. Therefore, we need $6$ level-set functions,
\begin{eqnarray*}
\phi_{i}\left(\vek x\right) & = & \left|\left|\vek x-\vek x_{i}^{\circ}\right|\right|-R_{i}^{\circ},\\
 & = & \sqrt{\left(x-x_{i}^{\circ}\right)^{2}+\left(y-y_{i}^{\circ}\right)^{2}}-R_{i}^{\circ},
\end{eqnarray*}
implying circular zero-level sets centered at $\vek x_{i}^{\circ}$
with radius $R_{i}^{\circ}$. See Table \ref{tab:SpannerArcs} for
the specific values of $\vek x_{i}^{\circ}$ and $R_{i}^{\circ}$.

\begin{table}
\centering%
\begin{tabular}{|l|r|r|r|r|}
\hline 
 & $x_{i}^{\star}$ & $y_{i}^{\star}$ & $n_{x,i}$ & $n_{y,i}$\tabularnewline
\hline 
\hline 
mouth, top & $0.000000$ & $30.000000$ & $0.000000$ & $1.000000$\tabularnewline
\hline 
mouth, bottom & $0.000000$ & $-30.000000$ & $0.000000$ & $1.000000$\tabularnewline
\hline 
handle, up & $0.000000$ & $36.100494$ & $-0.258819$ & $0.965926$\tabularnewline
\hline 
handle, bottom & $0.000000$ & $-13.592762$ & $-0.258819$ & $0.965926$\tabularnewline
\hline 
handle, end & $0.000000$ & $-929.486795$ & $0.965926$ & $0.258819$\tabularnewline
\hline 
\end{tabular}

\caption{\label{tab:SpannerSraightLines}Definition of the straight segments
of the spanner, measurements in $\mathrm{mm}$.}

\end{table}

\begin{table}
\centering%
\begin{tabular}{|l|r|r|r|}
\hline 
 & $x_{i}^{\circ}$ & $y_{i}^{\circ}$ & $R_{i}^{\circ}$\tabularnewline
\hline 
\hline 
circle 1 & $0.000000$ & $0.000000$ & $42.000000$\tabularnewline
\hline 
circle 2 & $0.000000$ & $0.000000$ & $80.000000$\tabularnewline
\hline 
circle 3 & $-125.361456$ & $66.697117$ & $65.100000$\tabularnewline
\hline 
circle 4 & $-90.856482$ & $-81.419283$ & $44.100000$\tabularnewline
\hline 
circle 5 & $-37.000000$ & $-21.000000$ & $80.000000$\tabularnewline
\hline 
circle 6 & $-37.000000$ & $21.000000$ & $80.000000$\tabularnewline
\hline 
\end{tabular}

\caption{\label{tab:SpannerArcs}Definition of the arc segments of the spanner,
measurements in $\mathrm{mm}$.}
\end{table}

The material is again composed by steel with $E=2.1\cdot10^{5}\,\unitfrac{N}{mm^{2}}$
and $\nu=0.3$. The beam is loaded by gravity acting as a body force
of $f_{y}=-78.5\cdot10^{-6}\,\unitfrac{N}{mm^{3}}$ and a traction
at the end of the handle. This traction acts in parallel direction
of the handle and is distributed linearly between $-100\unitfrac{N}{mm}$
at the bottom side and $+100\unitfrac{N}{mm}$ on the top side. It
loads the handle of the spanner with a resulting bending moment of
$M_{z}=38400\,\mathrm{Nmm}$. All nodes on the two straight, parallel
sides of the mouth are fixed. The deformed configuration is shown
in Fig.~\ref{fig:TestCase2dD_Sketch}(b) scaled by a factor of $50$.
The geometry of the spanner features $6$ reentrant corners, two at
the mouth, two on the top side between the handle and the front part
and two on the opposite side, see the red circles in Fig.~\ref{fig:TestCase2dD_Sketch}(f).
It is thus clear that singular stresses have to be expected.

\begin{figure}
\centering

\subfigure[approx.~error, no adapt.]{\includegraphics[height=4cm]{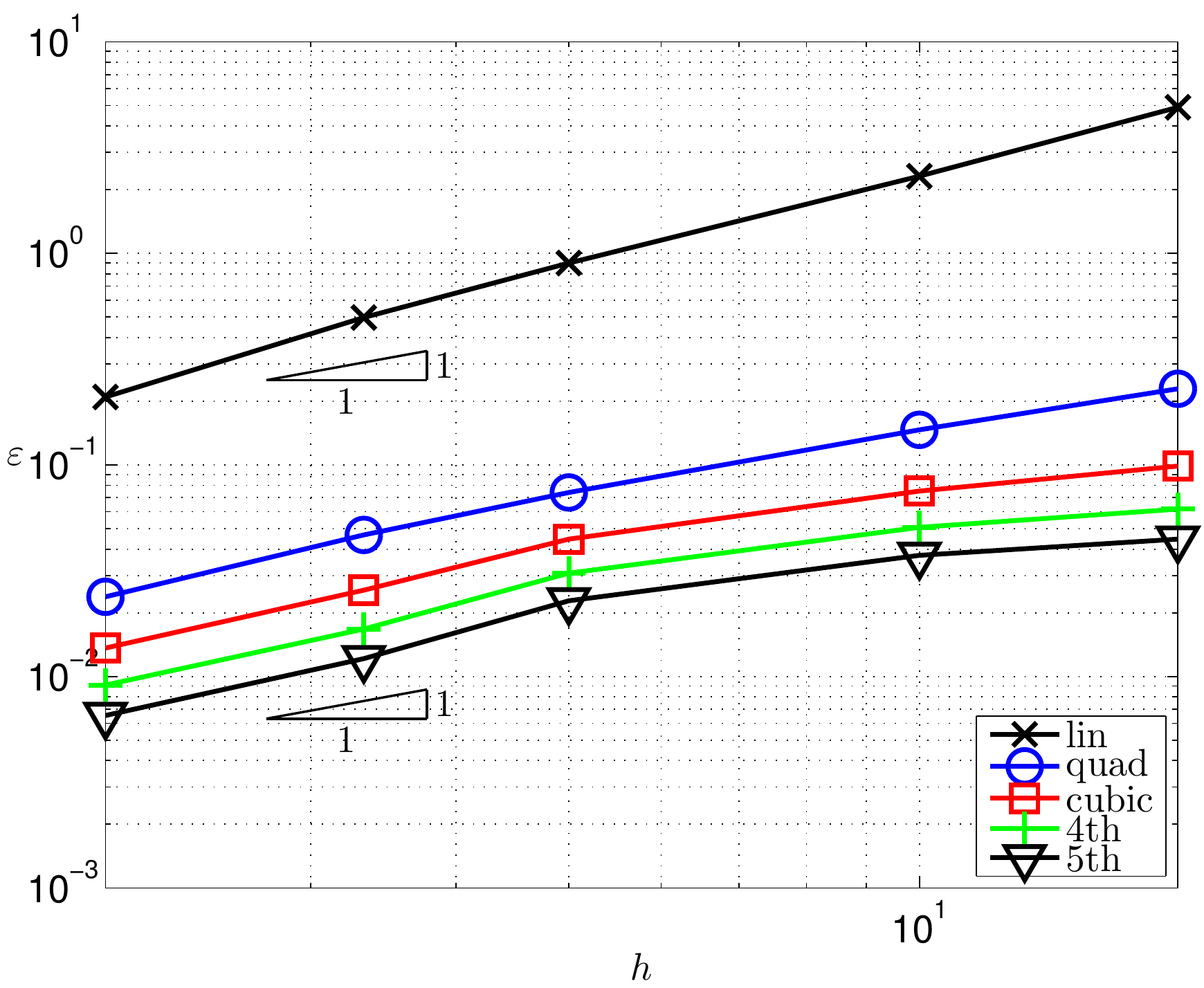}}\quad\subfigure[approx.~error, adapt.]{\includegraphics[height=4cm]{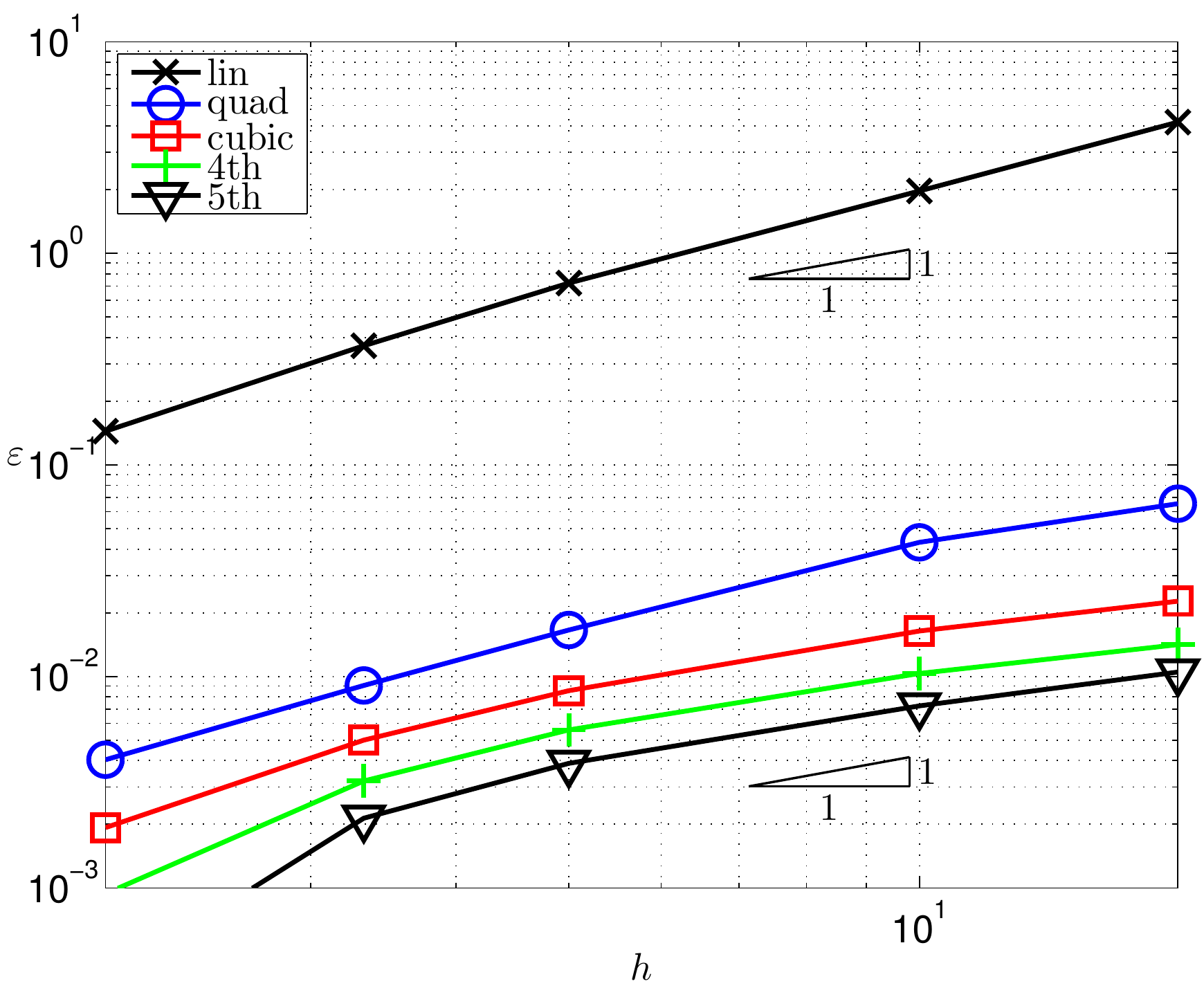}}\subfigure[condition number]{\includegraphics[height=4cm]{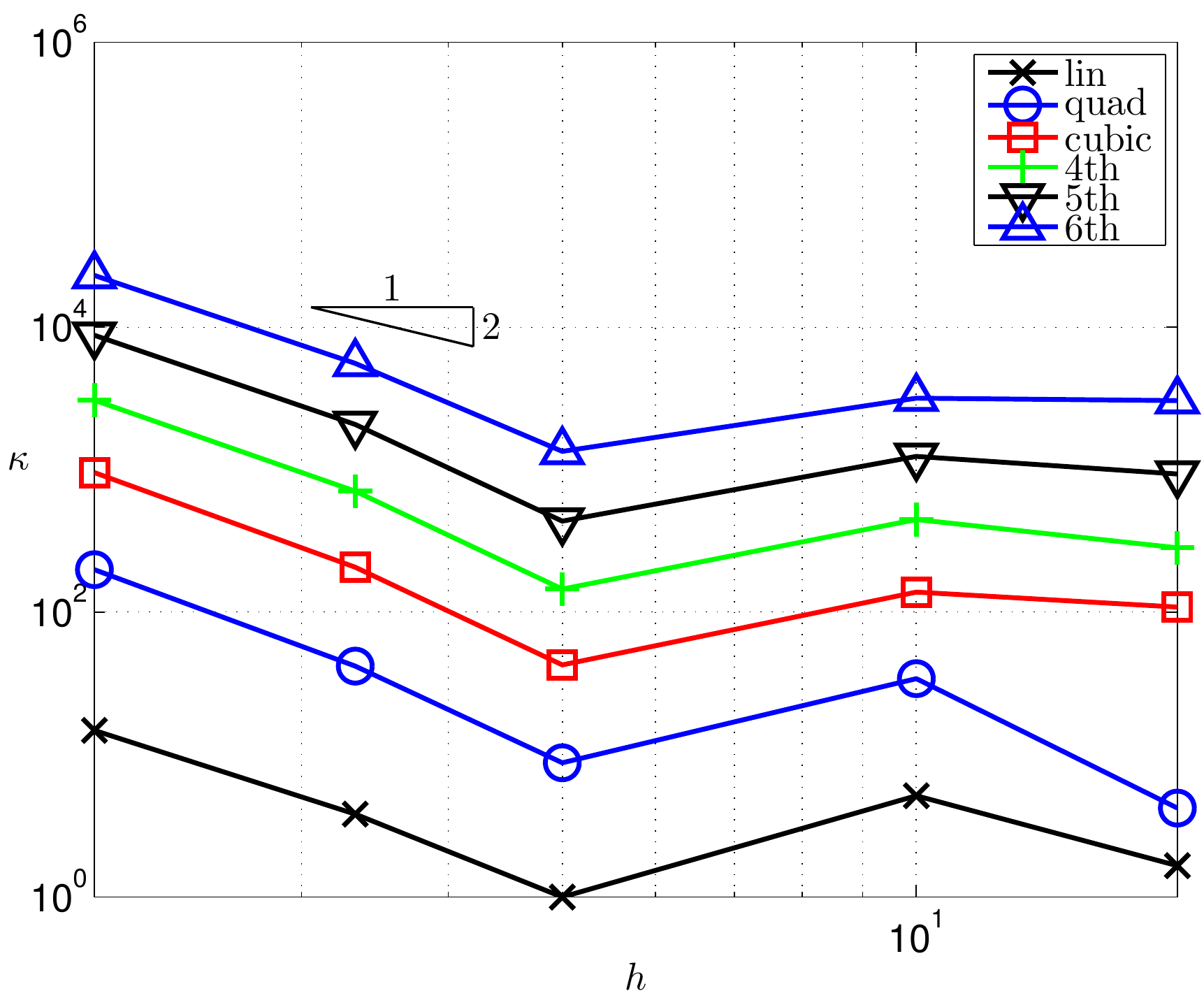}}

\caption{The approximation error (a) without adaptivity and (b) with adaptivity,
(c) shows the corresponding condition numbers of the system matrices
which are quite similar for the two cases.}

\label{fig:TestCase2dD_Res} 
\end{figure}

Convergence results are shown in Fig.~\ref{fig:TestCase2dD_Res}.
Again, there is no analytical solution available wherefore the stored
energy is used for the convergence study. An overkill solution yields
$\mathfrak{e}=61.49248\pm10^{-4}\mathrm{Nmm}$. Fig.~\ref{fig:TestCase2dD_Res}(a)
displays the convergence for automatically generated meshes \emph{without}
adaptive refinements at the singularities. Of course, only first order
convergence rates are achieved due to the singularities, however,
the error level is drastically improved for the higher-order meshes.
Fig.~\ref{fig:TestCase2dD_Res}(b) shows results for adaptively refined
meshes as seen in Fig.~\ref{fig:TestCase2dD_Sketch}(f). Only three
refinement steps improve the error by about one order of magnitude
for the higher-order elements. The only exception are linear elements
where the error in the bulk mesh is still too large to improve much
through a local refinement at the singularities only. Condition numbers
are seen in Fig.~\ref{fig:TestCase2dD_Res}(c) for the unrefined
meshes and look quite similar for the refined case as well.

\subsection{Cube with spherical hole\label{XX_CubeWithHole}}

\begin{figure}
\centering

\subfigure[domain with hole]{\includegraphics[height=4cm]{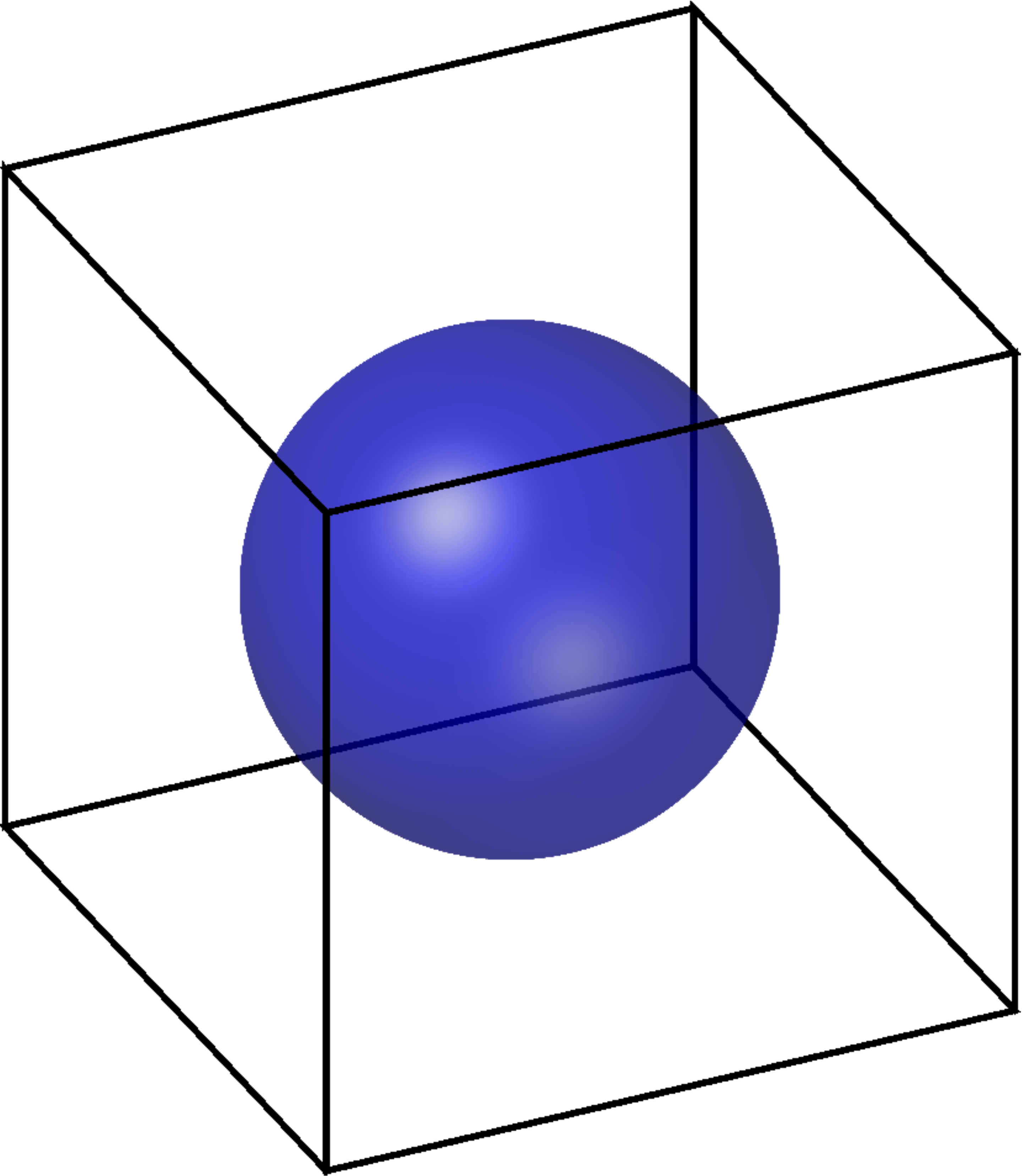}}\qquad\subfigure[background mesh]{\includegraphics[height=4cm]{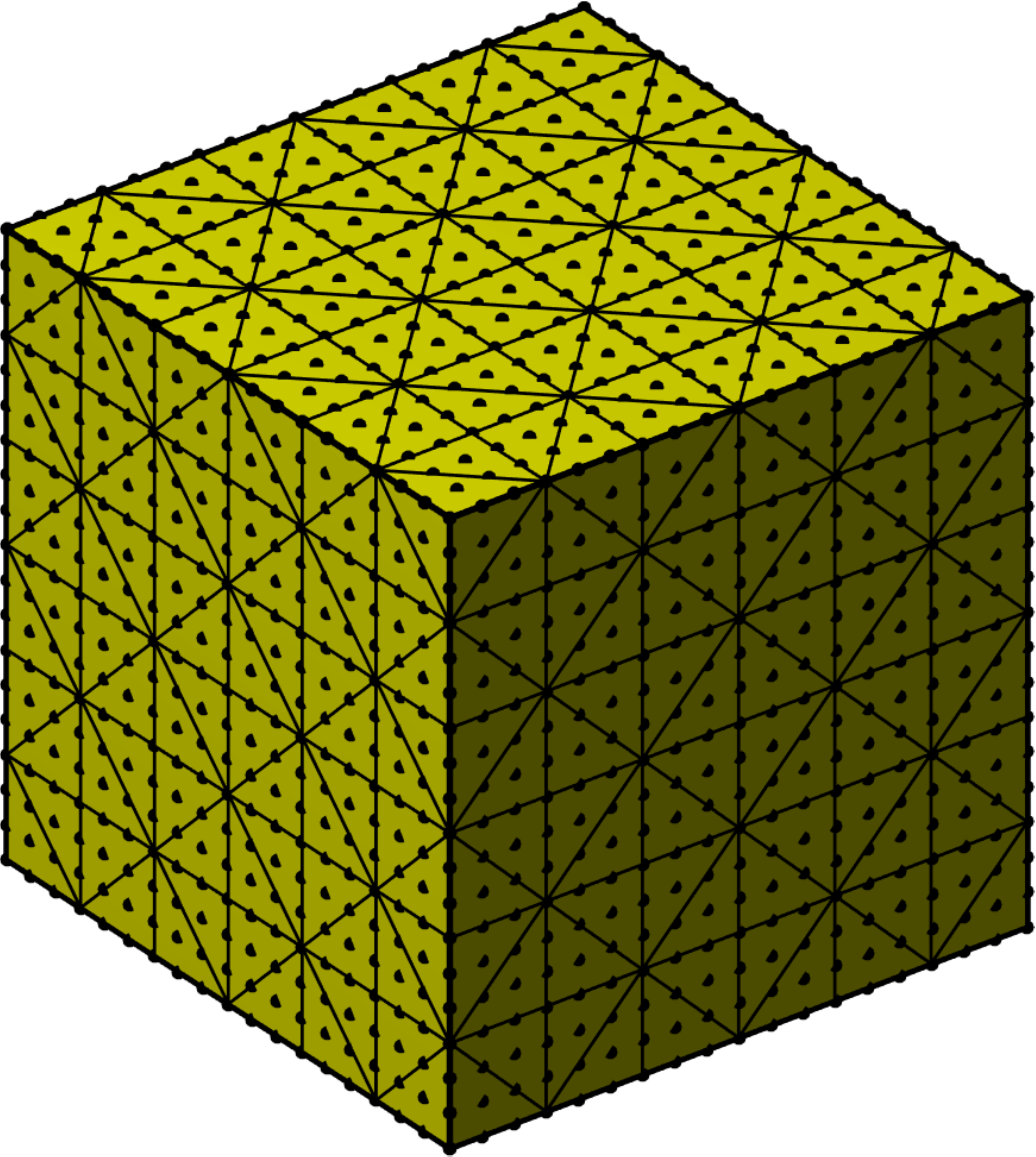}}

\subfigure[generated mesh]{\includegraphics[height=4cm]{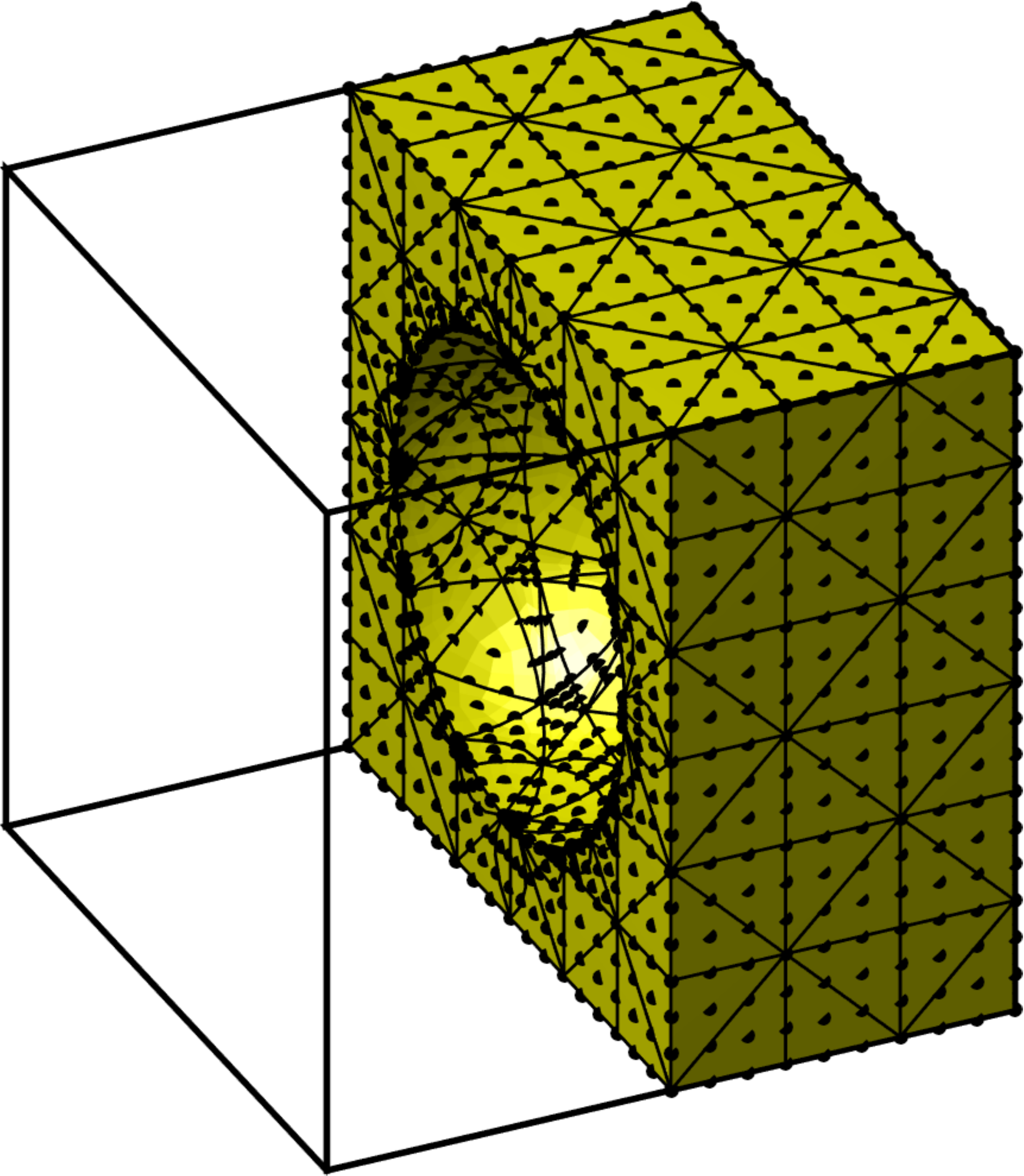}}\qquad\subfigure[exact solution]{\includegraphics[height=4cm]{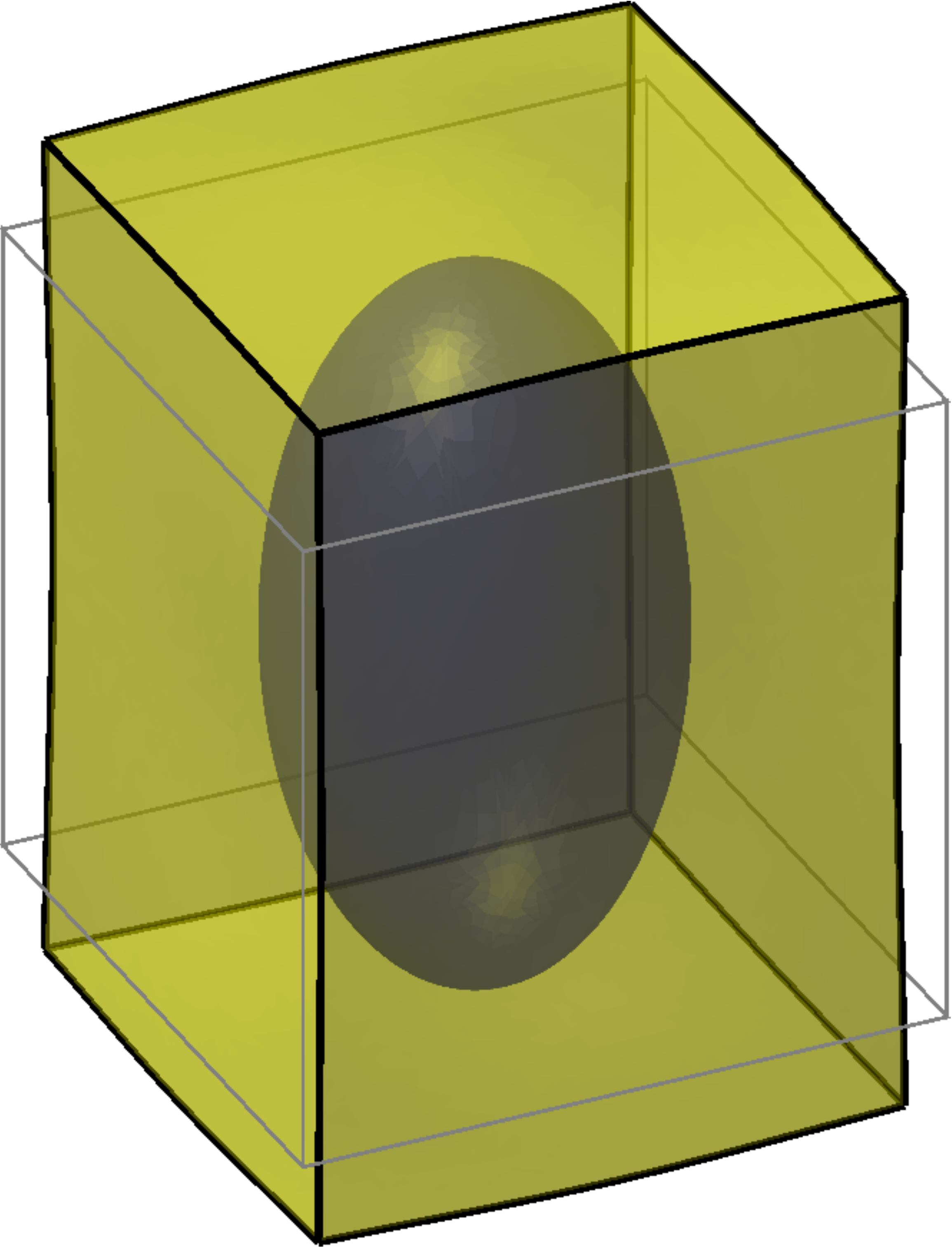}}

\caption{(a) Cube with spherical hole (b) background mesh, (c) part of the
generated higher-order mesh with hole, (d) deformed configuration.}

\label{fig:TestCase3dA_Sketch} 
\end{figure}

The next test cases feature three-dimensional geometries. The first
case is the extension to three dimensions of the square shell with
circular hole from Section \ref{XX_SquareWithHole} and the same material
properties are used here. The domain is $\left[-1,1\right]^{3}$ with
a spherical hole of radius $R=0.7123$, see Fig.~\ref{fig:TestCase3dA_Sketch}(a)
for a sketch of the situation. An example background mesh is seen
in Fig.~\ref{fig:TestCase3dA_Sketch}(b) and a resulting conforming
mesh with hole in (c). The deformed configuration according to the
exact solution from \cite{Goodier_1933a}, see the Appendix \ref{XX_ExactSolCubeWithHole},
is displayed in Fig.~\ref{fig:TestCase3dA_Sketch}(d).

\begin{figure}
\centering

\subfigure[approximation error]{\includegraphics[height=4cm]{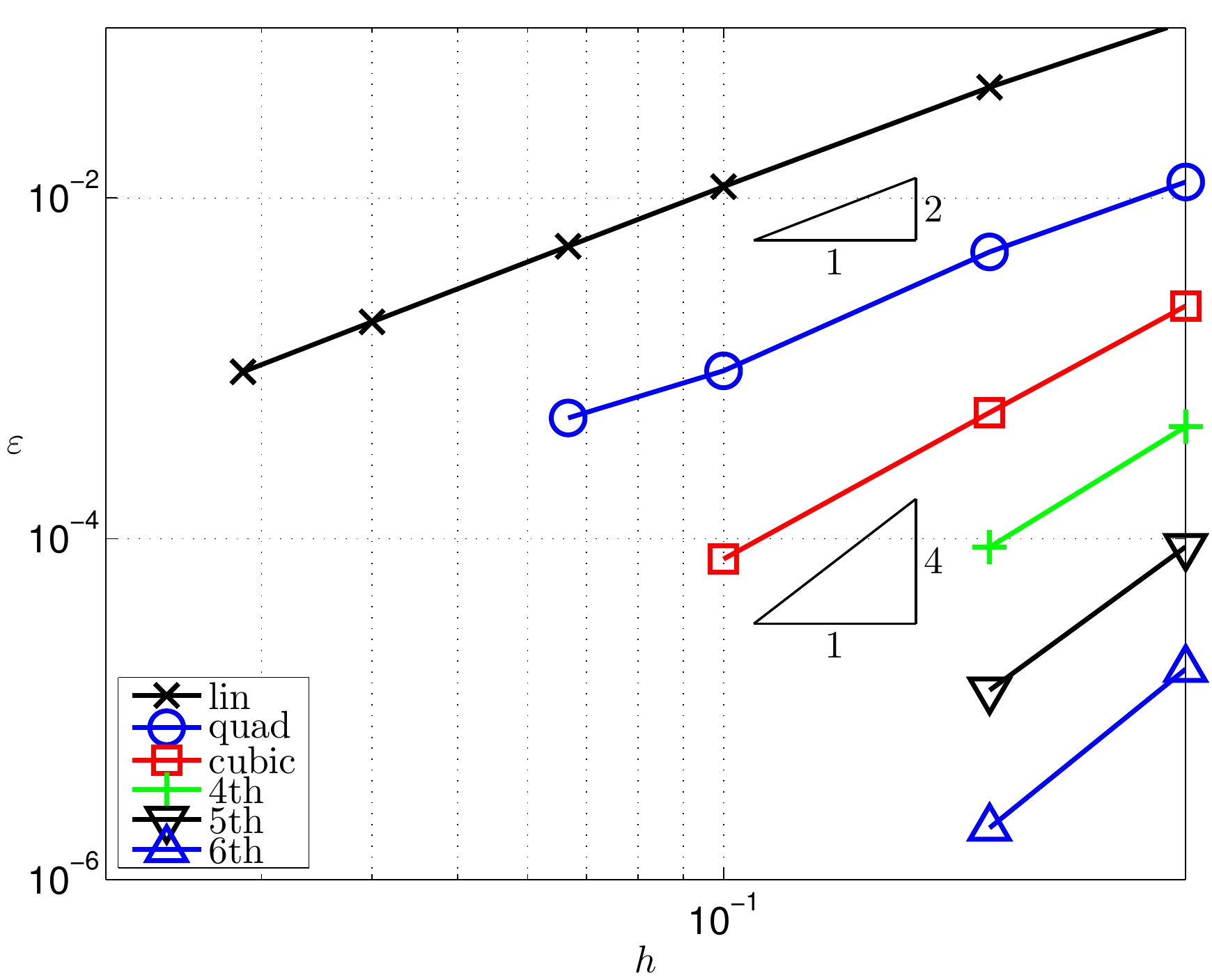}}\qquad\subfigure[condition number]{\includegraphics[height=4cm]{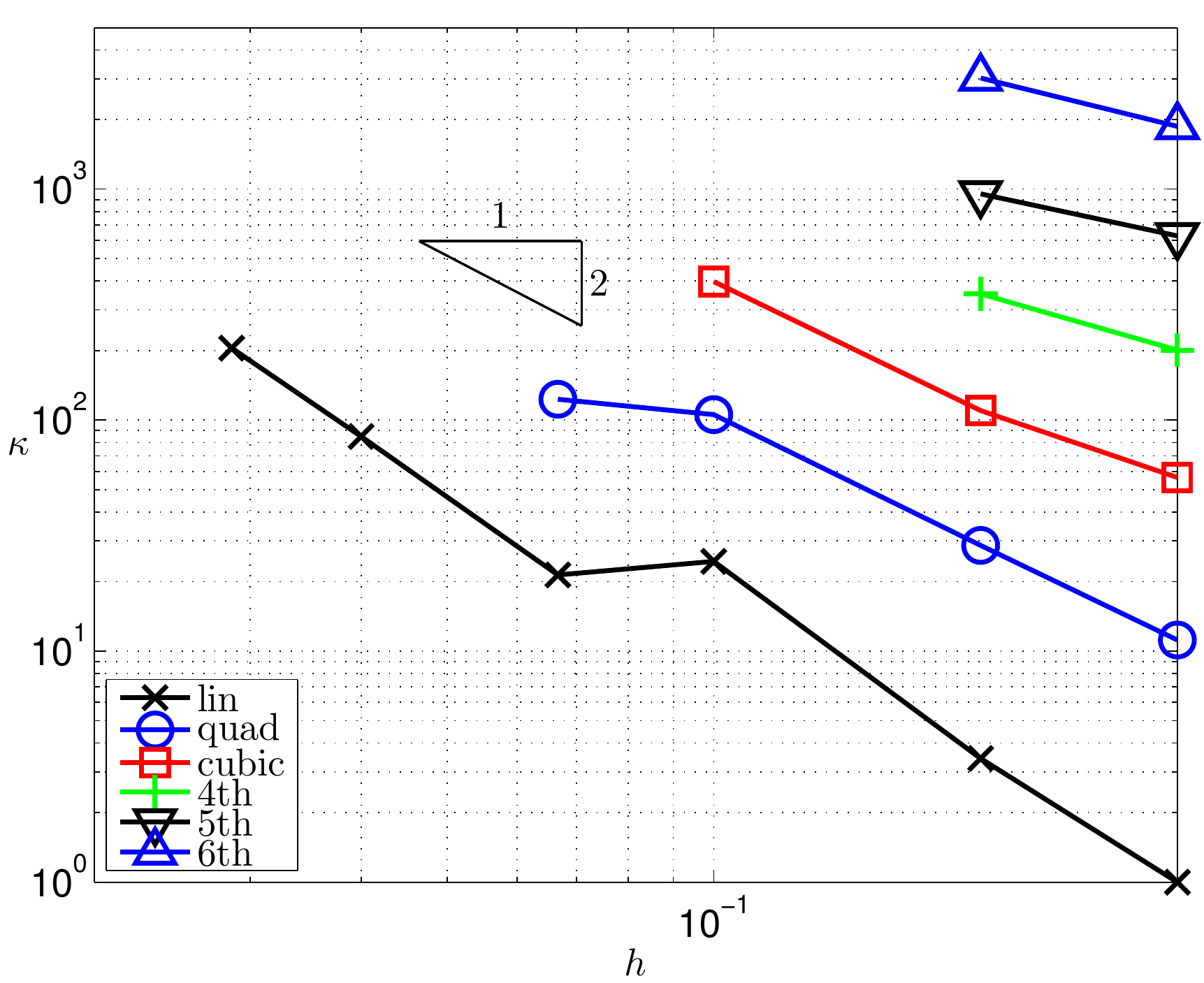}}

\caption{(a) The approximation error for the cube mesh with spherical hole,
(b) the corresponding condition numbers of the system matrices.}

\label{fig:TestCase3dA_Res} 
\end{figure}

For the convergence study, the number of elements per dimensions,
$n_{d}$, of the background mesh is systematically increased and $n_{d}=\left\{ 6,10,20,30,50,70\right\} $
elements are used with varying orders between $1$ and $6$. Only
meshes with less than $500.000$ nodes are considered, leading to
$1.5\cdot10^{6}$ degrees of freedom as three displacement components
at each node are present. Convergence results of the aproximation
error are displayed in Fig.~\ref{fig:TestCase3dA_Res}(a) and are
optimal as expected. Note that in \cite{Fries_2015a,Fries_2016a},
integration errors for this test case are investigated and in \cite{Fries_2016b},
interpolation errors. Nevertheless, these are the first higher-order
convergence results of \emph{approximation }errors achieved with the
CDFEM in three dimensions reported so far. Fig.~\ref{fig:TestCase3dA_Res}(b)
shows the corresponding condition numbers which are well-bounded and
prove the success of the node manipulations discussed in Section \ref{X_MeshGeneration}.

\subsection{Cube with spherical inclusion\label{XX_CubeWithSphereInclusion}}

\begin{figure}
\centering

\subfigure[generated mesh]{\includegraphics[height=4cm]{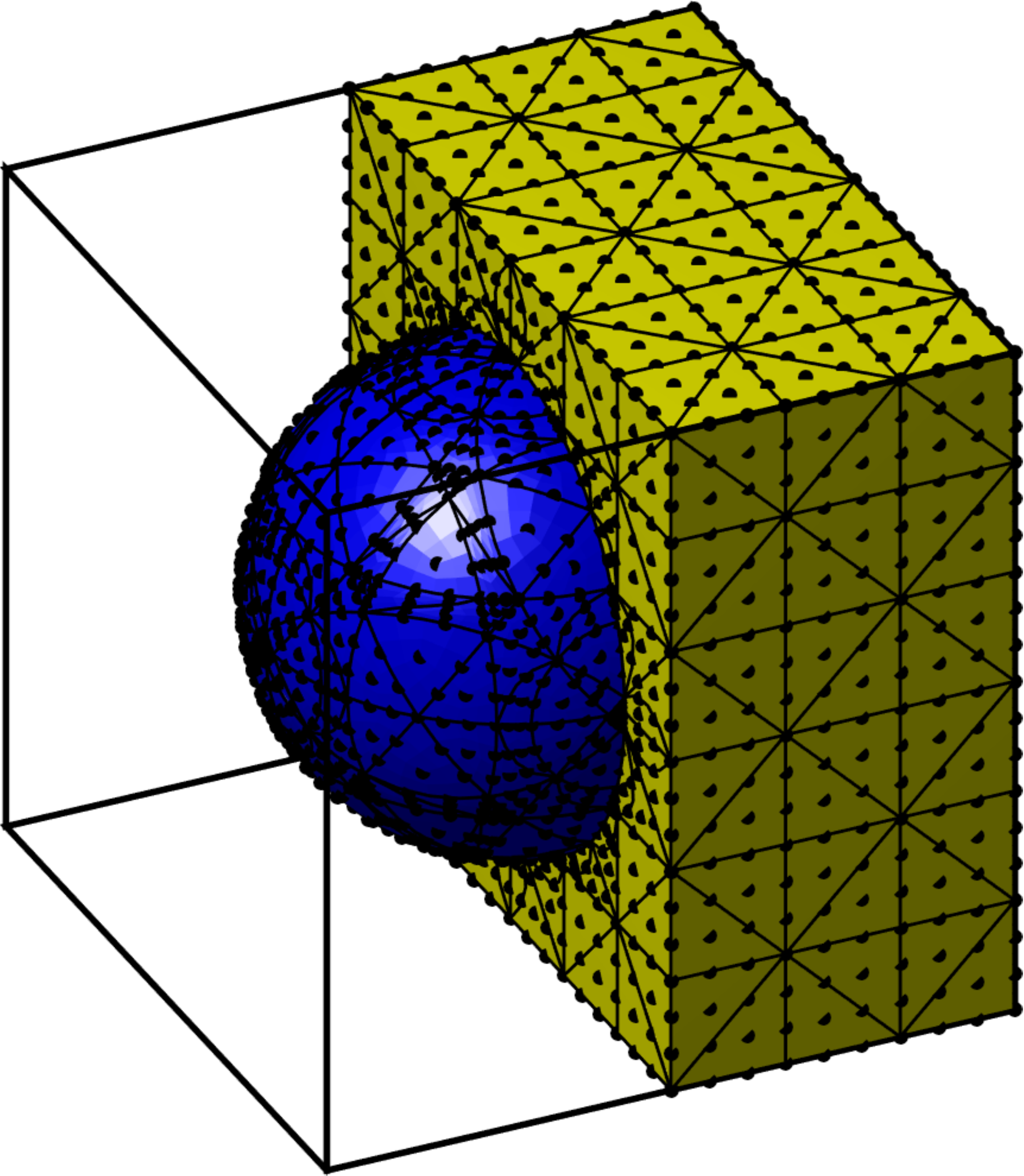}}\qquad\subfigure[exact solution]{\includegraphics[height=4cm]{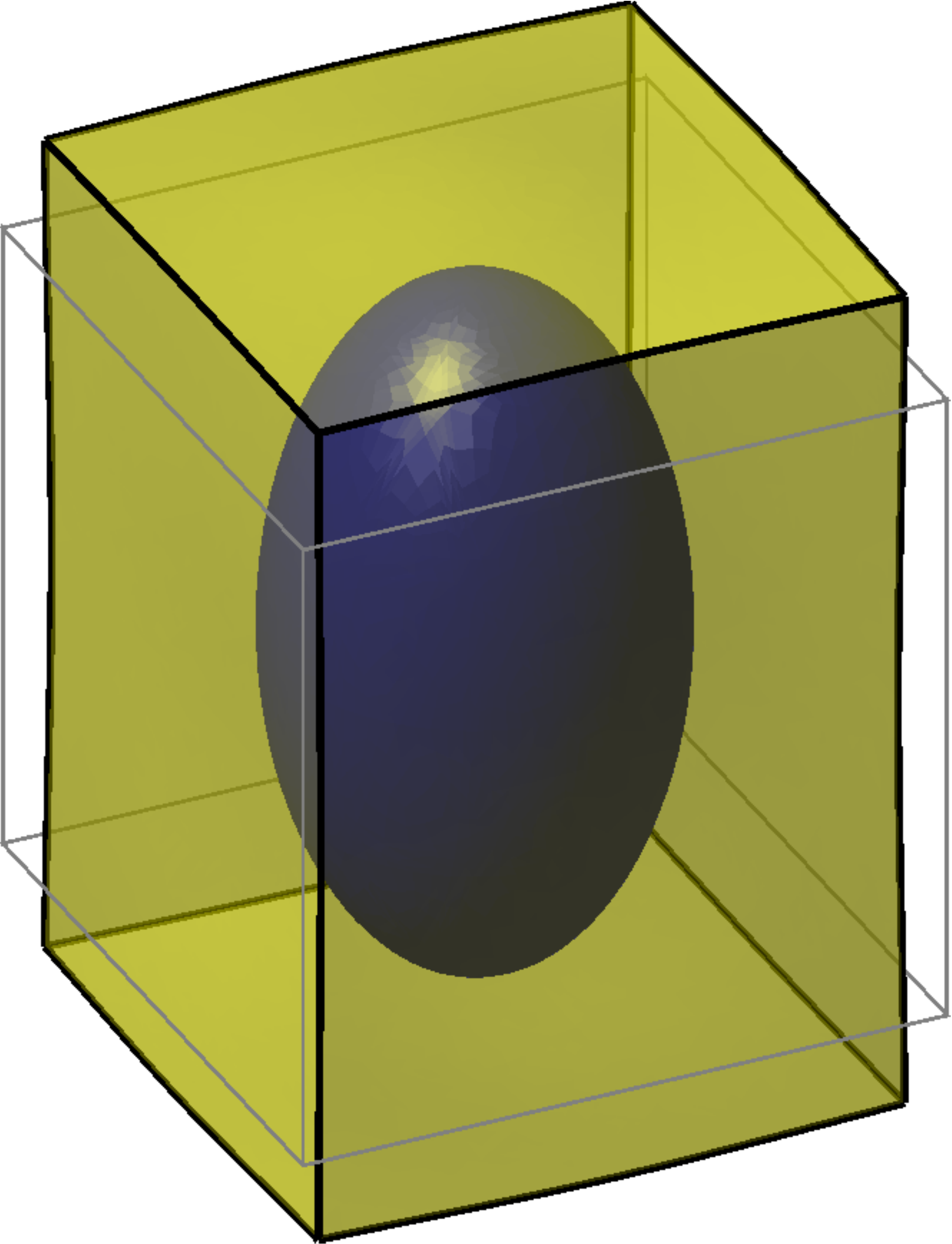}}

\caption{(a) part of the generated higher-order mesh with inclusion, (d) deformed
configuration.}

\label{fig:TestCase3dB_Sketch} 
\end{figure}

This test case is the extension of the square shell with circular
inclusion from Section \ref{XX_SquareWithCircInclusion} to three
dimensions. Fig.~\ref{fig:TestCase3dB_Sketch}(a) shows an example
for a resulting conforming mesh and elements inside the sphere are
plotted in blue. The exact solution is found in \cite{Goodier_1933a}
and also repeated in the Appendix \ref{XX_ExactSolCubeWithSphereInclusion}.
The corresponding deformed configuration is seen in Fig.~\ref{fig:TestCase3dB_Sketch}(b).
The convergence study is along the lines of Section \ref{XX_CubeWithHole}
and results are displayed in Fig.~\ref{fig:TestCase3dA_Res}. Again,
the convergence rates are optimal and condition numbers behave as
expected.

\begin{figure}
\centering

\subfigure[approximation error]{\includegraphics[height=4cm]{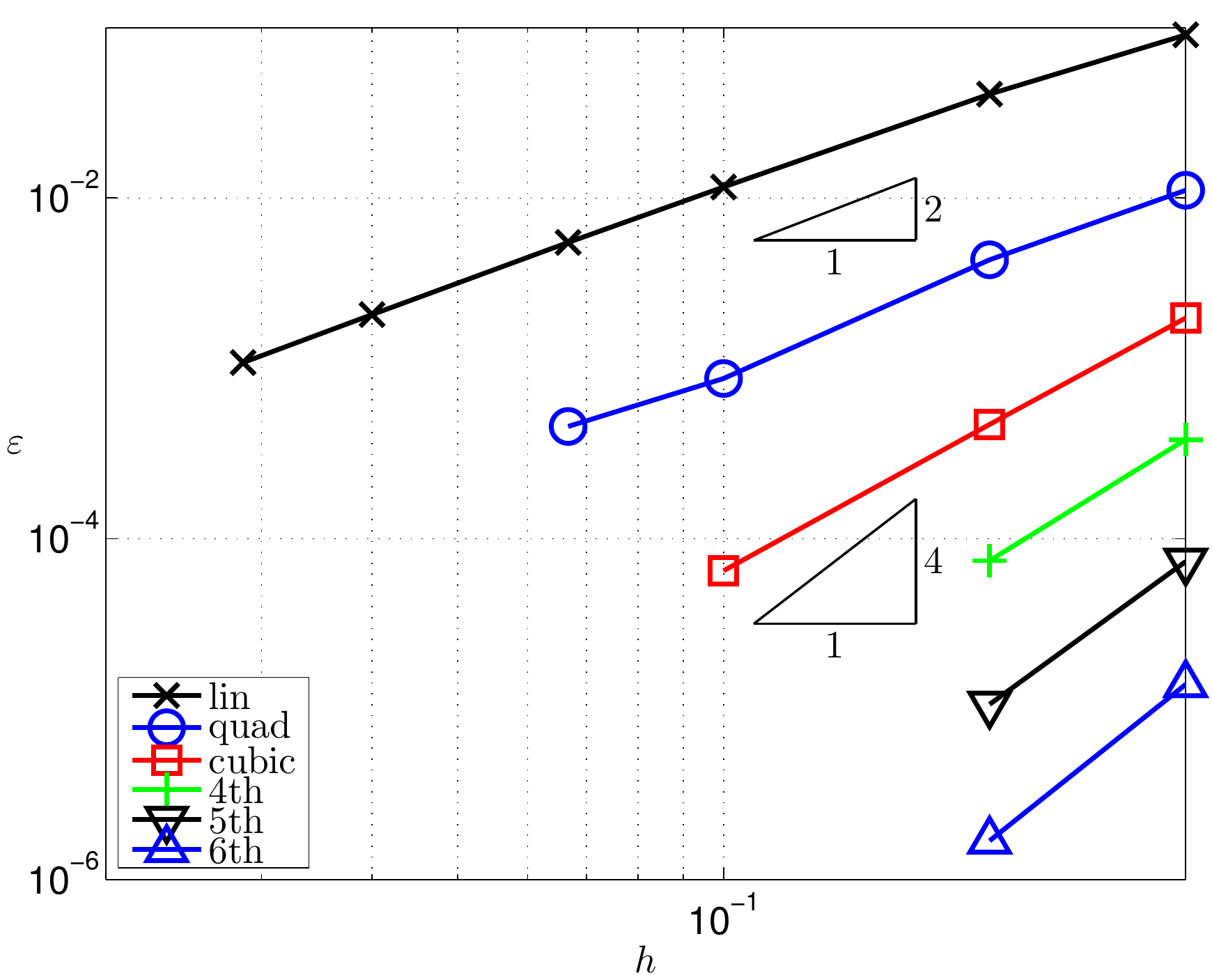}}\qquad\subfigure[condition number]{\includegraphics[height=4cm]{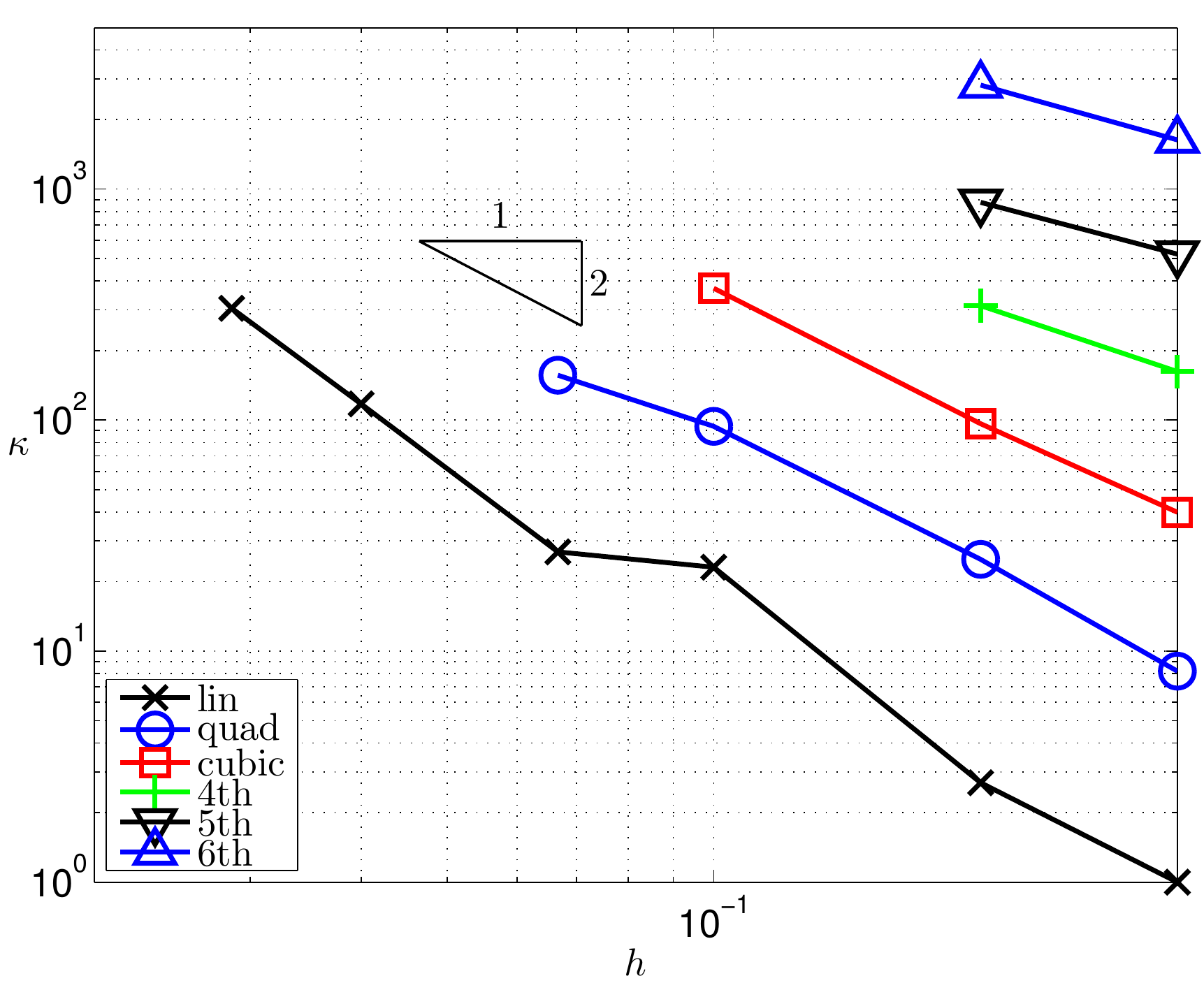}}

\caption{(a) The approximation error for the cube mesh with sphercial inclusion,
(b) the corresponding condition numbers of the system matrices.}

\label{fig:TestCase3dB_Res} 
\end{figure}

\subsection{Bi-material gyroid}

\begin{figure}
\centering

\subfigure[zero-level set]{\includegraphics[height=5cm]{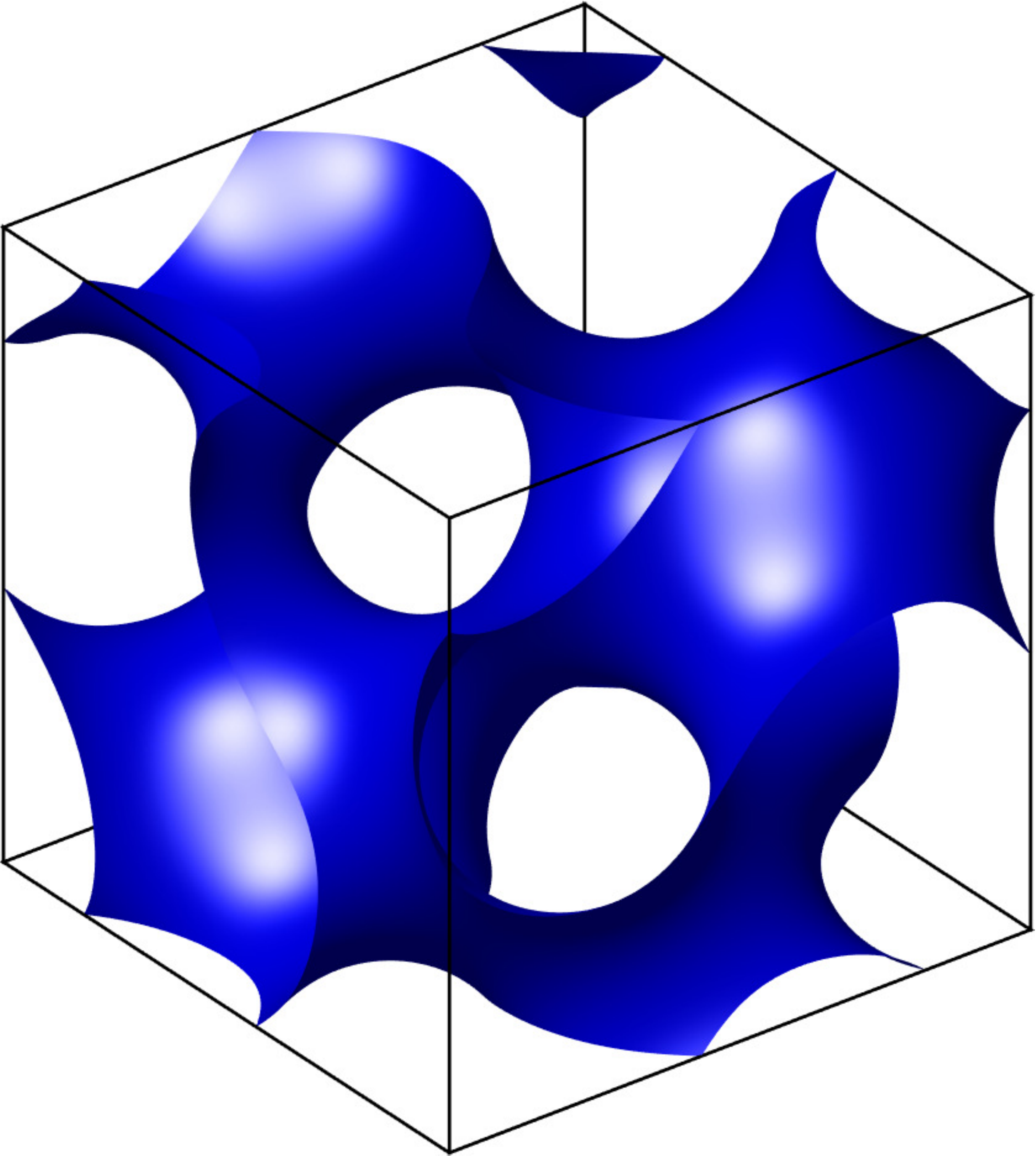}}\quad\subfigure[bi-material domain]{\includegraphics[height=5cm]{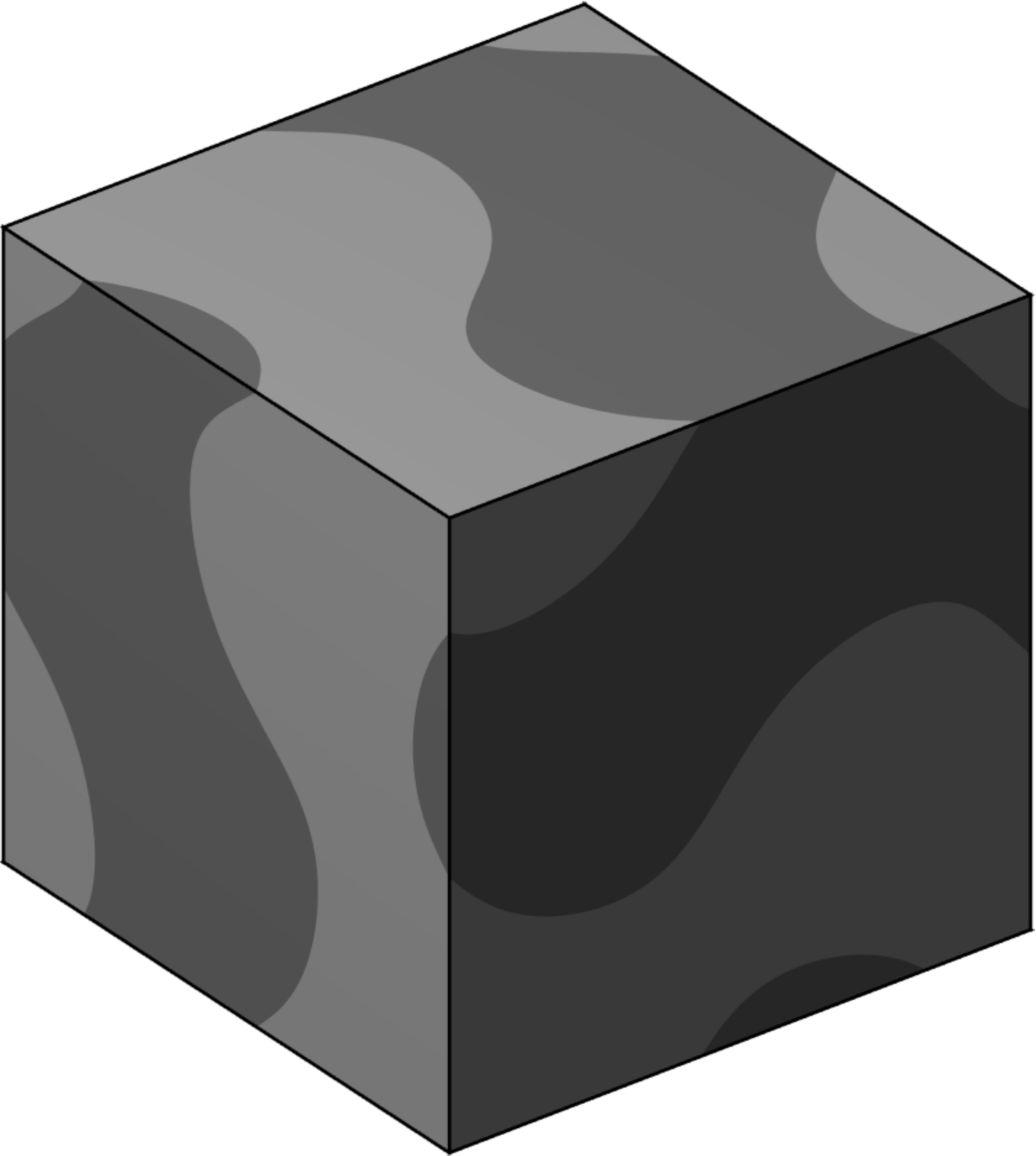}}\quad\subfigure[mesh]{\includegraphics[height=5cm]{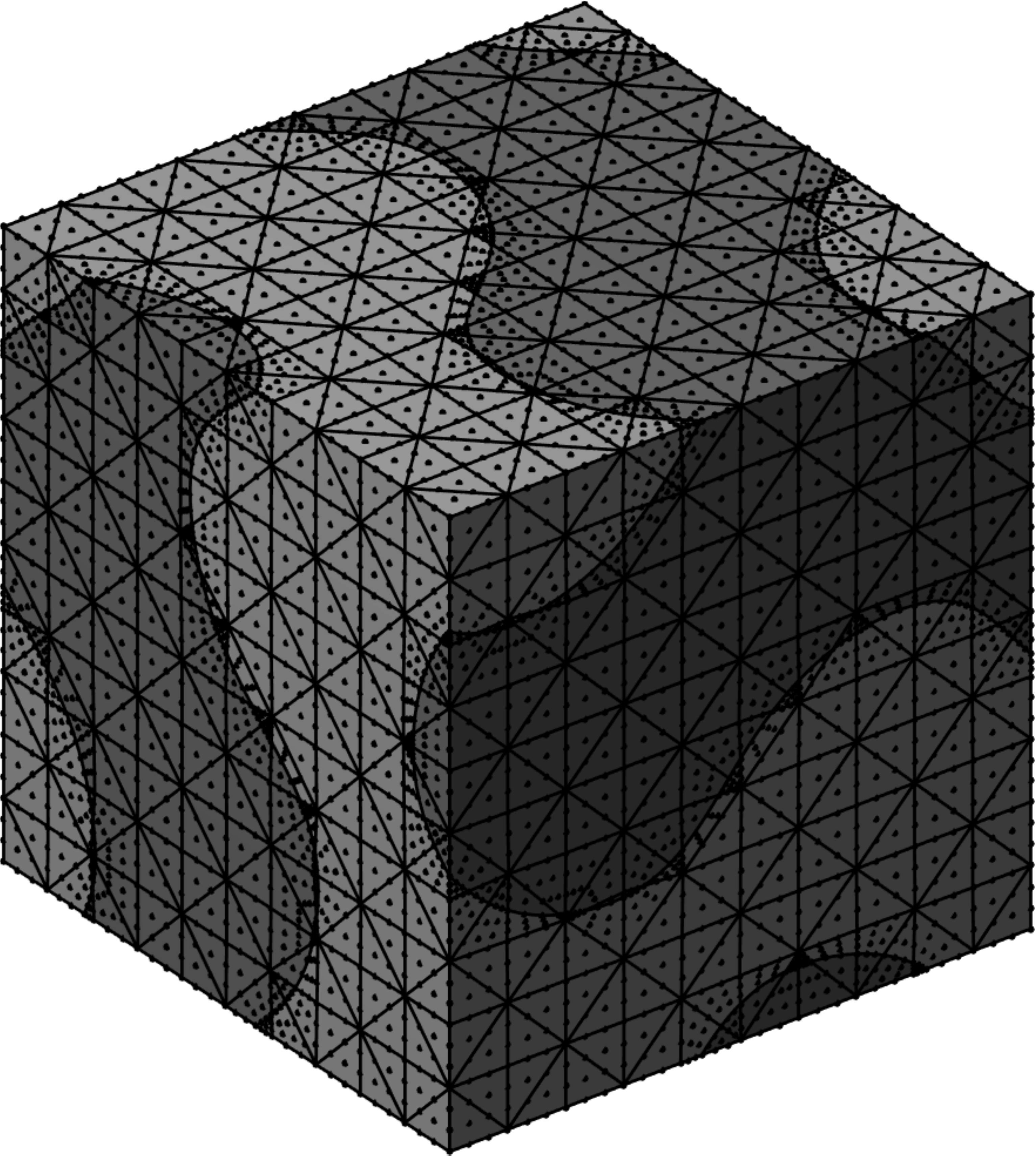}}

\caption{(a) The zero-level set, (b) the implied domain composed by two different
materials, (c) example for an automatically generated mesh.}

\label{fig:TestCase3dD_Sketch} 
\end{figure}

Another example of a domain in $\left[-1,1\right]^{3}$ composed by
two different materials is considered next where the materials are
seperated by the so-called ``gyroid'' surface, implied by the zero-level
set of
\[
\phi\left(\vek x\right)=\sin x^{\star}\cdot\cos z^{\star}+\sin y^{\star}\cdot\cos x^{\star}+\sin z^{\star}\cdot\cos y^{\star}
\]
with $x^{\star}=\pi\cdot\left(x+q\right)$, $y^{\star}=\pi\cdot\left(y+q\right)$,
$z^{\star}=\pi\cdot\left(z+q\right)$ with $q=0.123456$. See Fig.~\ref{fig:TestCase3dD_Sketch}(a)
and (b) for a representation of the zero-isosurface and the resulting
domain. Fig.~\ref{fig:TestCase3dD_Sketch}(c) shows an example of
an automatically generated mesh of order $3$. The same material parameters
as in Section \ref{XX_SquareWithCircInclusion} are chosen. All displacements
on the boundaries are fixed and the domain is loaded by a body force
of $f_{z}=-1$ in vertical direction.

\begin{figure}
\centering

\subfigure[approximation error]{\includegraphics[height=4cm]{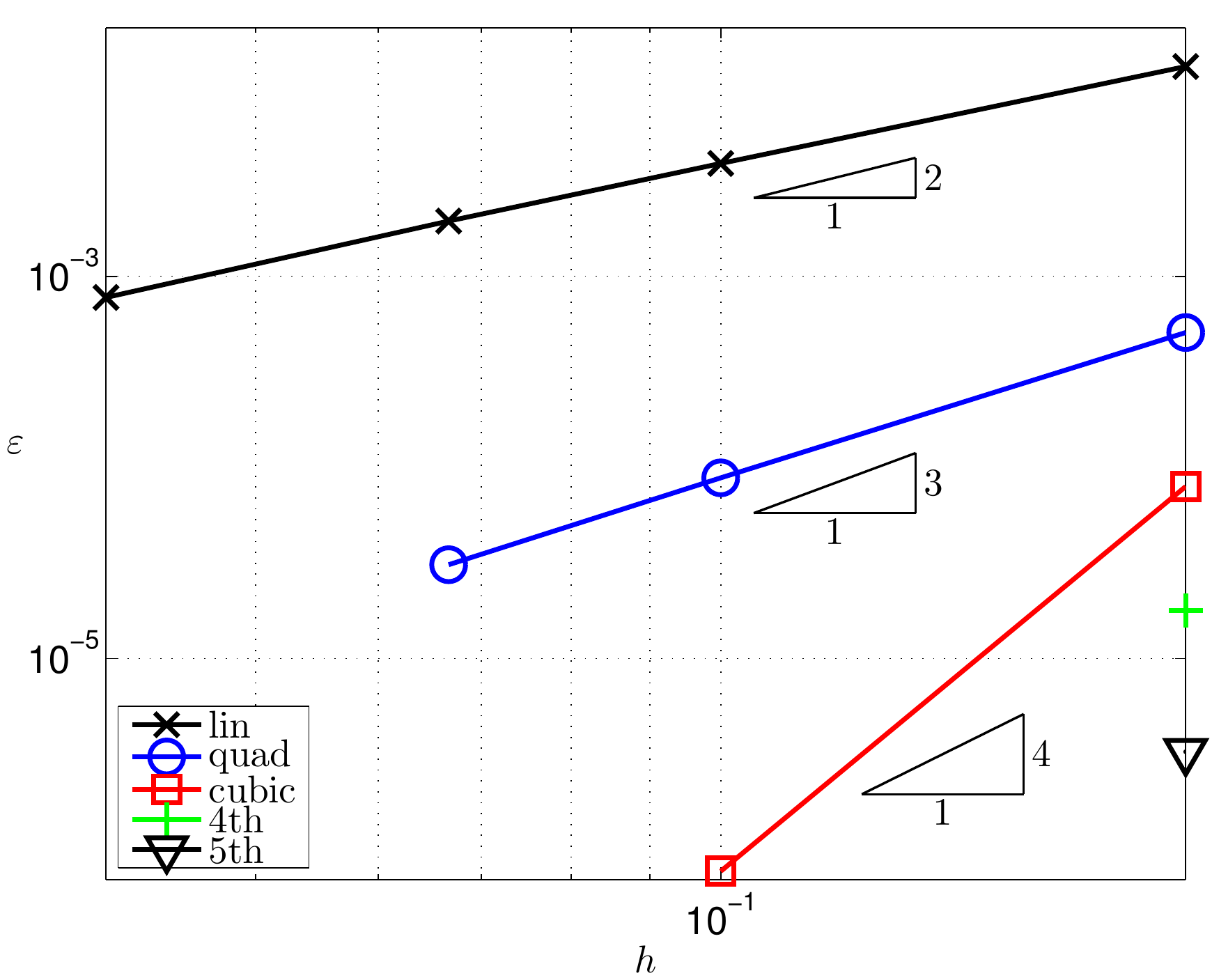}}\quad\subfigure[approximation error]{\includegraphics[height=4cm]{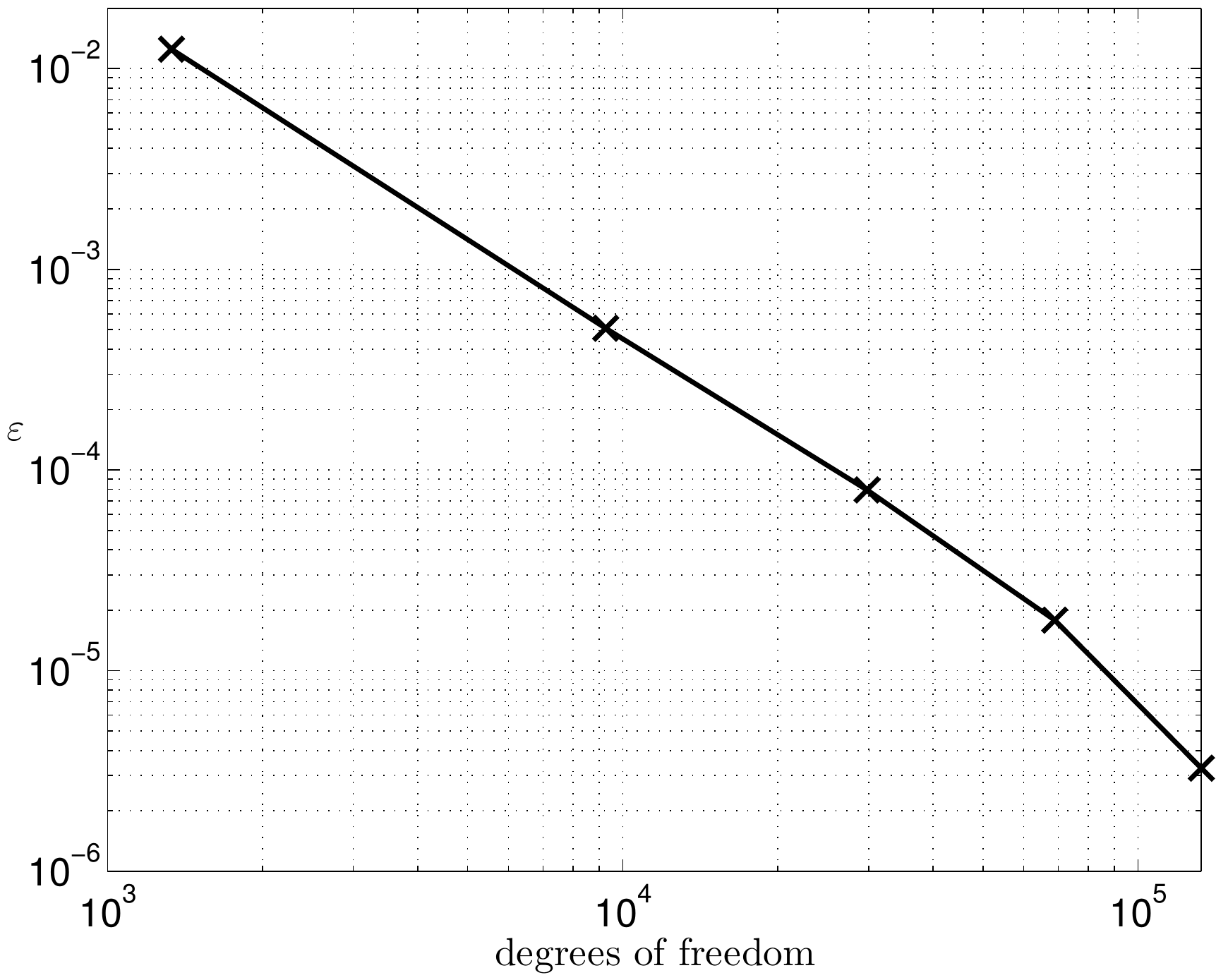}}\quad\subfigure[condition number]{\includegraphics[height=4cm]{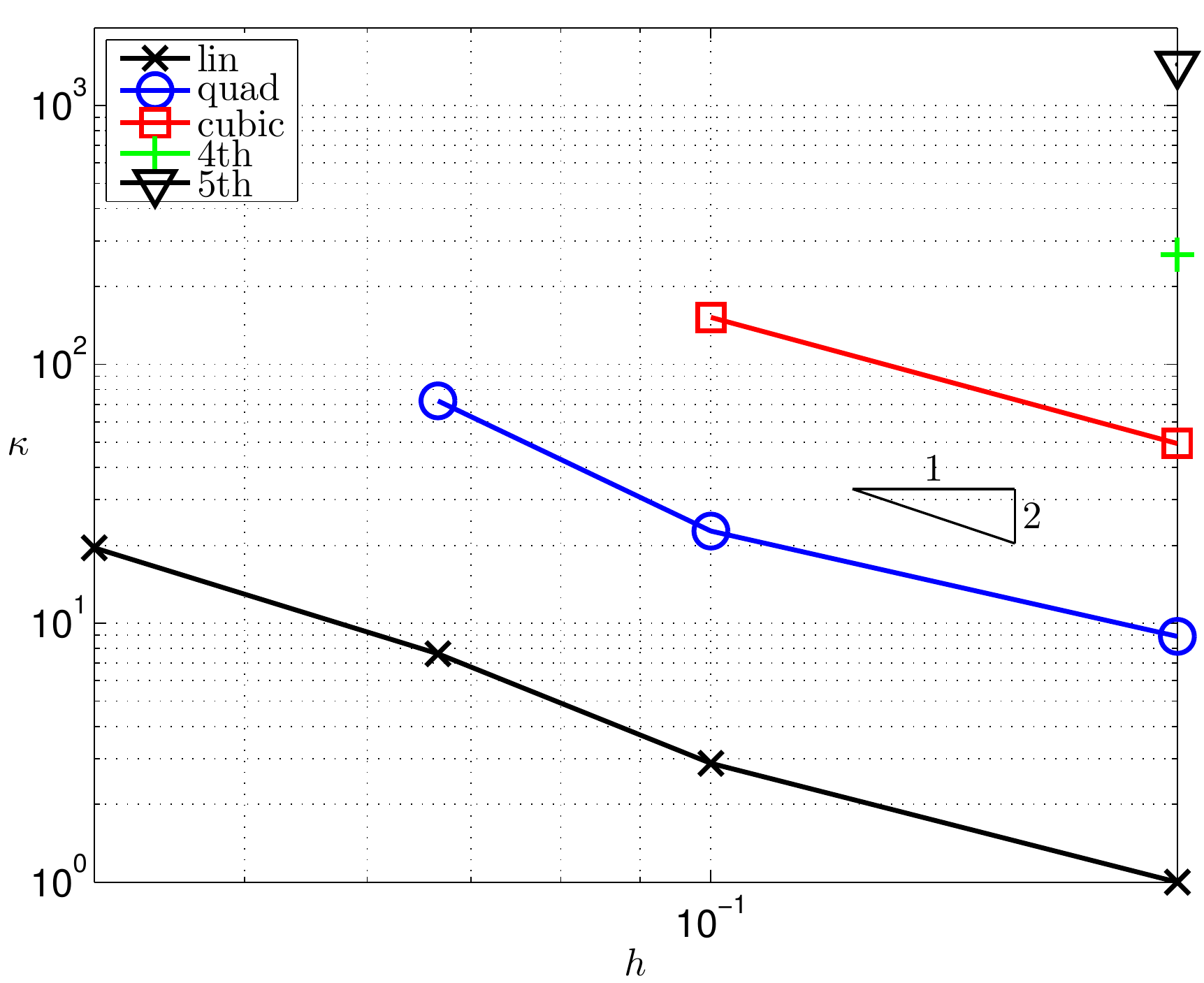}}

\caption{(a) The approximation error for the bi-material gyroid, (b) approximation
error over the number of degrees of freedom, (c) the corresponding
condition numbers of the system matrices.}

\label{fig:TestCase3dD_Res} 
\end{figure}

For the convergence studies, background meshes with $\{10,20,30,50\}$
elements per dimensions are chosen unless they lead to more than $10^{6}$
degrees of freedom. For this rather complex zero-level set, coarser
meshes with $<10$ elements per dimension lead to a significant number
of adaptive refinements, so that the element lengths $h$ vary too
much to obtain representative values for the convergence plots. Fig.~\ref{fig:TestCase3dD_Res}(a)
shows the convergence rates in the same style than before, however,
the limitation on the number of degrees of freedom leads to less data
points. Therefore, Fig.~\ref{fig:TestCase3dD_Res}(b) shows convergence
results on meshes with $10$ elements per dimension only but with
different element orders. The error is plotted with respect to the
degrees of freedom and expontential convergence is obtained. Finally,
the condition number is seen in Fig.~\ref{fig:TestCase3dD_Res}(c)
and is bounded as expected.

\subsection{Cantilever tube}

\begin{figure}
\centering

\subfigure[zero-level sets]{\includegraphics[width=6.5cm]{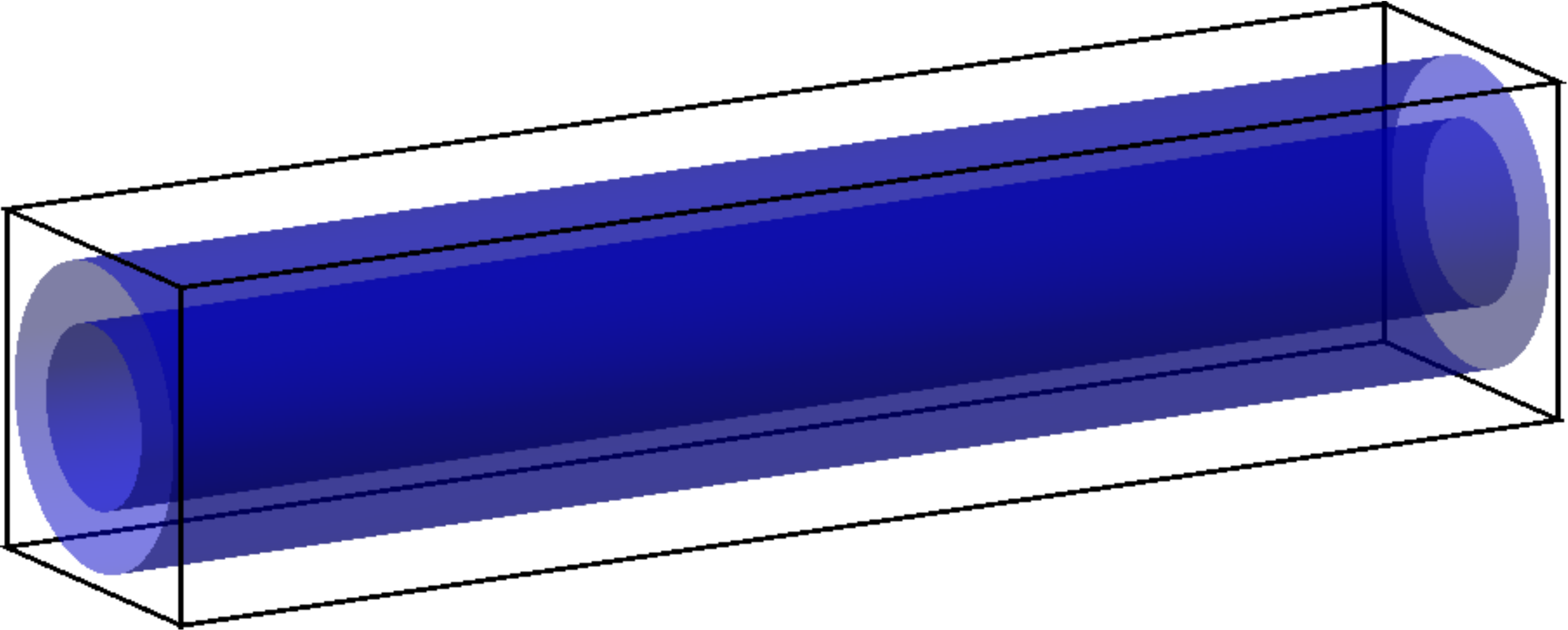}}\qquad\subfigure[background mesh]{\includegraphics[width=6.5cm]{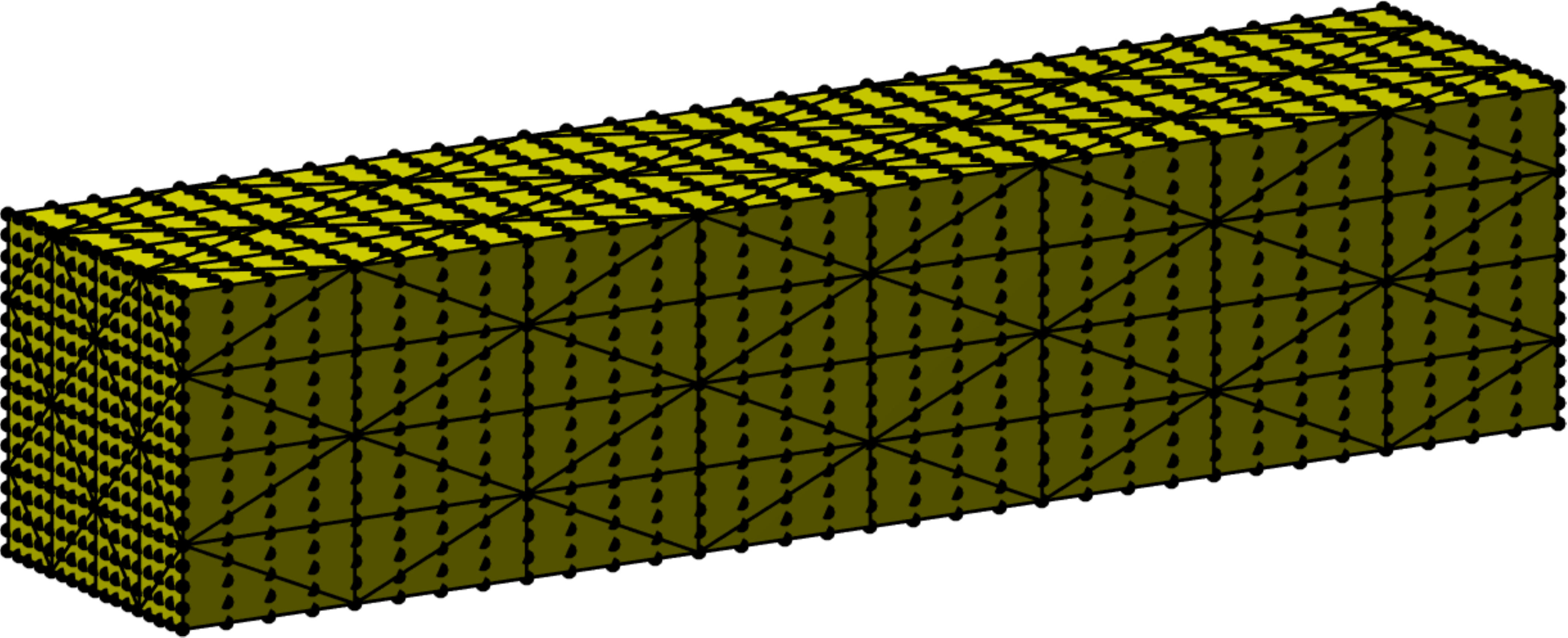}}

\subfigure[generated mesh]{\includegraphics[width=6.5cm]{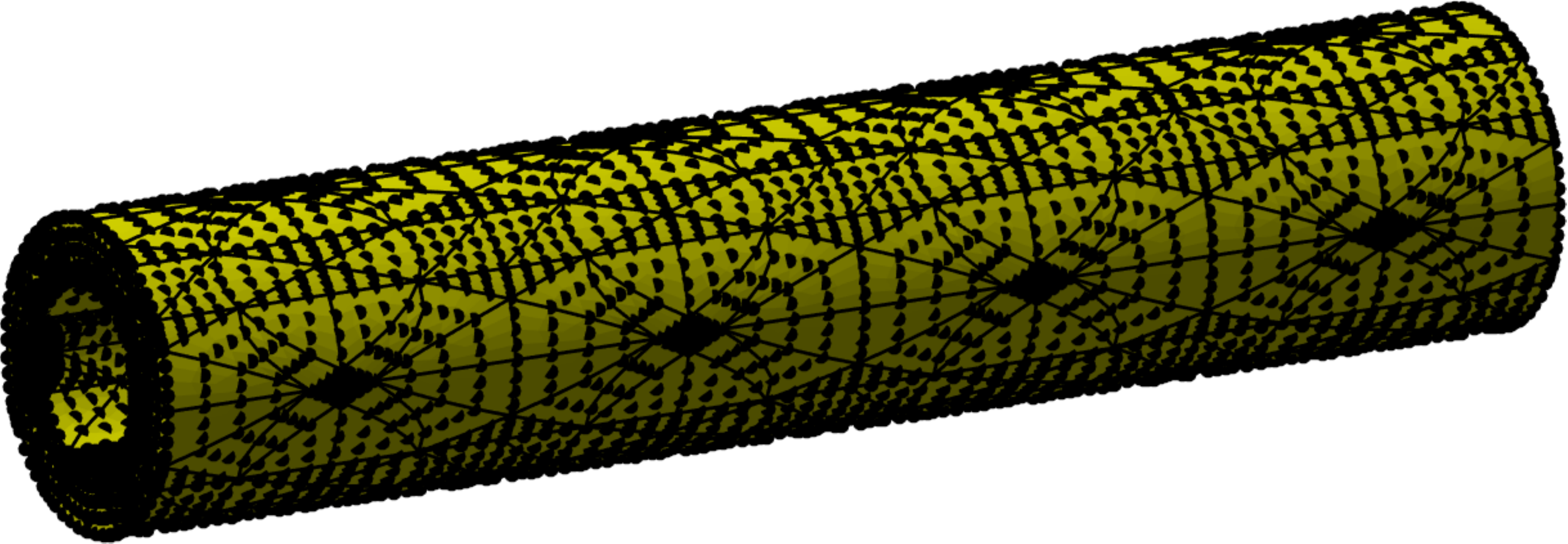}}\qquad\subfigure[exact solution]{\includegraphics[width=6.5cm]{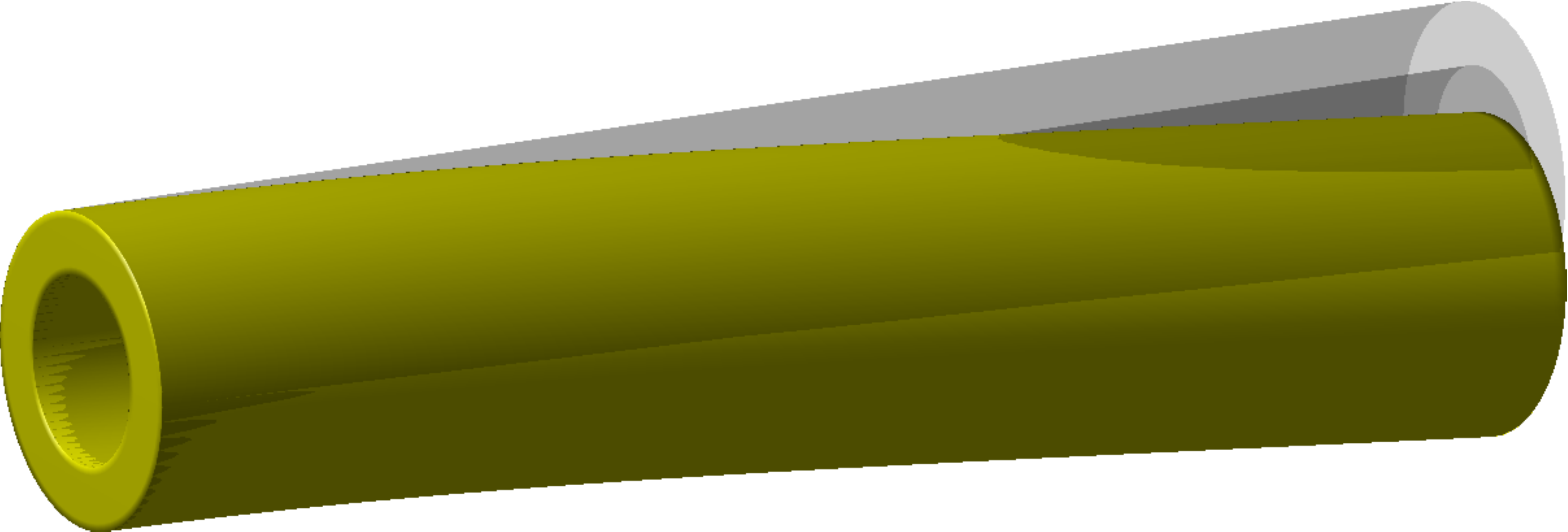}}

\caption{(a) Zero-level sets for the pipe test case, (b) an example of a background
mesh and (c) the resulting conforming mesh, (d) the deformed configuration.}

\label{fig:TestCase3dC_Sketch1}
\end{figure}

\begin{figure}
\centering

\subfigure[2D background mesh]{\includegraphics[height=3cm]{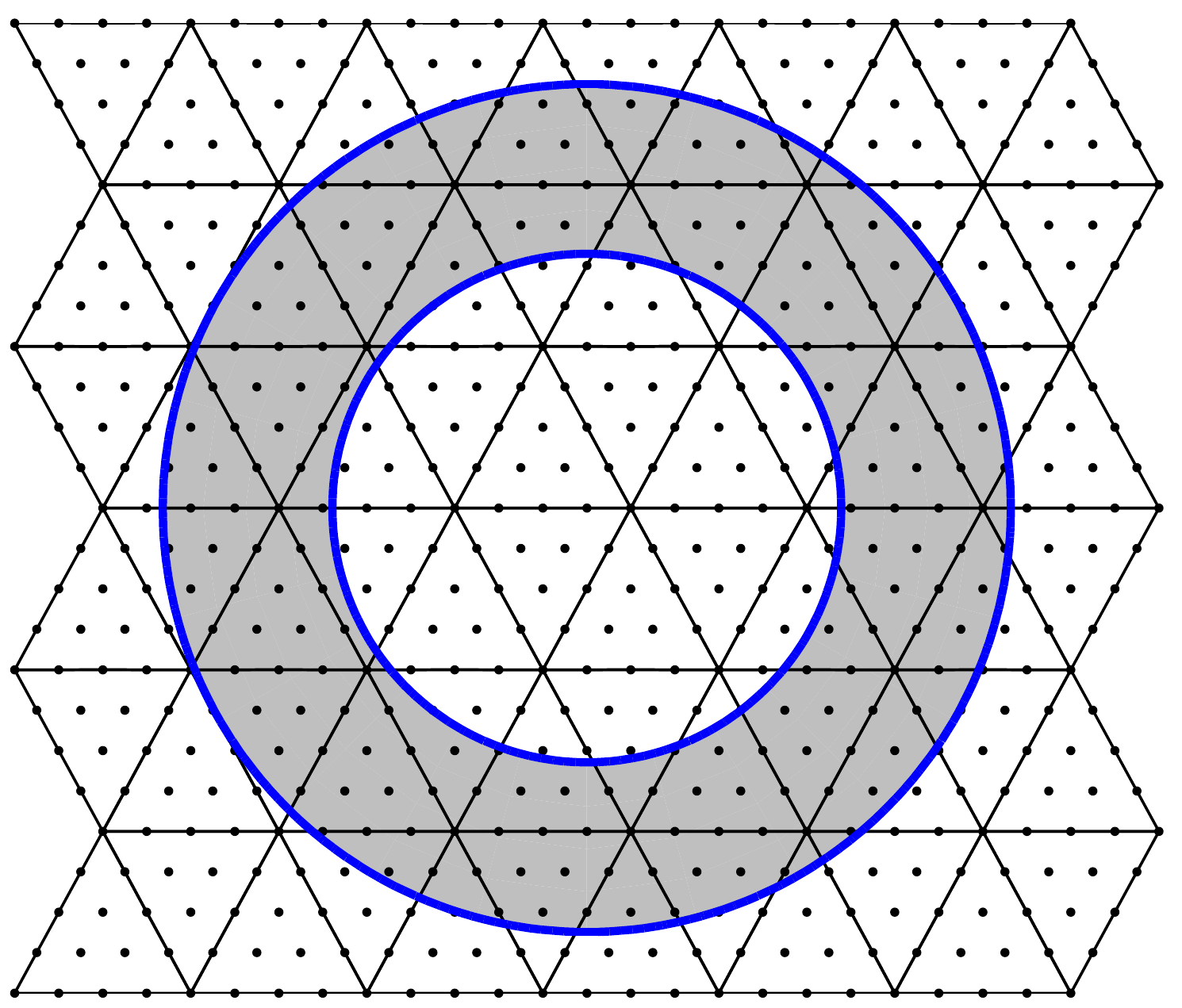}}\qquad\subfigure[generated 2D mesh]{\includegraphics[height=3cm]{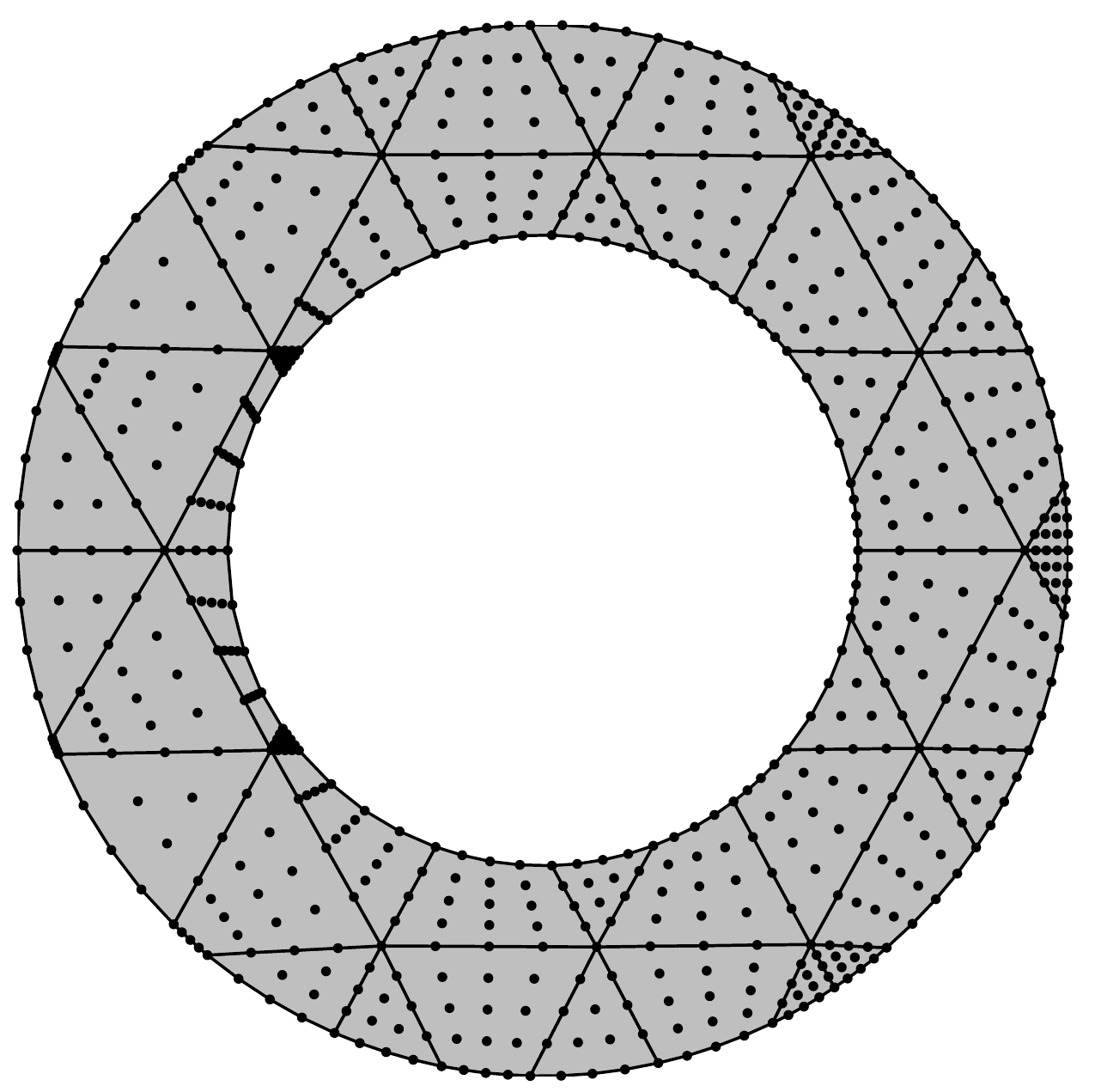}}\qquad\subfigure[extruded 3D mesh]{\includegraphics[width=6.5cm]{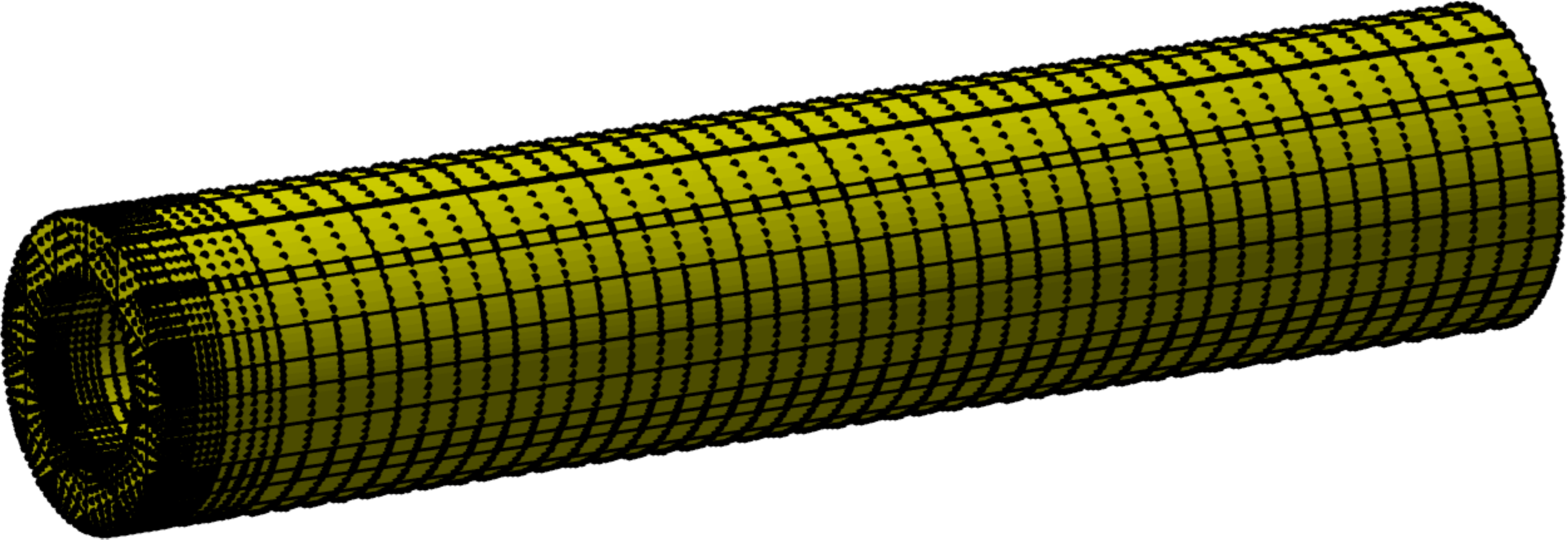}}

\caption{Mesh generation for the pipe test case based on the extrusion of a
2D mesh, (a) and (b) show the 2D background mesh and the resulting
conforming mesh, (c) extruded mesh in $x$-direction with refinement
at the clamped side.}

\label{fig:TestCase3dC_Sketch2}
\end{figure}

The next test case is more technical. We consider a pipe of length
$L=5\mathrm{m}$ with an inner radius $r_{\mathrm{i}}=0.3\mathrm{m}$
and an outer radius $r_{\mathrm{o}}=0.5\mathrm{m}$. The pipe is clamped
on one side, i.e., all displacement components are enforced to vanish
there. The material is composed by steel with $E=2.1\cdot10^{8}\,\unitfrac{kN}{m^{2}}$
and $\nu=0.3$. The beam is loaded by gravity acting as a body force
of $f_{z}=-78.5\,\unitfrac{kN}{m^{3}}$. The background mesh is given
in $\Omega_{\mathrm{BG}}=\left[0,5\right]\times\left[-0.55,0.55\right]^{2}$
and the pipe walls are implied by the two level-set functions
\[
\phi_{\mathrm{i/o}}=\sqrt{y^{2}+z^{2}}-r_{\mathrm{i/o}},
\]
see Fig.~\ref{fig:TestCase3dC_Sketch1}(a). An example background
mesh is seen in Fig.~\ref{fig:TestCase3dC_Sketch1}(b) and the resulting
conforming mesh in Fig.~\ref{fig:TestCase3dC_Sketch1}(c). The deformed
configuration scaled by a factor of $1000$ is shown in Fig.~\ref{fig:TestCase3dC_Sketch1}(d).
It is also useful to generate a mesh for the pipe by first generating
a conforming 2D mesh of the cross-sectional area based on a 2D background
mesh as shown in Figs.~\ref{fig:TestCase3dC_Sketch2}(a) and (b)
and then extruding this in $x$-direction generating prismatic and
hexahedral elements. This also allows for an efficient refinement
near the clamped side, see e.g., Fig.~\ref{fig:TestCase3dC_Sketch2}(c).

The situation is comparable to a clamped Bernoulli beam governed by
the differential equation $EI\cdot w(x)''''=q(x)$ with the second
moment of area $I=\nicefrac{\pi}{4}\cdot\left(r_{\mathrm{o}}^{4}-r_{\mathrm{i}}^{4}\right)$,
the line load $q(x)=-f_{z}\cdot A$ and the sectional area $A=\pi\cdot\left(r_{\mathrm{o}}^{2}-r_{\mathrm{i}}^{2}\right)$.
The line load $q$ and the deflection $w$ are positive in downward
direction, i.e., $w=-u_{z}$. With the left side ($x=0$) clamped
and the right side ($x=L$) fully free, the analytical solution for
the bending curve is
\[
w(x)=\dfrac{q}{EI}\cdot\left(\dfrac{L^{2}}{4}x^{2}-\dfrac{L}{6}x^{3}+\dfrac{1}{24}x^{4}\right).
\]
This yields a deflection on the right side of $w\left(L\right)=\frac{qL^{4}}{8EI}=0.343575\mathrm{mm}$
and a stored energy of $\mathfrak{e}=\nicefrac{1}{2}\int M\left(x\right)\cdot\varkappa\left(x\right)\, dx=\frac{EI}{2}\int\left[w''\left(x\right)\right]^{2}\, dx=\frac{q^{2}L^{5}}{40EI}\approx0.013557\mathrm{kNm}$.
An overkill FEM-solution of the three-dimensional problem yields $\max\left(\left|u_{z}\right|\right)=0.3605575\pm10^{-6}\mathrm{mm}$
and $\mathfrak{e}=0.01473635\mathrm{\pm10^{-7}kNm}$ which is quite
similar and taken for the convergence studies.

\begin{figure}
\centering

\subfigure[energy error]{\includegraphics[height=4cm]{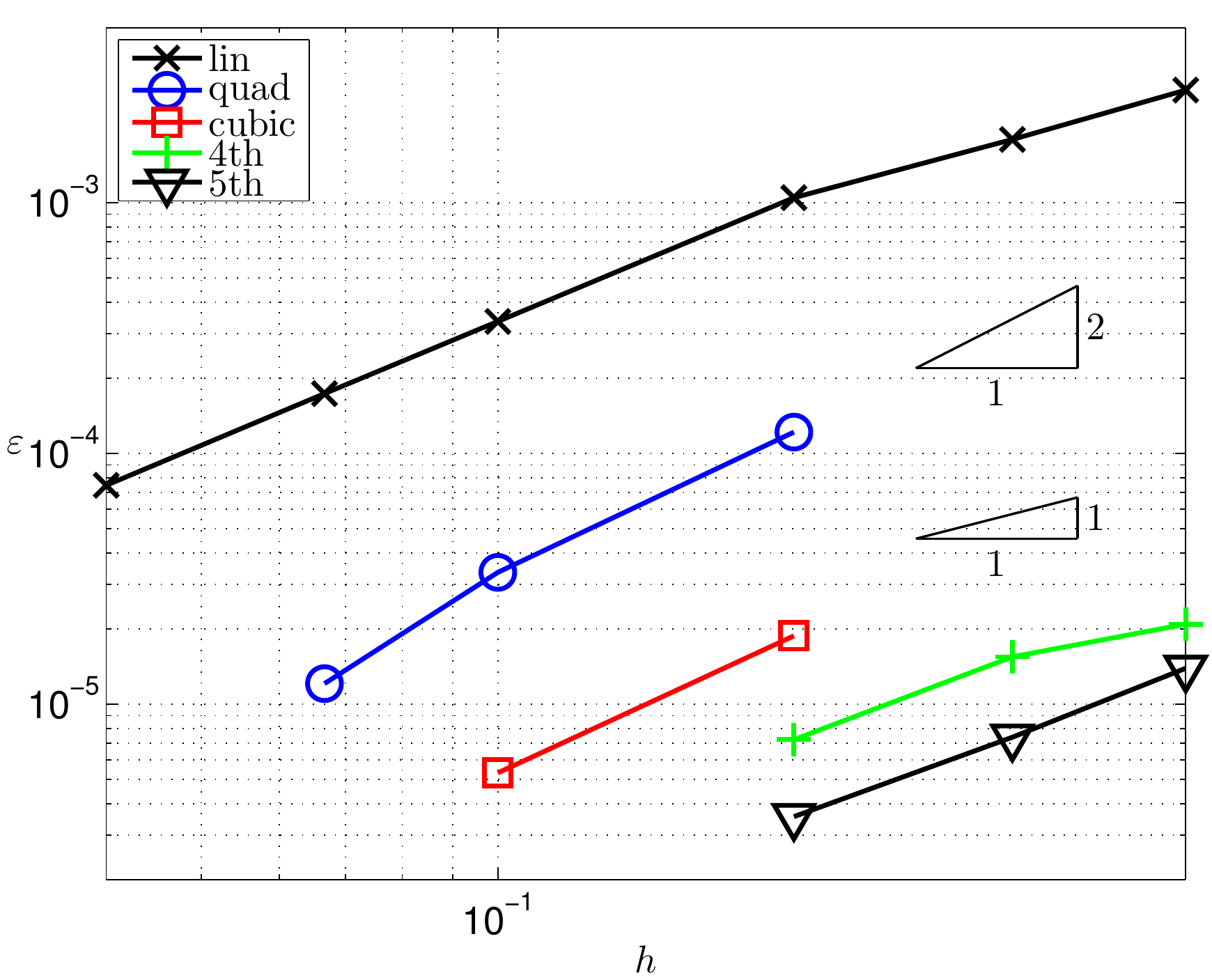}}\qquad\subfigure[bending error]{\includegraphics[height=4cm]{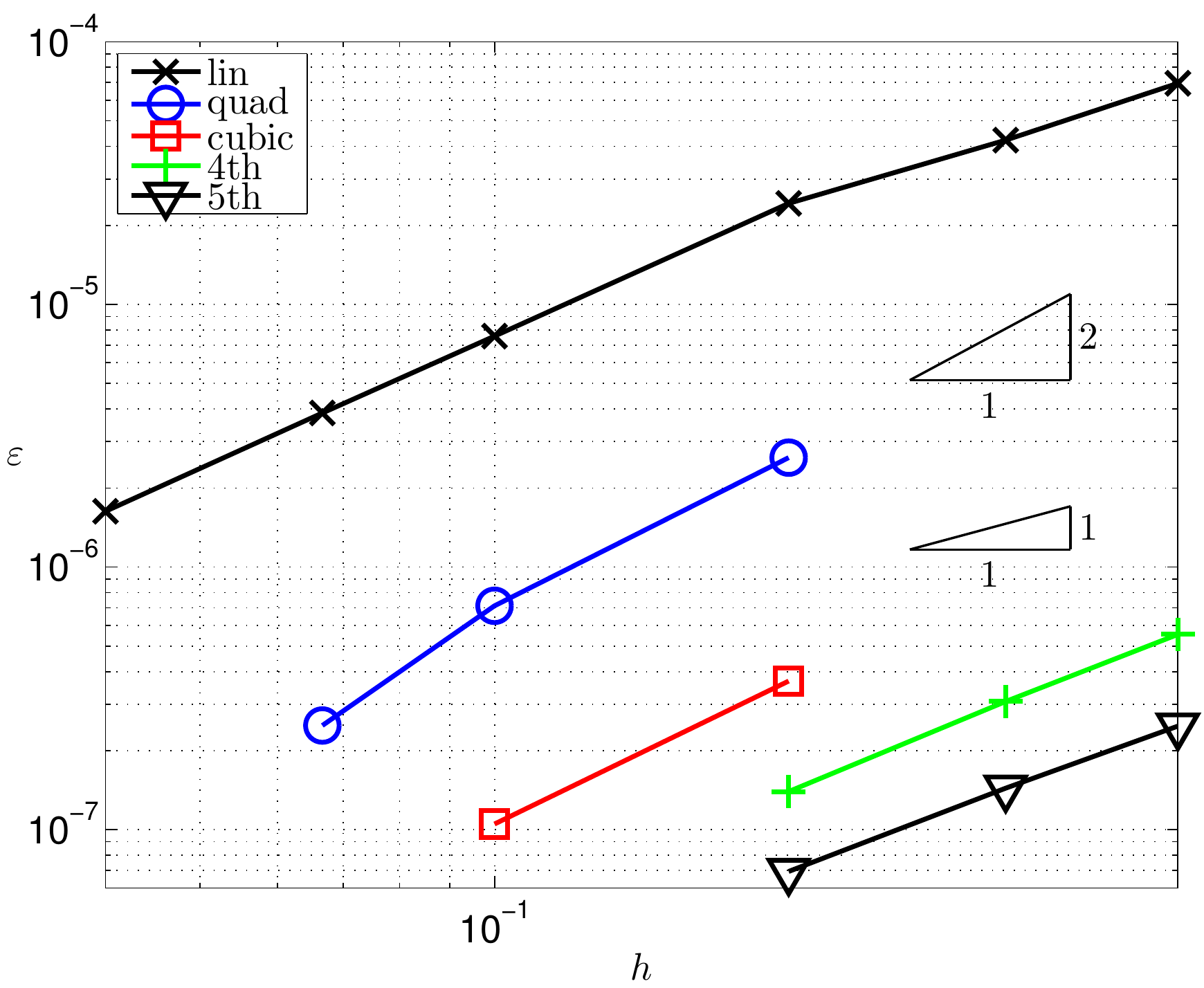}}

\subfigure[energy error]{\includegraphics[height=4cm]{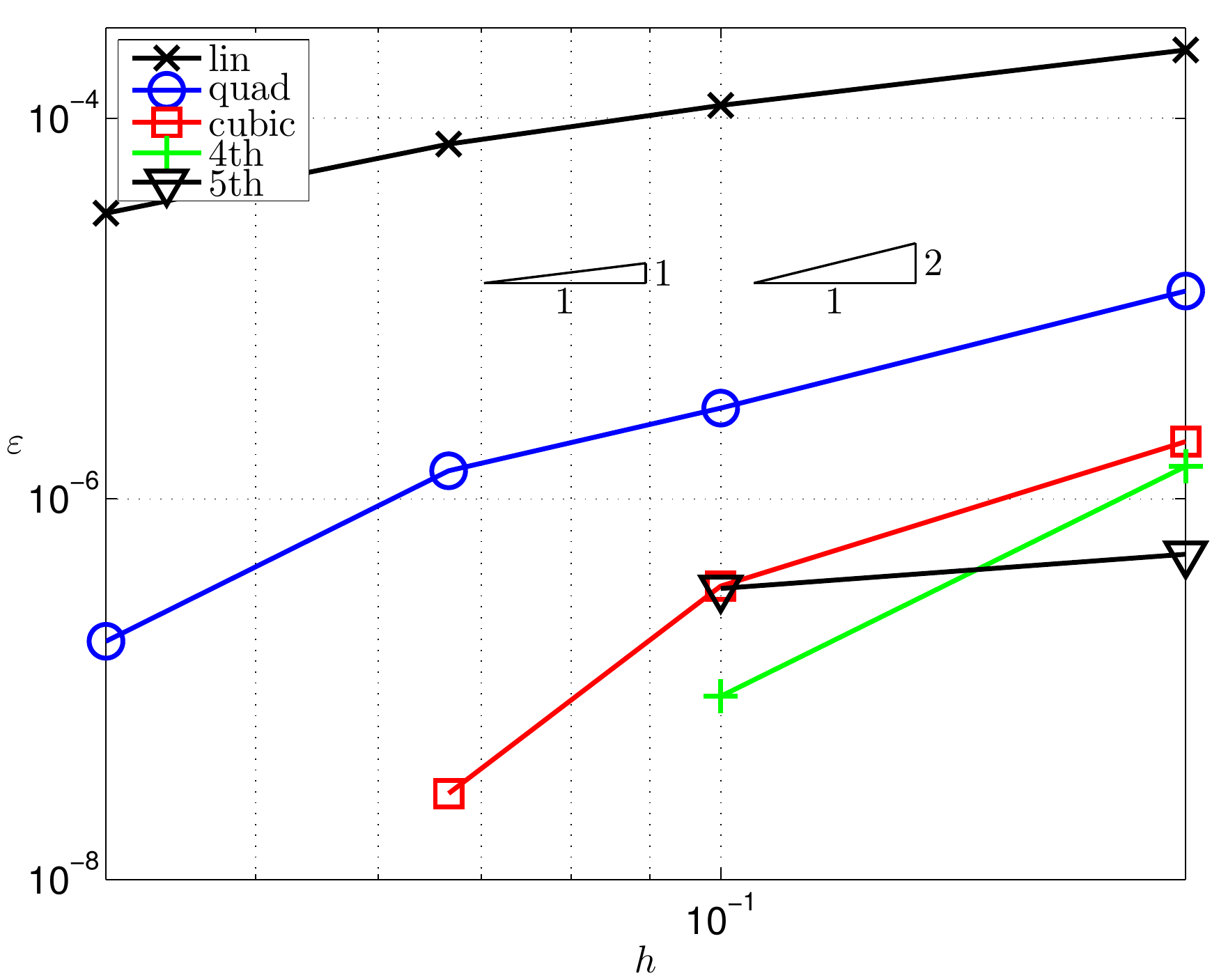}}\qquad\subfigure[bending error]{\includegraphics[height=4cm]{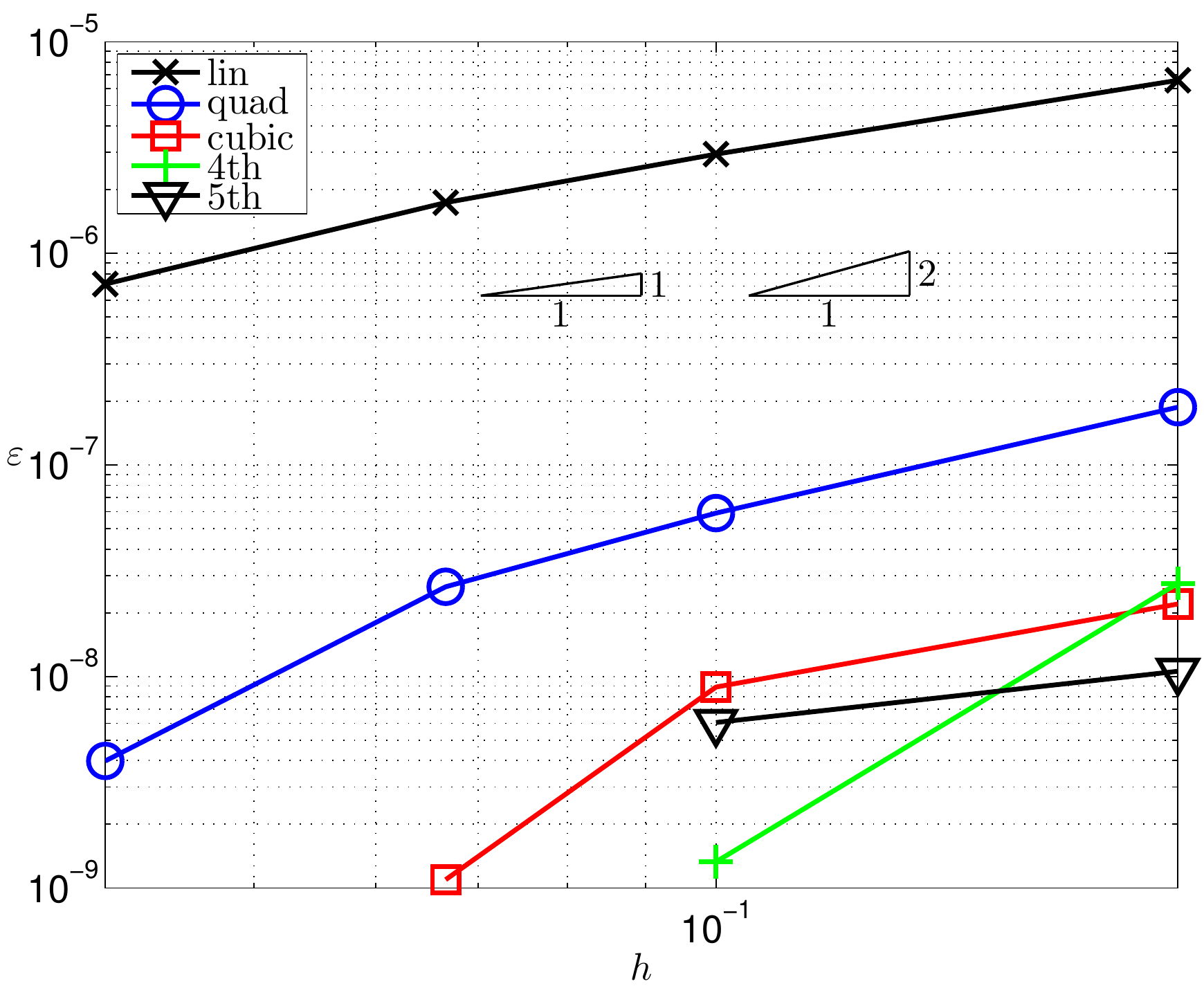}}

\caption{Convergence of the errors in the stored energy and deflection on the
right side for the pipe test case, (a) and (b) for the meshes generated
from 3D background meshes, (c) and (d) for meshes extruded from the
2D annulus meshes.}

\label{fig:TestCase3dC_Res} 
\end{figure}

Results are seen in Fig.~\ref{fig:TestCase3dC_Res}(a) and (b) for
meshes generated from 3D background meshes as in Fig.~\ref{fig:TestCase3dC_Sketch1}(b),
and in Fig.~\ref{fig:TestCase3dC_Res}(c) and (d) for meshes generated
from extruding 2D meshes as shown in Fig.~\ref{fig:TestCase3dC_Sketch2}.
The convergence rates are sub-optimal as expected due to the singularties
in the stresses at the clamped side. Nevertheless, it is seen that
higher-order elements are able to significantly improve the results.
Obviously, the refinement on the clamped side as for the extruded
meshes further improves the results.

\section{Conclusions\label{X_Conclusions}}

A higher-order CDFEM is proposed which automatically generates higher-order,
conforming meshes based on background meshes and level-set data. The
decomposition of cut elements into conforming sub-elements has been
described before, e.g., in \cite{Fries_2015a,Fries_2016a,Fries_2016b}
in the context of integration and interpolation. Herein, this idea
is extended to the approximation of BVPs. In addition to the fact
that the sub-elements must conform to the zero-level sets, this involves
the following additional challenges: The element set must be continuous
accross element boundaries so that hanging nodes are avoided. Therefore,
an adaptive procedure which guarantees the regularity of the background
mesh is suggested. Furthermore, the generated elements must be shape-regular
wherefore a node manipulation is proposed which (slightly) moves nodes
near the zero-level sets to guarantee bounded ratios of the areas/volumes
on the two sides of the cut elements. A suitable finite element mesh
composed by higher-order elements may then be generated from the element
set, i.e., the connectivity information is set up.

In particular, the combination of the automatic mesh generation and
adaptive refinements, not only for elements where the decomposition
fails, but also where a more accurate geometry description or approximation
of the BVP is desired, is found to be a key ingredient of the proposed
method. Numerical results are presented in the context of linear elasticity
without loss of generality of the method. Elements up to order $6$
are investigated herein and higher-order convergence rates are achieved. 

The resulting method is stable, efficent and achieves optimal, higher-order
convergence rates in two and three dimensions. As such, it is an attractive
alternative to FDMs and the XFEM. Remaining challenges include \emph{moving}
interfaces and iterative solvers. In a forthcoming part of this series
of publications, we shall report on a higher-order FDM where the shape
functions of the background mesh are used for the approximation and
comparisons to the CDFEM discussed herein will be made.

\section{Appendix\label{X_Appendix}}

\subsection{Exact solution for the square shell with circular hole\label{XX_ExactSolSquareWithHole}}

The exact solution for this problem of an infinite plate with a traction-free
circular hole under uniaxial tension with $\sigma_{0}$ is, e.g.,
found in \cite{Szabo_1991a,Liu_2002a}. It is given in polar coordinates
$\left(r,\theta\right)$ as
\begin{eqnarray*}
u_{x}\left(r,\theta\right) & = & \sigma_{0}\cdot\frac{R}{8\mu}\left[\frac{r}{R}\left(\varkappa+1\right)\cos\theta+\frac{2R}{r}\left(\left(1+\varkappa\right)\cos\theta+\cos3\theta\right)-\frac{2R^{3}}{r^{3}}\cos3\theta\right],\\
u_{y}\left(r,\theta\right) & = & \sigma_{0}\cdot\frac{R}{8\mu}\left[\frac{r}{R}\left(\varkappa-3\right)\sin\theta+\frac{2R}{r}\left(\left(1-\varkappa\right)\sin\theta+\sin3\theta\right)-\frac{2R^{3}}{r^{3}}\sin3\theta\right],
\end{eqnarray*}
where $R$ is the radius of the hole and $\varkappa$ is the Kolosov
constant defined as 
\[
\varkappa=\begin{cases}
3-4\nu & \quad\text{for plane strain},\\
\dfrac{3-\nu}{1+\nu} & \quad\text{for plane stress}.
\end{cases}
\]

\subsection{Exact solution for the square shell with circular inclusion\label{XX_ExactSolSquareWithCircInclusion}}

The exact solution for this test case is given, e.g., in \cite{Sukumar_2001c,Fries_2006b,Fries_2007b}.
The radius of the inclusion is $a=R$ and another scalar value $b>a$
defines a radius where a given traction is applied. The displacements
in direction of the polar coordinates are
\begin{eqnarray}
u_{r}\left(r,\theta\right) & = & \left\{ \begin{array}{cl}
\left((1-\frac{b^{2}}{a^{2}})\alpha+\frac{b^{2}}{a^{2}}\right)\cdot r, & \quad0\leq r\leq a,\\
\left((1-\frac{b^{2}}{r^{2}})\alpha+\frac{b^{2}}{r^{2}}\right)\cdot r, & \quad a<r\leq b,
\end{array}\right.\\
u_{\theta}\left(r,\theta\right) & = & 0.
\end{eqnarray}
The parameter $\alpha$ involved in these definitions is
\begin{equation}
\alpha=\frac{\left(\lambda_{1}+\mu_{1}+\mu_{2}\right)b^{2}}{\left(\lambda_{2}+\mu_{2}\right)a^{2}+\left(\lambda_{1}+\mu_{1}\right)\left(b^{2}-a^{2}\right)+\mu_{2}b^{2}}.
\end{equation}
It is trivial to transform these displacements into $x$- and $y$-direction
using the transformation matrix $T$, hence,
\[
\left[\begin{array}{c}
u_{x}\\
u_{y}
\end{array}\right]=\left[\begin{array}{cc}
\cos\theta & -\sin\theta\\
\sin\theta & \cos\theta
\end{array}\right]\cdot\left[\begin{array}{c}
u_{r}\\
u_{\theta}
\end{array}\right].
\]

\subsection{Exact solution for the cube with spherical hole\label{XX_ExactSolCubeWithHole}}

The solution is defined in spherical coordinates $\left(r,\theta,\varphi\right)$
related to the Cartesian coordinates $\left(x,y,z\right)$ by
\[
r=\sqrt{x^{2}+y^{2}+z^{2}},\quad\tan\theta=\nicefrac{y}{x},\quad\cos\varphi=\nicefrac{z}{r},
\]
with $0\leq\theta\leq2\pi$ and $0\leq\varphi\leq\pi$. Let $T$ be
the tension in direction $\theta=0$ acting as a load. Then, the solution
for a cavity with radius $R$ is given in spherical displacement components
as follows \cite{Goodier_1933a}:
\begin{eqnarray}
u_{r} & = & u_{r}^{0}-\dfrac{A}{r^{2}}-\dfrac{3B}{r^{4}}+\left(\dfrac{5-4\nu}{1-2\nu}\cdot\dfrac{C}{r^{2}}-\dfrac{9B}{r^{4}}\right)\cos\left(2\theta\right),\label{eq:SphereHoleDisplR}\\
u_{\theta} & = & u_{\theta}^{0}-\left(\dfrac{2C}{r^{2}}+\dfrac{6B}{r^{4}}\right)\sin\left(2\theta\right),\label{eq:SphereHoleDisplTheta}\\
u_{\varphi} & = & 0,\nonumber 
\end{eqnarray}
with
\begin{eqnarray}
u_{r}^{0} & = & \dfrac{Tr}{2E}\left[\left(1-\nu\right)+\left(1+\nu\right)\cos\left(2\theta\right)\right],\label{eq:SphereHoleDisplR0}\\
u_{\theta}^{0} & = & -\dfrac{Tr}{2E}\left(1+\nu\right)\sin\left(2\theta\right),\label{eq:SphereHoleDisplTheta0}
\end{eqnarray}
and the coefficients
\begin{eqnarray*}
A=-\dfrac{R^{3}T}{8\mu}\cdot\dfrac{13-10\nu}{7-5\nu}, & \quad B=\dfrac{R^{5}T}{8\mu}\cdot\dfrac{1}{7-5\nu}, & \quad C=\dfrac{R^{3}T}{8\mu}\cdot\dfrac{5\left(1-2\nu\right)}{7-5\nu}.
\end{eqnarray*}
This is easily converted to Cartesian coordinates based on the transformation
matrix $\mat T$,
\[
\left[\begin{array}{c}
u_{x}\\
u_{y}\\
u_{z}
\end{array}\right]=\mat T\cdot\left[\begin{array}{c}
u_{r}\\
u_{\theta}\\
u_{\varphi}
\end{array}\right]\qquad\mathrm{with}\;\mat T=\left[\begin{array}{ccc}
\sin\theta\cos\varphi & \cos\theta\cos\varphi & -\sin\varphi\\
\sin\theta\sin\varphi & \cos\theta\sin\varphi & \cos\varphi\\
\cos\theta & -\sin\theta & 0
\end{array}\right].
\]

\subsection{Exact solution for the cube with spherical inclusion\label{XX_ExactSolCubeWithSphereInclusion}}

The exact solution in spherical displacements for the spherical inclusion
problem of Section \ref{XX_CubeWithSphereInclusion} is given as \cite{Goodier_1933a}
\begin{eqnarray*}
u_{r} & = & \begin{cases}
Hr+Fr+3Fr\cdot\cos\left(2\theta\right) & \mathrm{for}\; r\leq R\\
u_{r}^{0}-\dfrac{A}{r^{2}}-\dfrac{3B}{r^{4}}+\left(\dfrac{5-4\nu_{1}}{1-2\nu_{1}}\cdot\dfrac{C}{r^{2}}-\dfrac{9B}{r^{4}}\right)\cos\left(2\theta\right) & \mathrm{for}\; r>R
\end{cases}\\
u_{\theta} & = & \begin{cases}
-3Fr\cdot\sin\left(2\theta\right) & \mathrm{for}\; r\leq R\\
u_{\theta}^{0}-\left(\dfrac{2C}{r^{2}}+\dfrac{6B}{r^{4}}\right)\sin\left(2\theta\right) & \mathrm{for}\; r>R
\end{cases}\\
u_{\varphi} & = & 0
\end{eqnarray*}
with
\begin{eqnarray}
u_{r}^{0} & = & \dfrac{Tr}{2E_{1}}\left[\left(1-\nu_{1}\right)+\left(1+\nu_{1}\right)\cos\left(2\theta\right)\right],\label{eq:SphereBiMatDisplR0}\\
u_{\theta}^{0} & = & -\dfrac{Tr}{2E_{1}}\left(1+\nu_{1}\right)\sin\left(2\theta\right),\label{eq:SphereBiMatDisplTheta0}
\end{eqnarray}
and the coefficients
\begin{eqnarray*}
A & = & R^{3}\cdot\Bigg[-\dfrac{T}{8\mu_{1}}\cdot\dfrac{\mu_{1}-\mu_{2}}{\left(7-5\nu_{1}\right)\mu_{1}+\left(8-10\nu_{1}\right)\mu_{2}}\cdot\\
 &  & \dfrac{(1-2\nu_{2})\cdot(6-5\nu_{1})2\mu_{1}+(3+19\nu_{2}-20\nu_{1}\nu_{2})\mu_{2}}{\left(1-2\nu_{2}\right)2\mu_{1}+\left(1+\nu_{2}\right)\mu_{2}}\\
 &  & +\dfrac{T}{4\mu_{1}}\cdot\dfrac{\left[\left(1-\nu_{1}\right)\left(1+\nu_{2}\right)/\left(1+\nu_{1}\right)-\nu_{2}\right]\mu_{2}-(1-2\nu_{2})\mu_{1}}{\left(1-2\nu_{2}\right)2\mu_{1}+\left(1+\nu_{2}\right)\mu_{2}}\Bigg],\\
B & = & \dfrac{R^{5}T}{8\mu_{1}}\cdot\dfrac{\mu_{1}-\mu_{2}}{\left(7-5\nu_{1}\right)\mu_{1}+\left(8-10\nu_{1}\right)\mu_{2}},\\
C & = & \dfrac{R^{3}T}{8\mu_{1}}\cdot\dfrac{5\left(1-2\nu_{1}\right)\left(\mu_{1}-\mu_{2}\right)}{\left(7-5\nu_{1}\right)\mu_{1}+\left(8-10\nu_{1}\right)\mu_{2}},\\
F & = & \dfrac{5T}{4}\cdot\dfrac{1-\nu_{1}}{\left(7-5\nu_{1}\right)\mu_{1}+\left(8-10\nu_{1}\right)\mu_{2}},\\
H & = & \dfrac{T\left(1-\nu_{1}\right)}{2\left(1+\nu_{1}\right)}\cdot\dfrac{1-2\nu_{2}}{\left(2-4\nu_{2}\right)\mu_{1}+\left(1+\nu_{2}\right)\mu_{2}}.
\end{eqnarray*}

\bibliographystyle{schanz}
\addcontentsline{toc}{section}{\refname}\bibliography{FriesRefs}
 
\end{document}